\documentclass[12pt]{article}
\usepackage{amssymb,amsmath,latexsym}
\usepackage{color}
\allowdisplaybreaks

\begin{document}

\begin{doublespace}

\def\1{{\bf 1}}
\def\ind{{\bf 1}}
\def\nn{\nonumber}
\def\bee{\begin{equation}}
\def\eee{\end{equation}}
\def\sA {{\cal A}} \def\sB {{\cal B}} \def\sC {{\cal C}}
\def\sD {{\cal D}} \def\sE {{\cal E}} \def\sF {{\cal F}}
\def\sG {{\cal G}} \def\sH {{\cal H}} \def\sI {{\cal I}}
\def\sJ {{\cal J}} \def\sK {{\cal K}} \def\sL {{\cal L}}
\def\sM {{\cal M}} \def\sN {{\cal N}} \def\sO {{\cal O}}
\def\sP {{\cal P}} \def\sQ {{\cal Q}} \def\sR {{\cal R}}
\def\sS {{\cal S}} \def\sT {{\cal T}} \def\sU {{\cal U}}
\def\sV {{\cal V}} \def\sW {{\cal W}} \def\sX {{\cal X}}
\def\sY {{\cal Y}} \def\sZ {{\cal Z}}

\def\bA {{\mathbb A}} \def\bB {{\mathbb B}} \def\bC {{\mathbb C}}
\def\bD {{\mathbb D}} \def\bE {{\mathbb E}} \def\bF {{\mathbb F}}
\def\bG {{\mathbb G}} \def\bH {{\mathbb H}} \def\bI {{\mathbb I}}
\def\bJ {{\mathbb J}} \def\bK {{\mathbb K}} \def\bL {{\mathbb L}}
\def\bM {{\mathbb M}} \def\bN {{\mathbb N}} \def\bO {{\mathbb O}}
\def\bP {{\mathbb P}} \def\bQ {{\mathbb Q}} \def\bR {{\mathbb R}}
\def\bS {{\mathbb S}} \def\bT {{\mathbb T}} \def\bU {{\mathbb U}}
\def\bV {{\mathbb V}} \def\bW {{\mathbb W}} \def\bX {{\mathbb X}}
\def\bY {{\mathbb Y}} \def\bZ {{\mathbb Z}}
\def\R {{\mathbb R}} \def\RR {{\mathbb R}} \def\H {{\mathbb H}}
\def\n{{\bf n}} \def\Z {{\mathbb Z}}

\newcommand{\expr}[1]{\left( #1 \right)}
\newcommand{\cl}[1]{\overline{#1}}
\newtheorem{thm}{Theorem}[section]
\newtheorem{lemma}[thm]{Lemma}
\newtheorem{defn}[thm]{Definition}
\newtheorem{prop}[thm]{Proposition}
\newtheorem{corollary}[thm]{Corollary}
\newtheorem{remark}[thm]{Remark}
\newtheorem{example}[thm]{Example}
\numberwithin{equation}{section}
\def\ee{\varepsilon}
\def\qed{{\hfill $\Box$ \bigskip}}
\def\NN{{\mathcal N}}
\def\AA{{\mathcal A}}
\def\MM{{\mathcal M}}
\def\BB{{\mathcal B}}
\def\CC{{\mathcal C}}
\def\LL{{\mathcal L}}
\def\DD{{\mathcal D}}
\def\FF{{\mathcal F}}
\def\EE{{\mathcal E}}
\def\QQ{{\mathcal Q}}
\def\SS{{\mathcal S}}
\def\RR{{\mathbb R}}
\def\R{{\mathbb R}}
\def\L{{\bf L}}
\def\K{{\bf K}}
\def\S{{\bf S}}
\def\A{{\bf A}}
\def\E{{\mathbb E}}
\def\F{{\bf F}}
\def\P{{\mathbb P}}
\def\N{{\mathbb N}}
\def\eps{\varepsilon}
\def\wh{\widehat}
\def\wt{\widetilde}
\def\pf{\noindent{\bf Proof.} }
\def\pff{\noindent{\bf Proof} }
\def\cp{\mathrm{Cap}}
\def\UU{\mathcal U}

\title{\Large \bf Potential theory of subordinate killed Brownian motion}

\author{{\bf Panki Kim}\thanks{This work was supported by National Research Foundation of Korea(NRF) grant funded by the Korea government(MSIP) (No. 2016R1E1A1A01941893).
}
\quad {\bf Renming Song\thanks{Research supported in part by a grant from
the Simons Foundation (\#429343, Renming Song).}} \quad and
\quad {\bf Zoran Vondra\v{c}ek}
\thanks{Research supported in part by the Croatian Science Foundation under the project 3526.}
}

\date{}

\maketitle

\begin{abstract}
Let $W^D$ be a killed Brownian motion in a domain $D\subset \R^d$ and $S$ an independent
subordinator with Laplace exponent $\phi$. The process $Y^D$ defined by $Y^D_t=W^D_{S_t}$ is called a subordinate killed Brownian motion.  
It is a Hunt process with infinitesimal generator $-\phi(-\Delta|_D)$, where
$\Delta|_D$ is the Dirichlet Laplacian. 
In this paper we study the potential theory of $Y^D$ under a weak scaling 
condition on the derivative of  $\phi$.
We first show that non-negative harmonic functions of $Y^D$ satisfy the scale invariant Harnack inequality.  
Subsequently we prove two types of scale invariant boundary Harnack principles with explicit decay rates for non-negative harmonic functions of $Y^D$. The first boundary Harnack principle  deals with a $C^{1,1}$ domain $D$ and non-negative functions which are harmonic near the boundary of $D$, while the second one is 
for a more general domain $D$
and non-negative functions which are harmonic near the boundary of an interior open subset of $D$. The obtained decay rates are not the same, reflecting different boundary and interior behaviors of $Y^D$.
\end{abstract}



\section{Introduction}\label{s:intro}
Let $W=(W_t, \P_x)$ be a Brownian motion in $\R^d$, $d\ge 1$, and let $S=(S_t)_{t\ge 0}$ be an independent subordinator with Laplace exponent $\phi$. The process $X=(X_t, \P_x)$ defined by $X_t=W_{S_t}$, $t\ge 0$,
is called a subordinate Brownian motion. It is an isotropic L\'evy process with characteristic exponent $\Psi(\xi)=\phi(|\xi|^2)$. 

In recent years, isotropic, and more generally, symmetric, L\'evy processes have
been intensively studied and many important results have been obtained. 
In particular, under certain weak scaling conditions on the characteristic exponent 
$\Psi$ (or the Laplace exponent $\phi$), it was shown that non-negative harmonic
functions with respect to these L\'evy processes satisfy the 
scale invariant Harnack inequality (HI) and the scale invariant boundary Harnack 
principle (BHP), 
see e.g.~\cite{BKK, Gr, KM, KM14, KSV09, KSV12, KSV12b, KSV14, KSV17}.

If $D$ is an open subset of $\R^d$, we can kill the process $X$ upon exiting $D$ and obtain a process $X^D$ known as a killed
subordinate Brownian motion. Functions that are harmonic in an open subset of $D$ with respect to $X^D$ are defined only on $D$, but by extending them to be identically zero on $\R^d\setminus D$, 
the HI and BHP follow directly from
those for $X$.
In fact, the HI (BHP, respectively)  for $X^D$ is
a special case of the HI (BHP, respectively) for $X$ in $D$.

By reversing the order of subordination and killing, one obtains a process different from $X^D$. Assume from now on that
$D$ is a  domain (i.e.,~connected open set) in $\R^d$, and let $W^D=(W^D_t, \P_x)$ be the Brownian motion $W$ killed upon
exiting $D$. The process $Y^D=(Y^D_t,\P_x)$ defined by $Y^D_t=W^D_{S_t}$, $t\ge 0$, is called a subordinate killed Brownian motion.
It is a symmetric Hunt process (see, for instance, \cite[Propositions 2.2 and 2.3]{SV08}) with 
infinitesimal generator  $-\phi(-\Delta_{|D})$, 
where $\Delta_{|D}$ is the Dirichlet Laplacian. 
This process is very natural and useful. For example, it was used  in
\cite{CS05} as a tool
to obtain two-sided estimates of the eigenvalues of the generator of $X^D$.
Despite its usefulness, the potential theory of subordinate killed Brownian motions
has been studied only sporadically, 
see e.g.~\cite{GRSS, SV03, GPRSSV} for stable subordinators, and \cite{SV06, SV09} for more general subordinators. In particular, \cite{SV09} contains versions of HI and BHP (with respect to the subordinate killed Brownian motion 
in a bounded Lipschitz domain $D$) 
which are very weak in the sense that the results are proved 
only for  non-negative functions which are harmonic in all of $D$. Those results are 
easy consequences of the fact that there is a one-to-one correspondence 
between non-negative harmonic functions (in all of $D$) with respect to $W^D$ and those with respect to $Y^D$. 
Additionally, some aspects of potential theory of subordinate killed Brownian motions 
in unbounded domains were recently studied in \cite{KSV16}.

In the partial differential equations literature, 
the operator $-(-\Delta|_D)^{\alpha/2}$, $\alpha\in (0, 2)$,
which
is the generator of the subordinate killed Brownian motion via an $\alpha/2$-stable
subordinator, also goes under the name of spectral fractional Laplacian, 
see \cite{BSV} and the references therein.  
This operator has been of interest to quite a few people in the PDE circle. 
For instance, a
version of  HI was also shown in \cite{SZ}.

The main goal of this paper is to show that the scale invariant BHP 
for non-negative functions harmonic (with respect to $Y^D$) near the boundary of $D$ is valid for a large class of subordinate killed Brownian motions  
$Y^D$ when $D$ is a bounded $C^{1,1}$ domain, or a $C^{1,1}$ domain 
with compact complement or a domain consisting of all the points above the 
graph of a bounded globally $C^{1,1}$ function. 
We also prove the scale invariant BHP for 
non-negative functions
harmonic (with respect to $Y^D$) near the boundary of a $C^{1,1}$ open set strictly contained in $D$. In the latter case $D$ need not have smooth boundary, but still has to satisfy certain geometric conditions.

We start by listing our assumptions on the underlying subordinator, 
or more precisely on the corresponding Laplace exponent.
Recall that the Laplace exponent of a subordinator is a Bernstein function vanishing at the origin, i.e.,~it has the representation
\begin{equation}\label{e:LKdecomp}
\phi(\lambda)=b\lambda +\int_{(0,\infty)}(1-e^{-\lambda t})\, \mu(dt)\, ,
\end{equation}
with $b\ge 0$ and $\mu$ a measure on $(0,\infty)$ satisfying
$\int_{(0,\infty)}(1\wedge t)\, \mu(dt)<\infty$, which is called the L\'evy measure of $S$.
A Bernstein function $\phi$ is called a complete Bernstein function if its L\'evy measure has a completely monotone density.
A Bernstein function $\phi$ is called a special Bernstein function if the
function $\lambda\mapsto \lambda/\phi(\lambda)$ is also a Bernstein function.
 It is well known that any complete Bernstein function is a special Bernstein function. 
For this and other properties of special and complete Bernstein functions, see \cite{SSV}.
The potential measure of the subordinator $S$ is defined by $U(A)=\int_0^\infty\P(S_t\in A)\,dt$.

In this the paper we will always assume the following three conditions:
\begin{itemize}
	\item[{\bf (A1)}]
	The potential measure $U$ of $S$ has a decreasing density $u$.
	\item[{\bf(A2)}]
	The L\' evy measure of $S$ is infinite and has a decreasing density $\mu$ that satisfies 
	\begin{align}
\label{e:CBmu0}
\mu(r) \le c \mu(r+1),  \quad r>1.
\end{align}

	\item[{\bf(A3)}] 
	There exist constants $\sigma>0$ and
$\delta \in (0, 1]$
 such that
\begin{equation*}
  \frac{\phi'(\lambda t)}{\phi'(\lambda)}\leq\sigma\, t^{-\delta}\ \text{ for all }\ t\geq 1\ \text{ and }\
 \lambda\geq1\, .
\end{equation*}
\end{itemize}

For some of the results in dimension $d=2$ we will need the following additional assumption: 
\begin{itemize}
	\item[{\bf(A4)}]
There exist $\sigma'>0$ and  $\delta'\in (0,2\delta)$
 such that
\begin{equation}\label{e:new231}
  \frac{\phi'(\lambda t)}{\phi'(\lambda)}\geq
\sigma'\,t^{-\delta'}\ \text{ for all
}\ t\geq 1\ \text{ and }\ \lambda\geq1\,.
\end{equation}
\end{itemize}

Assumptions {\bf(A1)}-{\bf(A4)} were introduced in \cite{KM} and were also used in 
\cite{KM14, KSV16}.
Assumption {\bf(A1)} implies that  $\phi$ is  a special Bernstein function, cf.~\cite[Theorem 11.3]{SSV}.
It follows from \cite[Lemma 2.1]{KSV12b} that \eqref{e:CBmu0} holds if $\phi$ is  a complete Bernstein function.
It immediately follows from {\bf(A3)} that 
$b=0$ in \eqref{e:LKdecomp}.
It is easy to check that
if $\phi$ is a complete Bernstein function satisfying the following weak scaling condition: There exist $a_1, a_2>0$ and
$\delta_1, \delta_2\in (0,1)$ satisfying
$$
a_1 \lambda^{\delta_1}\phi(t)\le \phi(\lambda t)\le a_2 \lambda^{\delta_2}\phi(t)\, ,\qquad \lambda \ge 1, t\ge 1\, ,
$$
then {\bf (A1)}--{\bf (A4)} (as well as  {\bf (A7)} below) are satisfied 
(see \cite[p.~4386]{KM14}).
One of the reasons for adopting the more general setup above is to cover the case
of geometric stable (see examples (5) and (6) below) and iterated geometric stable subordinators.

The following assumption is a necessary and sufficient condition for the transience of $X$, cf.~\cite[(5.46)]{SV09}. It is always satisfied when $d\ge 3$, hence imposes an additional restriction only in case $d\le 2$.

\begin{itemize}
\item[{\bf(A5)}] We assume that
\begin{equation*}
\int_{0}^{1}\frac{\lambda^{\frac{d}2-1}}{\phi(\lambda)}d \lambda<\infty.
\end{equation*}
\end{itemize}

In Sections \ref{s:jkgfe}--\ref{s:bhi}, 
when $d=2$ we assume that $D$ is bounded, while when $d \ge 3$ we allow $D$ to be unbounded. In case $D$ is unbounded, we assume the following extra condition:
\begin{itemize}
\item[{\bf(A6)}] 
There exist $\beta, \sigma_1 >0$ such that
\begin{align}\label{e:mas}
\frac{\mu(\lambda t)}{\mu(\lambda)} \ge  \sigma_1 t^{-\beta }
\quad \text{ for all } t\ge 1  \text{ and } \lambda>0\, .
\end{align}
\end{itemize}
Condition  \eqref{e:mas}  is used to control the large jumps of $Y^D$. The reason that we
assume $d\ge 3$ is that sharp two-sided Dirichlet heat kernel estimates
for $C^{1,1}$ domains with compact complement, which are used in Theorem \ref{t:Jfe:nn}, are only available for $d\ge3$, see \cite{Zh2}.

In the last section of the paper we will need an additional assumption used in \cite{KM14}. This assumption will 
not be used explicitly in the current paper, 
but will be used implicitly when citing results from \cite{KM14}.
After Lemma \ref{l:eqofforms} in Section \ref{s:bhi2} below, 
 if the constant $\delta$ in \textbf{(A3)} satisfies $0<\delta \le 1/2$, we will assume the following extra assumption:
\begin{itemize}
\item[{\bf(A7)}] 
There
exist $\sigma_2>0$ and $\gamma \in [\delta, 1)$ such that
\begin{align}\label{e:A5-fromKM14}
\frac{\phi(\lambda t)}{\phi(\lambda)}\ge \sigma_2 t^{1-\gamma} \quad \text{ for all } t\ge 1  \text{ and } \lambda \ge 1\, .
\end{align}
\end{itemize}

The following well-known examples of subordinators satisfy 
 {\bf (A1)}--{\bf (A3)}. 
Note that in all of the examples the Bernstein function $\phi$ is complete, cf.~\cite{SSV}, so {\bf (A1)}--{\bf (A2)} are satisfied. The parameter $\delta$ from {\bf (A3)} is written for future reference. 
\begin{description}
\item{(1)} Stable subordinator: $\phi(\lambda)=\lambda^{\alpha}$, $0<\alpha<1$, with  $\delta=1-\alpha$.
\item{(2)} Sum of two stable subordinators: $\phi(\lambda)=\lambda^\beta + \lambda^\alpha$, $0<\beta<\alpha<1$, with $\delta=1-\alpha$.
\item{(3)} Stable with logarithmic correction: $\phi(\lambda)=\lambda^\alpha(\log(1+\lambda))^\beta$, $0<\alpha<1$,
$0<\beta < 1-\alpha$, with $\delta=1-\alpha-\epsilon$ for every $\epsilon >0$.
\item{(4)}  Stable with logarithmic correction: $\phi(\lambda)=\lambda^\alpha(\log(1+\lambda))^{-\beta}$, $0<\alpha<1$,
$0<\beta < \alpha$, with $\delta=1-\alpha$.
\item{(5)} Geometric stable subordinator: $\phi(\lambda)=\log(1+\lambda^{\alpha})$, $0<\alpha< 1$, with $\delta=1$.
\item{(6)} Gamma subordinator: $\phi(\lambda)=\log(1+\lambda)$, with $\delta=1$. 
\item{(7)} Relativistic stable subordinator: $\phi(\lambda)=(\lambda+m^{1/\alpha})^{\alpha}-m$, $0<\alpha<1$ and $m>0$, with $\delta=1-\alpha$. 
\end{description} 
It is easy to see that {\bf (A4)} and {\bf (A7)} also hold for (1)--(4) and (7). (Note, since we will use {\bf (A7)} only if the constant $\delta$ in \textbf{(A3)} satisfies $0<\delta \le 1/2$, {\bf (A7)} will not be applied to examples (5) and (6).)
Condition  {\bf (A5)} is true when $\alpha<d/2$ in examples (1), (2) and (5), $\alpha+\beta<d/2$ in (3), $\alpha-\beta<d/2$ in (4), and $d>2$ in (6) and (7).  By using \cite[Sections 5.2.2 and 5.2.3]{SV09} one checks that \eqref{e:mas} is satisfied for examples (1)-(5), but not  for examples (6) and (7) since 
the corresponding L\'evy density has exponential decay.
Thus for  examples (6) and (7) we can only cover the case when $D$ is bounded. 

By using the tables at the end of \cite{SSV} one can come up with a lot of explicit examples of complete Bernstein functions such that conditions {\bf (A1)}--{\bf (A7)} are true.

Let $D\subset \R^d$ be an open set and let $Q\in \partial D$. We say that $D$ is $C^{1,1}$ near $Q$ if there exist 
a localization radius $R>0$, 
a $C^{1,1}$-function 
$\varphi=\varphi_Q:\R^{d-1}\to \R$ satisfying $\varphi(0)=0$, $\nabla \varphi(0)=(0,\dots, 0)$, $\| \nabla \varphi \|_{\infty}\le \Lambda$, $|\nabla \varphi(z)-\nabla \varphi(w)|\le \Lambda |z-w|$, and an orthonormal coordinate system $CS_Q$ with its origin at $Q$ such that
$$
B(Q,R)\cap D=\{y=(\wt{y},y_d)\in B(0,R) \textrm{ in } CS_Q:\, y_d>\varphi(\wt{y})\}\, ,
$$
where $\wt{y}:= (y_1, \dots, y_{d-1})$.
The pair $(R,\Lambda)$ will be called the $C^{1, 1}$ characteristics of $D$ at $Q$.
Recall that an open set $D\subset \R^d$ is said to be a (uniform) $C^{1,1}$ open set with characteristics $(R,\Lambda)$ if it is $C^{1,1}$ with characteristics $(R,\Lambda)$ near every boundary point $Q\in \partial D$.
For future reference, given $\Lambda>0$, we define
\begin{align}\label{e:kappa1}
\kappa=\kappa(\Lambda):=(1+(1+\Lambda)^2)^{-1/2}.
\end{align}

For any Borel set $B\subset D$, let $\tau_B=\tau_B^{Y^D}=\inf\{t>0:\, Y^D_t\notin B\}$ be the exit time of $Y^D$ from $B$. 
\begin{defn}\label{D:1.1}  \rm 
A non-negative Borel function defined on
$D$ is said to be  \emph{harmonic}
in an open set $V\subset D$ with respect to
$Y^D$ if for every bounded open set $U\subset\overline{U}\subset V$,
\begin{equation}\label{e:har}
f(x)= \E_x \left[ f(Y^D_{\tau_{U}})\right] \qquad
\hbox{for all } x\in U.
\end{equation}
A non-negative Borel function $f$ defined on $D$ is said to be \emph{regular harmonic} in an open set $V\subset D$ if 
$$
f(x)= \E_x \left[ f(Y^D_{\tau_{V}})\right] \qquad \hbox{for all } x\in V.
$$
\end{defn}

For any open set $U\subset \R^d$ and $x\in \R^d$, 
we use $\delta_U(x)$ to denote the distance between $x$ and the boundary $\partial U$.

The main result of this paper is the following scale invariant boundary Harnack inequality.
\begin{thm}\label{t:main}
Suppose that {\bf (A1)}--{\bf (A3)} hold. If the constant $\delta$ in {\bf (A3)} satisfies $0<\delta \le 1/2$, then we also assume 
that  {\bf (A7)} holds.
 If $d=2$, we assume that 
 {\bf (A4)} and {\bf (A5)} hold.
Let $D$ be a bounded $C^{1,1}$ domain, or a $C^{1,1}$ domain with compact complement or a domain consisting of all the points above the graph of a 
bounded globally $C^{1,1}$ function. 
In case of unbounded $D$ we consider  $d\ge 3$ only  and  assume that {\bf (A6)} holds.
Let $(R,\Lambda)$ be the $C^{1, 1}$ characteristics of $D$.
There exists a constant $C=C(d, \Lambda, R, \phi)>0$ 
such that for any $r \in (0, R]$, $Q\in \partial D$, and any non-negative function $f$ in $D$ which is harmonic  in $D \cap B(Q, r)$ with respect to $Y^D$ and vanishes continuously on $ \partial D \cap B(Q, r)$, we have
\begin{equation}\label{e:bhp_m}
\frac{f(x)}{\delta_D(x)}\,\le C\,\frac{f(y)}{\delta_D(y)} \qquad
\hbox{for all } x, y\in  D \cap B(Q, r/2).
\end{equation}
\end{thm}

In particular, we see from the theorem above that if a non-negative function which is harmonic with respect to $Y^D$  vanishes  near the boundary, then its rate of decay is proportional to the distance to the boundary (regardless of the particular subordinator as long as 
the conditions in Theorem \ref{t:main} hold). 
This shows that near the boundary of $D$, $Y^D$ behaves like the killed Brownian motion $W^D$.

We remark that Theorem \ref{t:main} is new even in the case of a stable subordinator. 
Recently, a BHP for a general discontinuous Feller process in metric measure space
has  been proved in \cite{BKK, KSV16} under some comparability assumptions on
the jumping kernel (see \cite[Assumption C]{BKK} and \cite[{\bf C1}$(z_0, r_0)$]{KSV16}) and a Urysohn-type property of the domain of the generator of the process (see \cite[Assumption D]{BKK} and \cite[{\bf B1}$(z_0, r_0)$]{KSV16}). 
However,  neither \cite[Theorem 3.4]{BKK} nor 
\cite[Corollary 4.2]{KSV17} can be applied to subordinate 
killed Brownian motions because \cite[Assumption C]{BKK} and
\cite[{\bf C1}$(z_0, r_0)$]{KSV17} do not hold. 
Moreover, one can observe from \eqref{e:BHP-yes} below that
the approximate factorization of
non-negative functions harmonic for the subordinate killed 
Brownian motion should be different
from those  in  \cite[(3.3)]{BKK} and
\cite[Theorem 4.1]{KSV17}. 

The proof of Theorem \ref{t:main} is probabilistic and is based on the 
``box method'', which was originally proposed in  \cite{BB}. 
This method has been extended and used for many discontinuous processes (for 
example, see \cite{BoB, BBC, KS, CKSV, KSV13}).
It  is not surprising that our method is similar to that 
for the censored stable process in \cite{BBC} -- both the censored process and the subordinate killed Brownian motion are intrinsically defined as processes living in a domain $D$ (hence no information outside $D$ is available). As a consequence, some of the domain monotonicity properties are not valid. 
Despite the fact that the general road map has been traced before, the proof of Theorem \ref{t:main} is still technically quite challenging. Some of the necessary  technical results such as Green function and jumping kernel estimates for $Y^D$ were developed in \cite{KSV12}. These estimates are improved and complemented in Section \ref{s:jkgfe}. 
In order to prove Theorem \ref{t:main}, we 
establish  a Carleson estimate in Section \ref{s:ce}
by following the ideas from \cite{BBC, G, CKSV, KSV13}.
Another key step to  prove Theorem \ref{t:main}  is  obtaining 
the correct explicit estimate on
the exit distribution (see \eqref{e:L:2} and \eqref{e:D2_43}). 
Unlike previous  papers, we do not use testing functions. 
Instead of using Dynkin's formula and applying to testing functions, we utilize  relations 
among the  killed subordinate Brownian motion, the subordinate killed Brownian motion and its killed processes,  
and  estimates of the Green function and jumping kernels  for these processes.

Theorem \ref{t:main} is a result about the boundary behavior of non-negative functions on $D$ that are harmonic (with respect to $Y^D$) 
near a portion of the boundary of $D$ and
vanish continuously at that part of the boundary. The corresponding decay rate is a reflection of the fact that, near the boundary, the subordinate killed Brownian motion $Y^D$ behaves like Brownian motion.  One can also study the decay rate of non-negative functions in $D$ which are (regular) harmonic near 
a portion of the boundary, strictly contained in $D$, of an open set $E\subset D$
and vanish in an appropriate sense. 
Intuitively, since the behavior of $Y^D$ in the interior of $D$ is similar to that of the killed subordinate Brownian motion $X^D$, one would expect the decay rate to be the same as the decay rate of functions harmonic with respect to $X^D$. This is confirmed in the second main result of this paper below. Furthermore, since such a result concerns only 
the interior of $D$, the smoothness of the 
boundary of $D$ is no longer necessary. 
Still, some geometric conditions for $D$ are needed. These conditions 
are related to the heat kernel $p^D(t,x,y)$ of the killed Brownian motion $W^D$
and its tail function $t\mapsto \P_x(t<\tau_D^W)$.

We will say that a decreasing function $f:(0,\infty)\to (0,\infty)$ satisfies the doubling property 
if, for every $T>0$, there exists a constant $c>0$ such that $f(t)\le cf(2t)$ for all $t\in (0, T]$. 
\begin{itemize} 
	\item[{\bf (B1)}] The function $t\mapsto \P_x(t<\tau_D^W)$ satisfies 
	the doubling property (with a doubling constant independent of $x\in D$).
	\item[{\bf (B2)}] There exist constants $c\ge 1$ and $M\ge 1$ such that for
	all  $t\le1$ and $x, y\in D$,
\begin{align}
\label{e:fac}
&c^{-1}\,\P_x(t<\tau^W_D)\,\P_y(t<\tau^W_D)\,t^{-d/2}e^{-\frac{M|x-y|^2}{t}}\nonumber\\
&\le
p^D(t, x, y)
\le c\,\P_x(t<\tau^W_D)\,\P_y(t<\tau^W_D)\,t^{-d/2}e^{-\frac{|x-y|^2}{Mt}}.
\end{align}
\end{itemize}

If $D$ is either a bounded Lipschitz domain or an unbounded domain
consisting of all the points above the graph of a globally Lipschitz function,
then {\bf (B1)} and {\bf (B2)} are satisfied, cf. \cite[(0.36) and (0.25)]{Varo}. 
It is  also easy to show that a $C^{1,1}$ domain with compact complement also 
satisfies {\bf (B1)} and {\bf (B2)}.
\begin{thm}\label{t:bhp-2}
Suppose that  {\bf (A1)}--{\bf (A3)} hold.
If the constant $\delta$ in {\bf (A3)} satisfies $0<\delta \le 1/2$, then we also assume 
that  {\bf (A7)} holds.
If $d=2$, we assume that 
 {\bf (A4)} and  {\bf (A5)} also hold.
Let $D\subset \R^d$ be a domain satisfying {\bf (B1)} and {\bf (B2)}. 
There exists a constant $b=b(\phi,d)>2$ such that, for 
every open set $E\subset D$ and every $Q\in \partial E\cap D$  such that $E$ is $C^{1,1}$ near $Q$ with characteristics
 $(R, \Lambda)$, the following holds: 
There exists a constant $C=C(\delta_D(Q)\wedge R,\Lambda,\phi,d)>0$ 
such that for every 
$r\le (\delta_D(Q)\wedge R)/(b+2)$ 
and every non-negative function $f$ on $D$ which is regular harmonic in $E\cap B(Q,r)$ with respect to $Y^D$ and vanishes on $E^c\cap B(Q,r)$, we have 
$$
\frac{f(x)}{\phi(\delta_E(x)^{-2})^{-1/2}}\le C \frac{f(y)}{\phi(\delta_E(y)^{-2})^{-1/2}}\, , \qquad x,y\in E\cap B(Q,2^{-6}\kappa^4 r)\, ,
$$
where $\kappa$ is the constant defined in \eqref{e:kappa1}.
\end{thm}

Again, Theorem \ref{t:bhp-2} is new even in the case of a stable subordinator. The method of proof of  Theorem \ref{t:bhp-2} is quite different from that of Theorem \ref{t:main}. It relies on a comparison of the Green functions 
of subprocesses of $Y^D$ and $X$ for small interior subsets of $D$, 
and on some already available potential-theoretic results for $X$ obtained in \cite{KM14}. To be more precise, let $U\subset D$ be a 
$C^{1,1}$ open set 
with diameter comparable to its distance to $\partial D$, and let $Y^{D,U}$ and  $X^U$ denote processes $Y^D$ and $X$ killed upon exiting $U$
respectively. We show that the Dirichlet form of $Y^{D,U}$ is equal to the Dirichlet form of a non-local Feynman-Kac-type transform of $X^U$, and that the corresponding conditional  gauge is bounded from below, cf.~Lemmas \ref{l:eqofforms} and \ref{l:bounded-gauge}. 
This immediately implies 
comparability of the Green functions. The proof of Theorem \ref{t:bhp-2} now mainly 
uses the corresponding result for the process $X$ and properties of the jumping kernel of $Y^D$.

Finally, one of the ingredients in the proof of Theorem \ref{t:main} is the scale invariant HI.
We will show in Theorem \ref{uhp} that when  {\bf (A1)}--{\bf (A3)} and {\bf (B1)}--{\bf (B2)} hold,
there exists a constant $C>0$ such that for any $r\in (0,1]$ and $B(x_0, r) \subset D$ and any function $f$ which is non-negative in $D$ and harmonic in $B(x_0, r)$ with respect to $Y^D$, we have
$$
f(x)\le C f(y) \qquad \mbox{ for all } x, y\in B(x_0, r/2).
$$
The proof of the HI is modeled after the powerful method developed in \cite{KM}.

Organization of the paper: In the next section we 
collect several results concerning subordinators satisfying 
{\bf (A1)}--{\bf (A3)}, the subordinate killed Brownian motion $Y^D$ and its relation
with the killed subordinate Brownian motion $X^D$.
In Section \ref{s:hi} we prove 
the scale invariant HI for $Y^D$.
As preparation for the subsequent sections, in Section \ref{s:jkgfe} 
we  give sharp  two-sided estimates
on the jumping kernel and Green function of $Y^D$.
The Carleson estimate, an important ingredient in proving the BHP, is obtained in Section \ref{s:ce}.
We continue in Section \ref{s:bhi} with the proof of the
BHP in $C^{1,1}$ domains with explicit decay rate.
The proof of Theorem \ref{t:bhp-2} is given in Section \ref{s:bhi2}. This last section can be read independently of Sections \ref{s:jkgfe}--\ref{s:bhi} and uses only Lemmas \ref{l:ilb4BMheatkernel}, \ref{l:ilb4jk} and 
Proposition \ref{p:estimate-of-J-away}.

Notation: We will use the following conventions in this paper.
Capital letters $C, C_i, i=1,2,  \dots$ will denote the constants
in the statements of results and the labeling of these constants starts anew
in each result. Lower case letters 
$c, c_i, i=1,2,  \dots$ are used to denote the constants in the proofs
and the labeling of these constants starts anew in each proof.
$c_i=c_i(a,b,c,\ldots)$, $i=0,1,2,  \dots$ denote  constants depending on $a, b, c, \ldots$.
Dependence of constants on the constants in  
{\bf (A1)}--{\bf (A7)} and {\bf (B1)}--{\bf (B2)} is implicit and will not be
mentioned explicitly.
For any two positive functions $f$ and $g$,
$f\asymp g$ means that there is a positive constant $c\geq 1$
so that $c^{-1}\, g \leq f \leq c\, g$ in their common domain of
definition.
We will use ``$:=$" to denote a
definition, which is read as ``is defined to be".
For $a, b\in \bR$,
$a\wedge b:=\min \{a, b\}$ and $a\vee b:=\max\{a, b\}$.
For any $x\in \R^d$, $r>0$ and $0<r_1<r_2$, we use $B(x, r)$ to denote the
open ball of radius $r$ centered at $x$ and use $A(x, r_1, r_2)$ to denote the
annulus $\{y\in \R^d: r_1\le |y-x|<r_2\}$.
For a set $V$ in $\R^d$, $|V|$ denotes  the Lebesgue measure of $V$ in $\R^d$.

\section{Preliminaries}\label{s:prelim}

In this section we collect several results concerning subordinators satisfying 
{\bf (A1)}--{\bf (A3)}, the subordinate killed Brownian motion $Y^D$ and its relation with 
the killed subordinate Brownian motion $X^D$.

Let $\phi$ be a Bernstein function and let $S$ be a subordinator with Laplace exponent $\phi$. 
Throughout this section, we always assume that {\bf (A1)}--{\bf (A3)} are in force.
The potential density  $u(t)$ of $S$ satisfies the following two estimates:
\begin{align}
u(t)\le (1-2e^{-1})^{-1}\frac{\phi'(t^{-1})}{t^2\phi(t^{-1})^2}\, ,\quad t>0\, ,\label{e:upper-estimate-u}
\end{align}
and, for every $M>0$ there exists $c_1=c_1(M)>0$ such that
\begin{align}
u(t)\ge c_{1} \frac{\phi'(t^{-1})}{t^2\phi(t^{-1})^2}\, ,\quad 0<t\le
M\, .\label{e:lower-estimate-u}
\end{align}
For the upper estimate see \cite[Lemma A.1]{KM},  
and for the lower one see \cite[Proposition 3.4]{KM}.

The density $\mu(t)$ of the L\'evy measure of $S$ satisfies the following two estimates:
\begin{align}
\mu(t)\le (1-2e^{-1})^{-1}t^{-2}\phi'(t^{-1})\, ,\quad t>0\, ,
\label{e:upper-estimate-mu}
\end{align}
and, for every $M>0$ there exists $c_2=c_2(M)>0$ such that
\begin{align}
\mu(t)\ge c_2 t^{-2}\phi'(t^{-1})\, ,\quad 0<t\le
M\, .\label{e:lower-estimate-mu}
\end{align}
For the upper estimate see \cite[Lemma A.1]{KM}, 
and for the lower one see \cite[Proposition 3.3]{KM}. 
From the last two inequalities it follows that $\mu(t)$ satisfies the doubling property near zero: For every $M>0$ there exists $c_3=c_3(M)>0$ such that
\begin{equation}\label{e:doubling-mu}
\mu(t)\le c_3 \mu(2t)\, , \qquad 0<t\le M\, .
\end{equation}

We will often use the next lemma, cf.~\cite[Lemma 2.1]{KSV16}. 
\begin{lemma}\label{l:pbf}
\begin{itemize}
\item[(a)] For every Bernstein function $\phi$,
\begin{equation}\label{e:wlsc-substitute}
1 \wedge  \lambda\le \frac{\phi(\lambda t)}{\phi(t)} \le 1 \vee \lambda\, ,\quad \textrm{for all }t>0, \lambda>0\, .
\end{equation}
\item[(b)]
If $\phi$ is a special Bernstein function, then for any $a \in [0,2]$,
 $\lambda\mapsto
\lambda^2\frac{\phi'(\lambda)}{\phi(\lambda)^a}$  is an increasing function.
Furthermore, for any $\gamma>2$, $\lim_{\lambda\to 0} \lambda^\gamma\frac{\phi'(\lambda)}{\phi(\lambda)^2}=0$.
\item[(c)]
If $\phi$ is a special Bernstein function, then for
every $d\ge 2$, $\gamma\ge 2$, $\lambda>0$, $b\in (0,1]$ and $a\in [1,\infty)$, it holds that
\begin{equation}\label{e:phi-prime}
\frac{b}{a^{d+\gamma+1}}\frac{\phi'(\lambda^{-2})}{\lambda^{d+\gamma}\phi(\lambda^{-2})^2}\le \frac{\phi'(t^{-2})}{t^{d+\gamma}\phi(t^{-2})^2} \le \frac{a}{b^{d+\gamma+1} }\frac{\phi'(\lambda^{-2})}{\lambda^{d+\gamma}\phi(\lambda^{-2})^2}\, ,\quad\textrm{for all }t\in [b\lambda, a\lambda]\, .
\end{equation}
\end{itemize}
\end{lemma}

Let $W$ be a Brownian motion in $\R^d$ with transition density
$$
p(t,x,y)=(4\pi t)^{-\frac{d}{2}}\exp\left(-\frac{|x-y|^2}{4t}\right)\,
,\qquad t>0,\  x,y\in \R^d\, .
$$
Let $D\subset \R^d$ be a domain and  $W^D$  the Brownian motion $W$ killed upon exiting $D$.
We denote by $p^D(t,x,y)$ the transition density of $W^D$, and by $(P_t^D)_{t\ge 0}$ the corresponding semigroup.
By the strong Markov property,
 $p^D(t,x,y)$ is given by the formula
\begin{align}
\label{e:pD}
p^D(t,x,y)=p(t,x,y)-\E_x[p(t-\tau_D, W_{\tau_D}, y), \tau_D<t]\, , \quad t>0, x,y\in D\, .
\end{align}

Suppose that $W$ is independent of the subordinator $S$.
Recall that $X_t=W_{S_t}$ is the subordinate Brownian motion and $(X^D_t)_{t\ge 0}$
is the subprocess of $X$ killed upon exiting $D$. Then $X$ has a transition density given by
$$
q(t,x,y)=\int_0^{\infty}p(s,x,y)\P(S_t\in ds)\, .
$$
When $X$  is transient, it admits a Green function $G^X(x,y)$ given by 
$$
G^X(x,y)=\int_0^{\infty}q(t,x,y)\, dt=\int_0^{\infty}p(t,x,y)u(t)\, dt\, .
$$
When $d\ge 3$, the Green function $G^X(x,y)$ enjoys the following estimate, 
cf.~\cite[Proposition 4.5]{KM} and \cite[(2.16)]{KSV16}: For
every $M>0$ there exists $c(M)\ge 1$ such that
\begin{equation}\label{e:GX-upper-estimate}
G^X(x,y)\le c(M) \frac{\phi'(|x-y|^{-2})}{|x-y|^{d+2}\phi(|x-y|^{-2})^2}\, ,\qquad \text{for all } x,y\in D, |x-y|\le M\, .
\end{equation}
When $d=2$, under the extra assumptions {\bf (A4)}--{\bf (A5)}, the Green
function estimate above is also valid, cf.~\cite[Proposition 4.5]{KM} and \cite[(2.16)]{KSV16}.
The transition semigroup of $X^D$ will be denoted by $(Q^D_t)_{t\ge 0}$.
Let $Y^D_t=W^D_{S_t}$ be the subordinate killed Brownian motion in $D$ with lifetime denoted by $\zeta$.  The transition semigroup  $(R^D_t)_{t\ge 0}$ of $Y^D$ admits a transition density given by
\begin{equation}\label{e:transition-YD}
r^D(t,x,y)=\int_0^{\infty}p^D(s,x,y)\P(S_t \in ds)\, .
\end{equation}
The subordinate killed Brownian motion $Y^D$ is a transient process, hence admits a Green function
\begin{equation}\label{e:green-YD}
G^{Y^D}(x,y)=\int_0^{\infty}r^D(t,x,y)\, dt=\int_0^{\infty}p^D(t,x,y)u(t)\, dt\, ,
\end{equation}
and clearly $G^{Y^D}(x,y)\le G^X(x,y)$, $x,y\in D$.  

Let $J^{Y^D}(x,y)$ be the jumping density of $Y^D$ 
given by
\begin{align}
\label{e:JY}
J^{Y^D}(x,y) =\int_0^{\infty}p^D(t,x,y)\mu(t)\, dt\, ,
\end{align}
and let $J^X(x, y)=j^X(|x-y|)$ be the L\'evy density of $X$ 
given  by
\begin{align}
\label{e:JX}
j^X(|x-y|)=\int^\infty_0p(t, x, y)\mu(t)\, dt.
\end{align}
Clearly $J^{Y^D}(x,y) \le j^X(|x-y|)$, $x, y \in D$. Furthermore, similarly as in \eqref{e:GX-upper-estimate} 
(cf.~\cite[Proposition 4.2]{KM}),
either when $d\ge 3$ or when $d=2$ and  {\bf (A4)} holds, there exists $c(M)\ge 1$ such that
\begin{equation}\label{e:JD-estimate}
c(M)^{-1} \frac{\phi'(|x-y|^{-2})}{|x-y|^{d+2}}\le 
j^X(|x-y|)\le c(M) \frac{\phi'(|x-y|^{-2})}{|x-y|^{d+2}}\, ,\qquad |x-y|\le M\, .
\end{equation}

We recall now some properties of $Y^D$ connected to its relationship with $X^D$. It is shown in \cite{SV08} that $Y^D$ can be realized as $X^D$ killed at a terminal time, and consequently, cf.~\cite[Proposition 3.1]{SV08}, the semigroup $(R^D_t)_{t\ge 0}$ is subordinate to the semigroup $(Q^D_t)_{t\ge 0}$ in the sense that $R^D_t f(x)\le Q^D_t f(x)$ for all Borel $f:D\to [0,\infty)$, all $t\ge 0$ and all $x\in D$. 
As a consequence,
 cf.~\cite[Proposition 4.5.2]{FOT} and \cite[Proposition 3.2]{SV03}, 
we have the following relation between 
the density $\kappa^{Y^D}$ of the killing measure of $Y^D$ and the
density $\kappa^{X^D}$ of the killing measure of $X^D$:
\begin{equation}\label{e:killing-functions-relation}
\kappa^{Y^D}(x) \ge \kappa^{X^D}(x)\, , \qquad x\in D\, .
\end{equation}
Furthermore, it follows from \cite{Sztonyk} that if $U\subset D$ is Lipschitz, then
\begin{equation}\label{e:exit-time-XD}
\P_x(X^D_{\tau_U}\in \partial U)=0.
\end{equation}
Here $\tau_U$ denotes the first exit time of the process $X^D$ from $U$. 
For simplicity, we will sometimes use the same notation for the exit time from $U$ for other processes when it will be clear from the context which process we have in mind. 
Since $Y^D$ can be realized as  $X^D$ killed at a terminal time, it follows immediately from \eqref{e:exit-time-XD} that we also have
\begin{equation}\label{e:exit-time-YD}
\P_x(Y^D_{\tau_U}\in \partial U)=0.
\end{equation}
Again, by using that $Y^D$ can be realized as $X^D$ killed at a terminal time, it follows from \cite[Corollary 4.2(i)]{SV08} that 
\begin{equation}\label{e:YD-killed-inside}
\P_x(Y^D_{\zeta-}\in D)=1\, \qquad \text{for all }x\in D\, ,
\end{equation}
i.e., the process $Y^D$ dies inside $D$ almost surely (and not at the boundary $\partial D$).

For any open set $U\subset D$,  let $Y^{D,U}$ be  the subprocess of $Y^D$ killed upon exiting $U$. Define 
\begin{equation}\label{e:hkos1}
r^{D,U}(t,x,y)  =  r^D(t,x,y) -  \E_x [ r^D(t - \tau_U,Y^D_{\tau_U},y) : \tau_U < t]\quad  
t>0,\  x,y \in U.
\end{equation}
Then, by the strong Markov property, $r^{D, U}(t,x,y)$ is the transition density of $Y^{D, U}$.
Let $G^{Y^D}_U(x,y)= \int_0^\infty  r^{D,U}(t, x,y) dt$ denote 
the Green function of $Y^{D, U}$.
Since $r^{D,U}(t,x,y)\le r^D(t,x,y)$ for all $t\ge 0$, we clearly have that $G^{Y^D}_U(x,y)\le G^{Y^D}(x,y)$ for all $x,y\in U$. 

Let again $U\subset D$ be an open set, $W^U$ the Brownian motion $W$ killed upon exiting $U$, and $Y^U_t:=W^U_{S_t}$ the corresponding subordinate killed Brownian motion with its Green function denoted by $G^{Y^U}(x,y)$, $x,y\in U$.  
Similar to the fact that the semigroup $(R^D_t)_{t\ge 0}$ of $Y^D$  is
subordinate to the semigroup $(Q^D_t)_{t\ge 0}$ of $X^D$, we also have that
the semigroup of $Y^U$ is subordinate to the semigroup of $Y^{D,U}$. 
In fact, this follows directly from \cite[Proposition 3.1]{SV08} by considering  $Y^D$ as the underlying process. After integrating over time, it follows that
\begin{equation}\label{e:UB-smaller-UDB}
G^{Y^U}(x,y)\le G^{Y^D}_U(x,y)\, \qquad \text{for all }x,y\in U\, .
\end{equation}

Note that
$$
p^D(t,z,w) \le p(t,x/2,y/2)=(4\pi t)^{-\frac{d}{2}}\exp\left(-\frac{|x-y|^2}{16t}\right), 
$$
for all $z,w$, $|z-w| \ge |x-y|/2$. 
Hence by the continuity of $(x,y)\mapsto p^D(t,x,y)$ and the dominated convergence theorem, we see from \eqref{e:transition-YD}  that $(x,y)\mapsto r^D(t,x,y)$ is continuous on $D\times D$. Consequently, $(R^D_t)_{t\ge 0}$  satisfies the strong Feller property.

Recall that for an open set $D\subset \R^d$, a point $x\in \partial D$ is said to be regular for $D^c$ with respect to Brownian motion
 if $\P_x(\tau^W_D=0)=1$,  where $\tau^W_D:=\inf\{t>0: W_t \notin D\}$ is the first exit time of $D$ for $W$. 
 It is well known that, if  all points on $\partial D$ are regular for $D^c$
 with respect to Brownian motion, 
 then $W^D$ is a Feller process on $D$ (see \cite[(7) in Theorem 2.4]{CZ}).
 By Phillips' theorem (see \cite[Proposition 13.1]{SSV}), in this case
$Y^D$ is also a Feller process, i.e.,  $(R^D_t)_{t\ge 0}$ is a strongly continuous contraction semigroup on $(C_0(D), \|\cdot \|_\infty)$. 
Thus $Y^D$ enjoys both the Feller and strong Feller property.
\begin{lemma}\label{l:spfp}
Suppose that $D\subset \R^d$ is a domain.
Let $K$ be a compact set and $G$ be an open  set with $K \subset G\subset
\overline G\subset D$.
Then for any $\varepsilon > 0$, there is $t_0>0$ such that
\begin{equation}
\sup_{x \in K} \P_x( \tau^{Y^D}_G \leq t_0 ) \leq \varepsilon.
\end{equation}
In particular, $\P_x( \tau^{Y^D}_G \leq t ) \rightarrow 0$ 
uniformly on $K$ as $t \rightarrow 0$.
\end{lemma}

\pf In the case when all points on $\partial D$ are regular for $D^c$ with respect to Brownian motion, $Y^D$ is a Feller process and this lemma follows from  \cite[Lemma 2]{Ch}. In the general case,  we can take an open set $U\subset D$ with all points on $\partial U$ being regular for 
$U^c$ such that $G\subset U$. Let $\tau^{Y^U}_G$ be the first time the process
$Y^U$ exits $G$. Then $\tau^{Y^U}_G\le \tau^{Y^D}_G$. Combining this with 
\cite[Lemma 2]{Ch} we immediately get the conclusion of the lemma.
\qed
\begin{prop}\label{p:cont}
Suppose that {\bf (A5)} holds,  $D\subset \R^d$ is a domain and that $U$ is 
an open set with $\overline U\subset D$.
Then for each $t>0$, 
$(x, y) \mapsto r^{D,U}(t,x,y)$ 
is continuous on $U \times U$, and  for each $y  \in U$,  $x \mapsto G^{Y^D}_U(x,y)$ is continuous on $U \setminus \{y\}$.
\end{prop}
\pf
Note that  for all $ \beta>0$, 
\begin{align}
\label{e:r: bound}
& \sup_{|x-y| \geq \beta, t>0} r^{D}(t,x,y)  \le \sup_{|x-y| \geq \beta, t>0} q (t,x,y) 
\le \sup_{ t>0}   \int_{(0, \infty)}\sup_{|x| \geq \beta}p(s,x)\, \P(S_t\in ds)\nn\\ &
\le \sup_{ t>0}   \int_{(0, \infty)}\sup_{s>0, |x|= \beta}p(s,x)\, \P(S_t\in ds)
= \sup_{s>0, |x|= \beta}p(s,x) <\infty.
\end{align}
Using Lemma \ref{l:spfp}, the strong Feller property of 
$(R^D_t)_{t\ge 0}$ and \eqref{e:r: bound}, 
one can follow the proof of \cite[Theorem 2.4]{CZ} line by line and 
show that, for each $t>0$, 
$(x, y) \to \E_x [ r^D(t - \tau_U,Y^D_{\tau_U},y) : \tau_U < t]$  
is continuous on $U \times U$. 
Thus, by \eqref{e:hkos1} and the already proven fact that 
$(x,y)\mapsto r^D(t,x,y)$ is continuous, we see that 
$(t,x, y) \mapsto r^{D,U}(x,y)$   
is continuous on $U \times U$. 

Fix $y \not=x_0  \in U$ and consider $x \in B(x_0, \eps) \subset  D \setminus B(y, \eps)$. Then,  since 
$|x-y| \ge |x_0-y|/2$,  we have $ r^{D,U}(t, x,y) \le q(t, x,y) \le q(t, x_0/2,y/2)$.
Moreover,  $ t \mapsto q(t, x_0/2,y/2) \in L^1(0, \infty)$ 
because $X$ is transient. Therefore, by applying the dominated convergence theorem to 
$
G^{Y^D}_U(x,y)=\int_0^\infty  r^{D,U}(t, x,y) dt,
$
  we conclude that $x \mapsto G^{Y^D}_U(x,y)$ 
 is continuous on $B(x_0, \eps)  \subset  D \setminus B(y, \eps)$.
 \qed


\section{Harnack inequality}\label{s:hi}
The goal of this section is to prove 
the following scale invariant HI for
non-negative functions harmonic with respect to $Y^D$ 
when the domain $D$ satisfies 
{\bf (B1)} and {\bf (B2)}.
\begin{thm}
[Harnack inequality]\label{uhp}
Assume that 	{\bf (A1)}--{\bf (A3)} hold and
that $D\subset \R^d$ is a domain satisfying {\bf (B1)}--{\bf (B2)}.
There exists a constant $C>0$ such that for any
$r\in (0,1]$ and $B(x_0, r) \subset D$ and any function $f$ which is non-negative in
$D$ and harmonic in $B(x_0, r)$ with respect to $Y^D$, we have
$$
f(x)\le C f(y), \qquad \text{ for all } x, y\in B(x_0, r/2).
$$
\end{thm}

As already mentioned in the introduction, the proof is modeled after the proof of 
the HI in \cite{KM}.
In the remaining part of this section we assume that 
the assumptions	{\bf (A1)}--{\bf (A3)}  hold true.
\begin{lemma}\label{l:ilb4BMheatkernel}
Let $D$ be a  domain in $\R^d$ satisfying {\bf (B2)} and $b> 0$ a constant. 
There exists $C>0$ depending on $b$ such that
\begin{equation}\label{e:ilb4BMheatkernel}
p^D(t, x, y)\ge Ct^{-d/2}e^{-M|x-y|^2/t}, \qquad \sqrt{t}\le b(\delta_D(x)\wedge\delta_D(y))
\wedge 1.
\end{equation}
\end{lemma}

\pf If $\sqrt{t}\le b\delta_D(x)$,  then
$
\P_x(t<\tau^W_D)\ge \P_x(t<\tau^W_{B(x, b^{-1}t^{1/2})})\ge c
$
for some constant $c=c(b)>0$ independent of $x$ and $t$. Now the conclusion of this
lemma follows immediately from {\bf (B2)}.
\qed
\begin{remark}\label{r:ilb4BMheatkernel}
{\rm 
By using a simple 
chaining  argument, the conclusion of this lemma is actually true under 
the following alternative (and seemingly weaker) assumption:
There exist 
$\lambda_1 \in [1, \infty)$ and $\lambda_2 \in (0, 1]$ such that for all $r \le 1$ and  $x,y\in D$  with 
$\delta_D(x)\wedge \delta_D(y)\ge r$ there exists a length parameterized rectifiable curve $l$ connecting $x$ to $y$ with the length $|l|$ of $l$  less than or equal to $\lambda_1|x-y|$ and $\delta_D(l(u))\geq\lambda_2 r$ for $u\in[0,|l|].$ 
Since in the proofs of Lemmas \ref{l:ilb4jk} and 
\ref{l:green-ball-lower-bound}--\ref{l:tauB},
and Propositions \ref{p:U-exit-upper-bound}--\ref{p:U-exit-lower-bound}, only 
\eqref{e:ilb4BMheatkernel} is used, we can replace {\bf (B2)} in these
results by the alternative assumption above.}
\end{remark}
\begin{lemma}\label{l:ilb4jk}
Let $D$ be  a domain satisfying {\bf (B2)} and $\varepsilon_0>0$ 
a constant.  There exists a constant $C=C(\varepsilon_0)\in (0, 1)$ such that for
every $x_0\in D$ 
and $r\le 1/2$ satisfying $B(x_0, (1+\varepsilon_0)r)\subset D$, we have
\begin{equation}\label{e:ilb4jk}
CJ^X(x, y)\le J^{Y^D}(x, y)\le J^X(x, y), \qquad x, y\in B(x_0, r).
\end{equation}
\end{lemma}

\pf  The second inequality in \eqref{e:ilb4jk} is obvious, so we will only prove
the first inequality. It follows from Lemma \ref{l:ilb4BMheatkernel} (with $b=2/\varepsilon_0$) that there exist $c_1, c_2>0$ such that for
any $t\le 4r^2$,
$$
p^D(t, x, y)\ge c_1t^{-d/2}e^{-c_2|x-y|^2/t}.
$$
Thus for $x, y\in B(x_0, r)$,
\begin{align*}
&J^{Y^D}(x, y)\ge\int^{|x-y|^2}_0p^D(t, x, y)\mu(t)dt
\ge \mu(|x-y|^2)\int^{|x-y|^2}_0p^D(t, x, y)dt\\
&\ge c_1\mu(|x-y|^2)\int^{|x-y|^2}_0t^{-d/2}e^{-c_2|x-y|^2/t}dt
\ge c_3|x-y|^{-d-2}\phi'(|x-y|^{-2})\\
&\ge c_4J^X(x, y),
\end{align*}
where in the second inequality of the second line we used \eqref{e:lower-estimate-mu}, and in the last inequality we used \eqref{e:JD-estimate}.
\qed
\begin{prop}\label{p:estimate-of-J-away}
Let $D$ be  a domain satisfying {\bf (B1)}--{\bf (B2)}. 
For every $\eps_0 \in (0,1]$, 
there exists a constant $C\ge 1$ depending on the constants 
from {\bf (B1)}--{\bf (B2)} and $\eps_0$
 such that for all $x_0\in D$ and all $r \le 1$ satisfying $B(x_0, (1+\eps_0)r)\subset D$, it holds that
\begin{equation}\label{e:estimate-of-J-away}
J^{Y^D}(z,x_1)\le C J^{Y^D}(z,x_2)\, ,\quad x_1,x_2\in B(x_0,r), \ \ z\in D\setminus B(x_0, (1+\eps_0)r)\, .
\end{equation}
\end{prop}

\pf Suppose that $r, \eps_0 \le 1$, $B(x_0, (1+\eps_0)r)\subset D$ and $x_1, x_2\in B(x_0, r)$.  
Then for $t<\eps_0^2r^2$, we have
$$
\P_{x_2}(t<\tau^W_D)\ge \P_{x_2}(t<\tau^W_{B(x_2, t^{1/2})})\ge c_1
$$
for some constant $c_1=c_1(\eps_0)>0$ independent of $x_2$ and $t$.
By combining with {\bf (B2)} we see that there exists $c_2>1$ such that for
$z\in D$ and $t<\eps_0^2r^{2}$, 
\begin{align*}
&p^D(t, z, x_1)\le c_2 \P_z(t<\tau^W_D)t^{-d/2}e^{-\frac{|z-x_1|^2}{Mt}},\\
&p^D(t, z, x_2)\ge c_2^{-1}\P_z(t<\tau^W_D)t^{-d/2}e^{-\frac{M|z-x_2|^2}{t}}.
\end{align*}
Now suppose that $z\in D\setminus B(x_0,(1+\eps_0)r)$ so that 
$$
\frac{\eps_0}{1+\eps_0}|z-x_0|\le |z-x_i|\le \left(1+\frac1{1+\eps_0}\right)|z-x_0|,
\qquad i=1, 2. 
$$
Then
\begin{align*}
&\int^{\eps_0^2r^2}_0p^D(t, z, x_1)\mu(t)dt\le c_2\int^{\eps_0^2r^2}_0\P_z(t<\tau^W_D)t^{-d/2}e^{-\frac{|z-x_2|^2}{16Mt}}\mu(t)dt\\
&=c_2(16M^2)^{d/2-1}\int^{16M^2\eps_0^2r^2}_0\P_z(\frac{t}{16M^2}<\tau^W_D)t^{-d/2}e^{-\frac{M|z-x_2|^2}{t}}\mu(t/(16M^2))dt.
\end{align*}
Using the doubling property of $t\mapsto \P_z(t<\tau^W_D)$  
({\bf (B1)}) and 
the doubling property \eqref{e:doubling-mu} of $\mu(t)$ near the origin,
we get that there exists $c_3>1$ such that
\begin{equation}\label{e:rs1}
\int^{\eps_0^2r^2}_0p^D(t, z, x_1)\mu(t)dt\le c_3\int^{16M^2\eps_0^2r^2}_0p^D(t, z, x_2)\mu(t)dt.
\end{equation}

Using the parabolic Harnack principle (see, for instance, \cite[Theorem 5.4]{FS}),
we get that there exists $c_4>1$ such that
\begin{align*}
&\int^\infty_{\eps_0^2r^2}p^D(t, z, x_1)\mu(t)dt=\sum^\infty_{n=1}\int^{(n+1)\eps_0^2r^2}_{n\eps_0^2r^2}p^D(t, z, x_1)\mu(t)dt\\
&\le c_4\sum^\infty_{n=1}\int^{(n+1)\eps_0^2r^2}_{n\eps_0^2r^2}p^D(t+\frac{\eps_0^2r^2}2, z, x_2)\mu(t)dt\\
&=c_4\sum^\infty_{n=1}\int^{(n+\frac32)\eps_0^2r^2}_{(n+\frac12)\eps_0^2r^2}p^D(t, z, x_2)\mu(t-\frac{\eps_0^2r^2}2)dt\, .
\end{align*}
If $t<1$, then by \eqref{e:doubling-mu} we have
$\mu(t-\eps_0^2r^2/2)\le \mu(t/2)\le c_5 \mu(t)$  with $c_5\ge 1$. 
If $t\ge 1$, then $\mu(t-\eps_0^2r^2/2)\le \mu(t-1/2)$. 
By \eqref{e:CBmu0}, there exists $c_6 \ge 1$ such that $\mu(s)\le c_6 \mu(s+1/2)$ for all $s>1/2$. Hence, $\mu(t-1/2)\le c_6\mu(t)$. With $c_7=c_5\vee c_6$,  we conclude that $\mu(t-\eps_0^2r^2/2)\le c_7\mu(t)$ for all $t\ge 3r^2/2$. Hence,
\begin{align}
&\int^\infty_{\eps_0^2r^2}p^D(t, z, x_1)\mu(t)dt\le c_4 c_7\sum^\infty_{n=1}\int^{(n+\frac32)\eps_0^2r^2}_{(n+\frac12)\eps_0^2r^2}p^D(t, z, x_2)\mu(t)dt\nonumber \\
&\le c_8\int^\infty_{(3\eps_0^2r^2)/2}p^D(t, z, x_2)\mu(t)dt\, .
\label{e:rs2}
\end{align}

Combining \eqref{e:rs1} and \eqref{e:rs2}, we get that there exists $c_9>1$ such that
$$
\int^\infty_0p^D(t, z, x_1)\mu(t)dt\le c_9\int^\infty_0p^D(t, z, x_2)\mu(t)dt\, ,
$$
which finishes the proof. \qed
\begin{lemma}\label{l:green-ball-lower-bound}
Suppose $d\ge 2$. Let $D$ be a  domain in $\R^d$ satisfying {\bf (B2)}. In the case
$d=2$, we also assume that {\bf (A4)} and {\bf (A5)} hold.
There exist $a\in (0,1/3)$  and $C>0$ such that for every $x_0\in D$ and every $r\in (0,1)$ satisfying $B(x_0,r)\subset D$, 
\begin{equation}\label{e:green-ball-lower-bound}
G^{Y^D}_{B(x_0,r)}(x,y)\ge C |x-y|^{-d-2}\frac{\phi'(|x-y|^{-2})}{\phi(|x-y|^{-2})^2}\, ,\quad x,y\in B(x_0,ar)\, .
\end{equation}
\end{lemma}
\pf 
It follows from Lemma \ref{l:ilb4BMheatkernel}
that there exist $c_1, c_2>0$ such that 
$$
p^D(t,x,y) \ge c_1 t^{-d/2} e^{-c_2|x-y|^2/t}, 
\quad \delta_D(x) \wedge \delta_D(y)\wedge 1 \ge \sqrt{t}. 
$$
Note that for $a\in (0,1)$, $(1-a)/2a > 1$ if and only if $a< 1/3$.
Choose $a\in (0,1/3)$ and let $x,y\in B(x_0,ar)$ with $r\in (0,1)$ satisfying $B(x_0,r)\subset D$. Then $\delta_D(x),\delta_D(y)\ge (1-a)r$ and $|x-y|\le 2ar$. 
Thus $\delta_D(x) \wedge \delta_D(y) \ge ((1-a)/2a) |x-y|$.
Therefore, by \eqref{e:lower-estimate-u},
\begin{align}\label{e:green-ball-lower-bound-2}
&G^{Y^D}(x,y)\ge 
\int_{0}^{|x-y|^2} p^D(t,x,y) u(t) dt \nn\\
&\ge c_1 u(|x-y|^2) \int_{0}^{|x-y|^2}  t^{-d/2} e^{-c_2|x-y|^2/t}  dt \nn \\
&\ge 
c_3|x-y|^{-d-2}\frac{\phi'(|x-y|^{-2})}{\phi(|x-y|^{-2})^2}\, .
\end{align}
For simplicity, let $B=B(x_0,r)$. For $y\in B(x_0,ar)$,  we have
$$
|Y^D_{\tau_B}-y|\ge |Y^D_{\tau_B}-x_0|-|y-x_0|\ge r-ar \ge \frac{1-a}{2a}|x-y|\ge |x-y|\, .
$$
By the above inequalities, \cite[Proposition 4.5]{KM} and
Lemma \ref{l:pbf}(b), we have
\begin{eqnarray}
G^{Y^D}(Y^D_{\tau_B}, y) &\le &G^X(Y^D_{\tau_B}, y) \le 
G^X\left(\frac{1-a}{2a}x, \frac{1-a}{2a}y\right)\nn\\
&\le &c_4\left(\frac{1-a}{2a}|x-y|\right)^{-d+2}  (\frac{1-a}{2a}|x-y|)^{-4}\frac{\phi'(( \frac{1-a}{2a}|x-y|)^{-2})}{\phi( (\frac{1-a}{2a}|x-y|)^{-2})^2} \label{e:neq2}\\
&\le &c_4\left(\frac{1-a}{2a}|x-y|\right)^{-d+2} |x-y|^{-4}\frac{\phi'(|x-y|^{-2})}{\phi(|x-y|^{-2})^2}\nn\\
&=& c_4 \left(\frac{2a}{1-a}\right)^{d-2} |x-y|^{-d-2}\frac{\phi'(|x-y|^{-2})}{\phi(|x-y|^{-2})^2}\, .
 \label{e:neq3}
\end{eqnarray}

When $d \ge 3$, choose $a\in (0,1/3)$ small enough so that 
$c_3-c_4\left(\frac{2a}{1-a}\right)^{d-2}\ge \frac{1}{2c_4}$. 
Then, by \eqref{e:neq3}, for $x,y\in B(x_0,ar)$ we have
\begin{align*}
&G^{Y^D}_{B(x_0,r)}(x,y)=G^{Y^D}(x,y)-\E_x G^{Y^D}(Y^D_{\tau_B},y)\\
&\ge
\left(c_3 -  c_4 \left(\frac{2a}{1-a}\right)^{d-2} \right)  |x-y|^{-d-2}
\frac{\phi'(|x-y|^{-2})}{\phi(|x-y|^{-2})^2} \\
&\ge 
\frac{1}{2c_4} |x-y|^{-d-2}\frac{\phi'(|x-y|^{-2})}{\phi(|x-y|^{-2})^2}\, .
\end{align*}

When $d = 2$, 
by \eqref{e:neq2}, for $x,y\in B(x_0,ar)$, we have
\begin{align}
&G^{Y^D}_{B(x_0,r)}(x,y) 
\ge
 c_3 \left(1 - \frac{c_4}{c_3}
\frac{\psi\left(\left(\tfrac{2a}{1-a}
\right)^2|x-y|^ { -2 } \right)}{\psi(|x-y|^{-2})} \right)  
|x-y|^{-4}\frac{\phi'(|x-y|^{-2})}{\phi(|x-y|^{-2})^2},
\label{e:Ulow2} \end{align}
 where
$$
 \psi(\lambda)=\lambda^2\tfrac{\phi'(\lambda)}{\phi(\lambda)^2}, \qquad \lambda>0\,.
 $$
Since the $\delta'$ in {\bf (A4)} is in $(0, 2\delta)$, we can choose $\varepsilon>0$
small enough so that  $2\delta-\delta'-2\varepsilon>0$.
By \cite[Lemma 3.2(ii)]{KM}, there exists $c_\eps>0$ such that
\begin{align}\label{e:l3211}
	\frac{\phi(\lambda s)}{\phi(\lambda)}\leq c_\eps\,s^{1-\delta +\varepsilon}\
\text{ for all }\ \lambda\geq 1 \ \text{ and }\ s\geq 1\,.
\end{align}
Choose $a<\frac{1}{3}$ small enough so that 
$$
\frac{c_4}{c_3}\frac1{\sigma'}c_\eps^2
\left(\frac{2a}{1-a}\right)^{2\delta-\delta'-2\varepsilon}\leq \frac{1}{2},
$$
where $\sigma'$ is the constant in {\bf (A4)}.
Then,  by {\bf (A4)} and \eqref{e:l3211}, we have 
\begin{align}
\frac{c_4}{c_3}
\frac{\psi\left(\left(\tfrac{2a}{1-a}\right)^2|x-y|^{-2}\right)}
{\psi(|x-y|^{-2})}
=& \frac{c_4}{c_3}
 \left(\tfrac{2a}{1-a}\right)^{4}
 \frac{\phi'\left(\left(\tfrac{2a}{1-a}\right)^2|x-y|^{-2}\right)\phi(|x-y|^{-2}
)^2}{\phi'(|x-y|^{-2}
)\phi\left(\left(\tfrac{2a}{1-a}\right)^2|x-y|^{-2}\right)^2}\nn\\
 \le& 
\frac{c_4}{c_3}\frac1{\sigma'}
c_\eps^2  \left(\frac{2a}{1-a}\right)^{2\delta-\delta'-2\varepsilon} 
\le \frac12. \label{eq:gr-39}
\end{align}
The lemma for $d = 2$ now follows from \eqref{e:Ulow2} and  \eqref{eq:gr-39}. 
\qed
\begin{lemma}\label{l:tauB}
Suppose $d\ge 2$. Let $D$ be a  domain in $\R^d$ satisfying {\bf (B2)}. When 
$d=2$, we also assume 
that {\bf (A4)} and {\bf (A5)} hold.
There exists $C>0$ such that for every $x_0\in D$ and every $r\in (0,1)$ 
with $B(x_0,r)\subset D$,
$$
\E_x \tau^{Y^D}_{B(x_0,r)} 
\ge C \phi(r^{-2})^{-1}\, ,\quad x\in B(x_0,ar/2)\, ,
$$
where $a\in (0,1/3)$ is the constant from Lemma \ref{l:green-ball-lower-bound}.
\end{lemma}
\pf Using Lemmas \ref{l:pbf} and \ref{l:green-ball-lower-bound}, 
the proof is exactly the same as that of \cite[Proposition 5.2]{KM}. \qed

For all $x\in D$ and $r>0$ with $B(x,r)\subset D$, and all non-negative functions 
$f$, we define 
$$
(\UU_r f)(x)=\frac{\E_x[f(Y^D(\tau_{B(x,r)}))]-f(x)}{\E_x\tau^{Y^D}_{B(x,r)}}.
$$
Then Example 5.4, Remark 5.5 and Proposition 5.6 from \cite{KM} are valid for $Y^D$.
Thus, with $\eta(z):=\E_z\tau^{Y^D}_{B(x,r)}$,
\begin{align}\label{e:new2}
\sU_s \eta(y)=-1 \quad \text{ for any } y\in B(x,r)  \text{ and } s<r-|y-x|,
\end{align}
and, for any function $h\colon D \rightarrow [0,\infty)$ which is  harmonic in 
a bounded open set $U\subset D$, 
\begin{align}\label{rem:harm} 
 (\sU_sh)(x)=0\ \text{ for all }\ x\in U.
\end{align}
If  $f(x_0)\leq f(x)$ for all $x\in \R^d$ then $(\sU_r f)(x_0)\geq 0$. Thus we have the following type of maximum principle.
\begin{prop}[Maximum principle]\label{prop:max_pr}
Assume that there exist $x_0\in D$ and $r>0$ 
such that $(\sU_rf)(x_0)<0$. Then
\begin{equation}\label{eq:gr-40}
 f(x_0)>\inf_{x \in D}f(x)\,.
\end{equation}
\end{prop}

The next three results are  also valid when $d=2$ if we further assume that
{\bf (A4)} and {\bf (A5)} hold. Since we will only need these results for $d\ge 3$, we state them for $d\ge3$ without the extra assumptions. 
From now until the proof of Theorem \ref{uhp}, we will assume that
$d\ge 3$.

By using  $G^{Y^D}(x,y)\le G^X(x,y)$ and 
the estimate \eqref{e:GX-upper-estimate} with $M=2$ (see also \cite{Gr})
  one gets that  there exists $c_1>0$ such that for every $x_0\in D$ 
and every $r\in (0,1)$ with $B(x_0,r)\subset D$,  we have
\begin{equation}\label{e:exit-time-est}
\E_x \tau^{Y^D}_{B(x_0,r)} \le c_1 \phi(r^{-2})^{-1}\, ,\quad x\in B(x_0,r)\, .
\end{equation}

Let 
\begin{equation}\label{e:defofgandj}
g(r):=c_2 r^{-d-2}\frac{\phi'(r^{-2})}{\phi(r^{-2})^2}\, , \quad j(r):=c_3 r^{-d-2} \phi'(r^{-2})\, , 
\qquad r>0\, ,
\end{equation}
where $c_2$ is the constant $c(2)$ from \eqref{e:GX-upper-estimate} and $c_3$ is the
constant $c(2)$ from \eqref{e:JD-estimate}. Then for all $x,y\in D$ 
with $|x-y|<2$,
\begin{equation}\label{e:newrs1}
G^{Y^D}(x,y)\le G^X(x, y) \le   g(|x-y|)\, , \qquad J^{Y^D}(x,y)\le  j^X(|x-y|)  \le   j(|x-y|)\, .
\end{equation}
In particular, for any domain $D$, any $x_0\in \R^d$, any $s\in (0, 1)$ with $B(x_0, s)\subset D$,
\begin{equation}\label{e:newrs}
G^{Y^D}(x,y)\le g(|x-y|)\le g\left(\frac{s}{8}\right),
 \qquad x\in B\left(x_0, \frac{s}2\right),\  y\in B\left(x,\frac{s}{8}\right)^c.
\end{equation}

Using 
\eqref{e:exit-time-YD},  \eqref{e:new2},  \eqref{rem:harm}, \eqref{e:exit-time-est},  \eqref{e:newrs1}, \eqref{e:newrs},
 Lemma \ref{l:tauB}, 
and Propositions \ref{p:cont} and \ref{prop:max_pr}, we can repeat
the argument in the proof of  \cite[Proposition 5.7]{KM} with functions 
$
\eta(z)=\E_z\tau^{Y^D}_{B(x_0,r)}$  and $ h(z)=G^{Y^D}_{B(x_0,r)}(x,z)
$
 to get the following result.
Since the argument is basically the same, we omit the proof.
\begin{prop}\label{p:U-exit-upper-bound}
Suppose $d\ge 3$.
Let $D$ be a  domain in $\R^d$ satisfying {\bf (B2)}.
There exists a constant $C>0$ such that for all $r\in (0,1)$ and all $x_0\in D$ 
satisfying $\overline B(x_0,r)\subset D$,
\begin{equation}\label{e:U-exit-upper-bound}
G^{Y^D}_{B(x_0,r)}(x,y)\le C r^{-d-2}\frac{\phi'(r^{-2})}{\phi(r^{-2})}\E_y \tau^{Y^D}_{B(x_0,r)}\, ,\quad x\in 
B(x_0,ar/4), y\in A(x_0,2^{-1}ar,r)\, ,
\end{equation}
where  $a\in (0,1/3)$ is the constant from Lemma \ref{l:green-ball-lower-bound} and $b=a/2$.
\end{prop}

Using \eqref{e:JD-estimate}, \eqref{e:exit-time-YD}, \eqref{e:new2}, \eqref{rem:harm} and  \eqref{e:exit-time-est},
Lemmas \ref{l:pbf}, \ref{l:ilb4jk}, \ref{l:green-ball-lower-bound} and \ref{l:tauB},  
and Propositions \ref{prop:max_pr} and \ref{p:U-exit-upper-bound}, we can repeat
the argument in the proof of  \cite[Proposition 5.8]{KM} to get the following result.
Since the argument is basically the same, we omit the proof.
\begin{prop}\label{p:U-exit-lower-bound}
Suppose $d\ge3$.
Let $D$ be a  domain in $\R^d$ satisfying {\bf (B2)}.
Then for every $\eps_0 \in (0,1]$, there exist constants $C>0$ and $b\in (0,1)$ such that for any $r\in (0,1)$ and any $x_0\in D$ 
satisfying $B(x_0,(1+\eps_0) r)\subset D$,
\begin{equation}\label{e:U-exit-lower-bound}
G^{Y^D}_{B(x_0,r)}(x,y)\ge C r^{-d-2}\frac{\phi'(r^{-2})}{\phi(r^{-2})}\E_y \tau^{Y^D}_{B(x_0,r)}\, ,\quad x\in B(x_0,br), y\in B(x_0,r)\, ,
\end{equation}
where  $a\in (0,1/3)$ is the constant from Lemma \ref{l:green-ball-lower-bound}.
\end{prop}
\begin{corollary}\label{c:U-exit-bounds}
Suppose $d\ge 3$.
Let $D$ be a  domain in $\R^d$ satisfying {\bf (B2)}.
For every $\eps_0 \in (0,1]$, there exist
constants $C>1$ and $b_1,b_2\in (0,1/3)$ satisfying $2b_1<b_2$ such that for all $r\in (0,1)$ and all $x_0$ satisfying that $B(x_0, (1+\eps_0)r)\subset D$,
$$
C^{-1} r^{-d-2}\frac{\phi'(r^{-2})}{\phi(r^{-2})}\E_y \tau^{Y^D}_{B(x_0,r)} \le G^{Y^D}_{B(x_0,r)}(x,y)
\le C r^{-d-2}\frac{\phi'(r^{-2})}{\phi(r^{-2})}\E_y \tau^{Y^D}_{B(x_0,r)}\, ,
$$ 
for all $ x\in B(x_0,b_1 r)$ and all $ y\in A(x_0, b_2 r, r)$.
\end{corollary}

For any open set $U\subset D$, let
$$
K^{D,U}(x,z):=\int_U G^{Y^D}_U(x,y)J^{Y^D}(y,z)\, dy\, ,\quad x\in U, z\in \overline{U}^c\cap D
$$
be the Poisson kernel of $Y^D$ on $U$.

\noindent
{\bf Proof of Theorem \ref{uhp}:}
We will first assume that $d\ge 3$.
Without loss of generality  (by considering $(1-\eps_0)r$ instead of 
$r$ if necessary),  we may assume that $B(x_0, (1+\eps_0)r)\subset D$ for some
small $\eps_0 \in (0,1)$.
Let $b_1, b_2\in (0,1/3)$ be as in Corollary \ref{c:U-exit-bounds}.
We will first show that there exists a constant 
$c>0$, independent of $x_0$ and $r$,
such that for all $x_1,x_2\in B(x_0, b_1r)$ and $z\in \overline{B(x_0,r)}^c\cap D$,
\begin{equation}\label{e:pk-est}
K^{D,B(x_0,r)}(x_1,z)\le c K^{D,B(x_0,r)}(x_2,z)\, .
\end{equation}
Note that
\begin{align*}
&K^{D,B(x_0,r)}(x_1,z)\\
&=\int_{B(x_0,b_2 r)}G^{Y^D}_{B(x_0,r)}(x_1,y)J^{Y^D}(y,z)\, dy+\int_{A(x_0,b_2 r,r)}G^{Y^D}_{B(x_0,r)}(x_1,y)J^{Y^D}(y,z)\, dy\\
&=I_1+I_2\, .
\end{align*}
In order to estimate $I_2$ we use Corollary \ref{c:U-exit-bounds} to get 
\begin{align*}
&I_2 \le  c_1\frac{r^{-d-2}\phi'(r^{-2})}{\phi(r^{-2})}\int_{A(x_0,b_2r, r)}\E_y \tau^{Y^D}_{B(x_0,r)}J^{Y^D}(y,z)\, dy\\
&\le  c_1^2 \int_{A(x_0,b_2r, r)}G^{Y^D}_{B(x_0,r)}(x_2,y)J^{Y^D}(y,z)\, dy \le c_1^2 K^{D,B(x_0,r)}(x_2,z)\, .
\end{align*}
To estimate $I_1$ we argue as follows: First, since $z\in B(x_0,r)^c \subset B(x_0,2b_2 r)^c$, it follows 
from Proposition \ref{p:estimate-of-J-away} (with $\eps_0=1$)
that for $y\in B(x_0,b_2 r)$, $J^{Y^D}(x_0,z)\asymp J^{Y^D}(y,z)$. Hence, by \eqref{e:exit-time-est},  Lemma \ref{l:pbf}(a) and Lemma \ref{l:tauB},
\begin{eqnarray*}
I_1&\le & c_2 J^{Y^D}(x_0,z)\int_{B(x_0,b_2 r)}G^{Y^D}_{B(x_0,r)}(x_1,y)\, dy \le c_2 J^{Y^D}(x_0,z) \E_{x_1}\tau^{Y^D}_{B(x_0,r)}\\
&\le &c_3 J^{Y^D}(x_0,z)\phi(r^{-2})^{-1}\le c_4 J^{Y^D}(x_0,z)\phi((b_2r)^{-2})^{-1}\\
&\le & c_5 J^{Y^D}(x_0,z) \E_{x_2}\tau^{Y^D}_{B(x_0, b_2 r)}=c_5 J^{Y^D}(x_0,z)\int_{B(x_0,b_2 r)}G^{Y^D}_{B(x_0,b_2 r)}(x_2,y)\, dy\\
&\le & c_5 J^{Y^D}(x_0,z)\int_{B(x_0,b_2 r)}G^{Y^D}_{B(x_0,r)}(x_2,y)\, dy\\
&\le & c_6 \int_{B(x_0,b_2 r)}G^{Y^D}_{B(x_0,r)}(x_2,y)J^{Y^D}(y,z)\, dy \le c_6 K^{D,B(x_0,r)}(x_2,z)\, .
\end{eqnarray*}
Together with the previous display, this proves \eqref{e:pk-est}. Now, let $f$ be a 
non-negative function in $D$ which is harmonic with respect to $Y^D$ in $B(x_0,r)$. 
Then, by the L\'evy system formula and \eqref{e:exit-time-YD},
$$
f(x)=\int_{\overline{B(x_0,r)}^c}K^{D,B(x_0,r)}(x,z)f(z)\, dz\, ,\quad x\in B(x_0,r)\, .
$$
Hence for $x_1,x_2 \in B(x_0,b_1 r)$ we have
$$
f(x_1)=\int_{\overline{B(x_0,r)}^c}K^{D,B(x_0,r)}(x_1,z)f(z)\, dz\le c_6 \int_{\overline{B(x_0,r)}^c}K^{D,B(x_0,r)}(x_2,z)f(z)\, dz=c_6 f(x_2)\, .
$$
For $x_1, x_2\in B(x_0,r/2)$, the inequality follows by a standard chain argument.

Now we assume that $d\le 2$. We will prove this case by reducing it to the high dimensional case. Let $\widetilde D=D\times \R^2$ and let $\widetilde{Y}^{\widetilde{D}}$ be the subordinate killed Brownian motion in $\widetilde D\subset
\R^{d+2}$. For any function $h$ defined in $D$, we define a function
$$
\widetilde h(x, y)=h(x), \qquad x\in D, y\in \R^2.
$$
One can easily check by using the definitions that if $h$ is harmonic in $U\subset D$ with
respect to $Y^D$, then $\widetilde h$ is harmonic in $U\times \R^2$ with respect 
to $\widetilde{Y}^{\widetilde{D}}$. Now the theorem for $d\le 2$ follows from the case for $d\ge 3$.
\qed


\section{Jumping kernel and Green function estimates}\label{s:jkgfe}
 
In this section we will always assume that $d \ge 2$ and that {\bf (A1)}--{\bf (A3)} hold.  If $d=2$, we assume that 
 {\bf (A4)} and {\bf (A5)} also hold.
Throughout this section, 
$D\subset \R^d$ is either a bounded $C^{1,1}$ domain, or
a $C^{1,1}$ domain with compact complement or a
domain consisting of all the points above the graph of a bounded globally 
$C^{1,1}$ function. 
As already mentioned in the introduction, unbounded sets are considered only in 
dimension $d\ge 3$, and then we assume that {\bf (A6)} holds.
We will use $(R, \Lambda)$ to denote the $C^{1, 1}$ characteristics of $D$
in all three cases.

We first recall the following result from \cite{KSV16}. 
\begin{prop}[{\cite[Theorem 3.1]{KSV16}}]\label{p:gfe0}
Suppose that $d \ge 2$.
For every
$M>0$,
there exists a constant $C=C(M)\ge 1$ such that
for all $x,y \in D$ with $|x-y|\le M$,
\begin{eqnarray*}
&&C^{-1}\left(\frac{\delta_D(x)}{|x-y|}\wedge 1
\right)\left(\frac{\delta_D(y)}{|x-y|}\wedge 1
\right)
\frac{\phi'(|x-y|^{-2})}
{|x-y|^{d+2}\phi(|x-y|^{-2})^2}\nonumber\\
&& \le G^{Y^D}(x,y)
\le C \left(\frac{\delta_D(x)}{|x-y|}\wedge 1
\right)\left(\frac{\delta_D(y)}{|x-y|}\wedge 1
\right)\frac{\phi'(|x-y|^{-2})}
{|x-y|^{d+2}\phi(|x-y|^{-2})^2}\, .
\end{eqnarray*}
\end{prop}

We choose a
$C^{1,1}$-function $\varphi: \bR^{d-1}\to \bR$ satisfying $\varphi
(\wt 0)= 0$, $\nabla\varphi (\wt 0)=(0, \dots, 0)$, $\| \nabla \varphi
 \|_\infty \leq \Lambda$, $| \nabla \varphi (\wt y)-\nabla \varphi (\wt w)|
\leq \Lambda |\wt y-\wt w|$, and an orthonormal coordinate system $CS_z$ with
its origin at $z \in \partial D$ such that
$$
B(z, R)\cap D=\{ y=(\wt y, y_d) \in B(0, R) \mbox{ in } CS_z: y_d >
\varphi (\wt y) \}.
$$

Define $\rho_z (x) := x_d -  \varphi (\wt x),$ 
where $(\wt x, x_d)$ are the coordinates of $x$ in $CS_z$. Note that for every $z
\in \partial D$ and $ x \in B(z, R)\cap D$, we have
\begin{equation}\label{e:d_com}
(1+\Lambda^2)^{-1/2} \,\rho_z (x) \,\le\, \delta_D(x)  \,\le\,
\rho_z(x).
\end{equation}
We define for $r_1, r_2>0$,
$$
D_z( r_1, r_2) :=\left\{ y\in D: r_1 >\rho_z(y) >0,\, |\wt y | < r_2
\right\}.
$$

Recall that  $\kappa=(1+(1+ \Lambda)^2)^{-1/2}$.
It is well known (see, for instance \cite[Lemma
2.2]{So}) that there exists $L=L(R, \Lambda, d)>0$ such that for
every $z \in\partial D$ and $r \le \kappa R$, one can find a $C^{1,1}$
domain $V_z(r)$ with characteristics $(rR/L, \Lambda L/r)$
such that $D_z( r/2, r/2) \subset V_z(r)  \subset  D_z( r, r) $. 
In this and the following two sections,
given a $C^{1, 1}$ domain $D$, 
$V_z(r)$ always refers to the $C^{1, 1}$ 
domain above.

It is easy to see that for every $z \in\partial D$ and $r \le \kappa R$,
\begin{align}
\label{e:Uzr}
V_z(r) \subset  D_z( r, r) \subset D\cap B(z,  r /\kappa).
\end{align}
In fact, for all $y \in D_z(r,r)$, 
 \bee\label{e:lsd}
|y|^2 = |\wt y|^2+ |y_d|^2 <r^2 +(|y_d- \varphi(\wt y)|+ |\varphi(\wt
y)|)^2 < (1+(1+ \Lambda)^2) r^2.
\eee

By the fact that $r^{-1}V_z(r)$ is a $C^{1,1}$
domain with characteristics $(R/L, \Lambda L)$, we have  
 the following short-time
estimates (cf.~\cite{Zh}):
There exist positive constants
$c_1, c_2$
such that for any $t\in (0,  4\kappa^{-2}]$ and any $x,y\in r^{-1}V_z(r)$,
$$
p^{r^{-1}V_z(r)}(t,x,y)\ge c_1\left(\frac{\delta_{r^{-1}V_z(r)}(x)}{\sqrt{t}} \wedge 1 \right)
\left(\frac{\delta_{r^{-1}V_z(r)}(y)}{\sqrt{t}} \wedge 1 \right)
\, t^{-d/2}\exp\left(-\frac{c_2 |x-y|^2}{t}\right)\, .
$$
Thus by the scaling property  
$p^{V_z(r)}(t,x,y)=r^{-d} p^{r^{-1}V_z(r)}(r^{-2}t,r^{-1}x,r^{-1}y)$,
we have for any $t\in (0, 4\kappa^{-2} r^2]$ and any $x,y\in V_z(r)$,
\begin{equation}\label{e:lower-bound-for-pr}
p^{V_z(r)}(t,x,y)\ge c_1\left(\frac{\delta_{V_z(r)}(x)}{\sqrt{t}} \wedge 1 \right)
\left(\frac{\delta_{V_z(r)}(y)}{\sqrt{t}} \wedge 1 \right)
\, t^{-d/2}\exp\left(-\frac{c_2 |x-y|^2}{t}\right)\, .
\end{equation}

We will use the following bound 
below: 
By the change of variables
$s=c|x-y|^2/t$, for every $c>0$ and $a  \in \R$, and any open set $U$, we have
\begin{align}\label{e:UD-lower-2}
&\int_0^{|x-y|^2}
\left(\frac{\delta_U(x)}{\sqrt{t}} \wedge 1 \right)
\left(\frac{\delta_U(y)}{\sqrt{t}} \wedge 1 \right)
t^{-a/2}\exp\left(
-\frac{c |x-y|^2}{t}\right) \, dt \nonumber\\
&= \int_{ c}^{\infty}
\left(\frac{\sqrt{s/c}\, \delta_U(x)}{|x-y|} \wedge 1 \right)
\left(\frac{\sqrt{s/c}\, \delta_U(y)}{|x-y|} \wedge 1 \right)
\left(\frac{c|x-y|^2}{s}\right)^{-a/2}e^{-s}
\frac{c|x-y|^2}{s^2}\, ds\nonumber\\
&\ge c^{1-(a/2)} \left(\frac{\delta_U(x)}{|x-y|} \wedge 1 \right)
\left(\frac{\delta_U(y)}{|x-y|} \wedge 1 \right)
|x-y|^{-a+2} \int_{ c}^{\infty} s^{a/2-2}
e^{-s}\, ds.
\end{align}
\begin{prop}\label{p:gfe}
There exists a constant 
$C=C(R, \Lambda)\in (0, 1]$ 
such that
for all $z \in \partial D$, $r \le \kappa R$ and  $x,y \in   V_z(r)$,
\begin{eqnarray}
G^{Y^{V_z(r)}} (x,y)
\ge 
C \left(\frac{\delta_{V_z(r)}(x)}{|x-y|}\wedge 1
\right)\left(\frac{\delta_{V_z(r)}(y)}{|x-y|}\wedge 1
\right)\frac{\phi'(|x-y|^{-2})}
{|x-y|^{d+2}\phi(|x-y|^{-2})^2}\, .\label{e:gfe}
\end{eqnarray}
\end{prop}
\pf
Note that, by \eqref{e:Uzr}, we see that for $x, y\in V_z(r)$, $|x-y| < 2\kappa^{-1}r$. Thus, by \eqref{e:lower-bound-for-pr},
$$
G^{Y^{V_z(r)}}(x,y) \ge c_{1}\int_0^{|x-y|^2}
\left(\frac{\delta_{V_z(r)}(x)}{\sqrt{t}} \wedge 1 \right)
\left(\frac{\delta_{V_z(r)}(y)}{\sqrt{t}} \wedge 1 \right)
 t^{-d/2}\exp\left(
-\frac{c_2 |x-y|^2}{t}\right) u(t)\, dt. 
$$
Using the fact that $u$ is decreasing and \eqref{e:lower-estimate-u}, we now have 
\begin{align}\label{e:UD-lower-1}
&
G^{Y^{V_z(r)}}(x,y)\nonumber\\
&\ge c_{1}u(|x-y|^2)\int_0^{|x-y|^2}
\left(\frac{\delta_{V_z(r)}(x)}{\sqrt{t}} \wedge 1 \right)
\left(\frac{\delta_{V_z(r)}(y)}{\sqrt{t}} \wedge 1 \right)
t^{-d/2}\exp\left(
-\frac{c_2 |x-y|^2}{t}\right) \, dt\nonumber\\
&\ge c_{3}\frac{\phi'(|x-y|^{-2})}
{|x-y|^{4}\phi(|x-y|^{-2})^2} \int_0^{|x-y|^2}\left(\frac{\delta_{V_z(r)}(x)}{\sqrt{t}} \wedge 1 \right)\nonumber\\
&\qquad\times
\left(\frac{\delta_{V_z(r)}(y)}{\sqrt{t}} \wedge 1 \right)
t^{-d/2}\exp\left(
-\frac{c_2 |x-y|^2}{t}\right) \, dt.
\end{align}By combining this with \eqref{e:UD-lower-2}
we arrive at
\begin{align*}
G^{Y^{V_z(r)}} (x,y)&\ge c_{4}
\left(\frac{\delta_{V_z(r)}(x)}{|x-y|} \wedge 1 \right)
\left(\frac{\delta_{V_z(r)}(y)}{|x-y|} \wedge 1 \right)
\frac{\phi'(|x-y|^{-2})}
{|x-y|^{d+2}\phi(|x-y|^{-2})^2}.
\end{align*}
\qed

We now  consider estimates for $J^{Y^D}$. We first recall
\begin{prop}[{\cite[Proposition 3.5]{KSV16}}]
\label{p:J^{Y^D}(z,y)}
For every $M>0$,
there exists a constant $
C=C(M, R, \Lambda)
\ge 1$ such that such that for all $x,y \in D$ with $|x-y|\le M$,
\begin{align*}
&C^{-1}\left(\frac{\delta_D(x)}{|x-y|}\wedge 1
\right)
\left(\frac{\delta_D(y)}{|x-y|}\wedge 1
\right)
\frac{\phi'(|x-y|^{-2})}
{|x-y|^{d+2}}
\\
& \le J^{Y^D}(x,y)
\le C \left(\frac{\delta_D(x)}{|x-y|}\wedge 1
\right)
\left(\frac{\delta_D(y)}{|x-y|}\wedge 1
\right)
\frac{\phi'(|x-y|^{-2})}
{|x-y|^{d+2}}\, .
\end{align*}
\end{prop}

Propositions \ref{p:gfe0} and \ref{p:J^{Y^D}(z,y)} imply global two-sided estimates
on $G^{Y^D}$ and $J^{Y^D}$ for bounded $D$, but only give ``local'' two-sided estimates for unbounded $D$.
Now we assume $d \ge 3$ and give global two-sided estimates for $J^{Y^D}$ for our
two types of unbounded $C^{1,1}$ domains.
The proof of the next theorem is very similar to that of \cite[Theorem 3.2]{KSV16}, where two-sided global estimates were proved for $G^{Y^D}$.
\begin{thm}\label{t:Jfe:nn}
Suppose that $d \ge 3$ and that  {\bf (A1)}--{\bf (A3)} and {\bf (A6)} hold true.

\noindent
(1) Let $D\subset \R^d$ be a
domain consisting of all the points above the graph of a bounded globally $C^{1,1}$
function.
There exists a constant $
C_1=C_1(R, \Lambda) \ge 1$ such that for all $x,y \in D$,
\begin{align*}
&C_1^{-1}\left(\frac{\delta_D(x)}{|x-y|} \wedge 1 \right)
\left(\frac{\delta_D(y)}{|x-y|} \wedge 1 \right)\frac{\mu(|x-y|^2)}
{|x-y|^{d-2}} \le J^{Y^D}(x,y)\nn\\
&\qquad \le C_1 \left(\frac{\delta_D(x)}{|x-y|} \wedge 1 \right)
\left(\frac{\delta_D(y)}{|x-y|} \wedge 1 \right)\frac{\mu(|x-y|^2)}
{|x-y|^{d-2}}\, .
\end{align*}
(2) Let $D\subset \R^d$ be a
$C^{1,1}$ domain with compact complement. There exists a constant $
C_2=C_2(R, \Lambda)\ge 1$ such that
for all $x,y \in D$,
\begin{align*}
&C_2^{-1}\left(\frac{\delta_D(x)}{|x-y|\wedge 1} \wedge 1 \right)
\left(\frac{\delta_D(y)}{|x-y|\wedge 1} \wedge 1 \right)\frac{\mu(|x-y|^2)}
{|x-y|^{d-2}} \le J^{Y^D}(x,y) \nn\\
&\qquad \le C_2 \left(\frac{\delta_D(x)}{|x-y|\wedge 1} \wedge 1 \right)
\left(\frac{\delta_D(y)}{|x-y|\wedge 1} \wedge 1 \right)\frac{\mu(|x-y|^2)}
{|x-y|^{d-2}}\, .
\end{align*}
\end{thm}

\pf
Using 
\cite[(3.4)]{KSV16} and \cite[(3.10)--(3.13)]{KSV16} instead of \cite[(3.14)--(3.19)]{KSV16},
the proof of (1) is similar to (2).
 Thus, we give the proof of (2) only.

\noindent
{\it Upper bound}:
Using \eqref{e:mas} and the fact  $\mu$ is decreasing, we have
from \cite[(3.14)]{KSV16} that
\begin{align*}
&J^{Y^D}(x,y)\,=\,\int_0^{\infty}p^D(t,x,y)\mu(t)\, dt\\
&\le c_1
\int_0^{\infty}\left(\frac{\delta_D(x)}{\sqrt{t}\wedge 1} \wedge 1 \right)
\left(\frac{\delta_D(y)}{\sqrt{t}\wedge 1} \wedge 1 \right)t^{-d/2}\exp
\left(-\frac{c_2|x-y|^2}{t}\right)\mu(t)\, dt
\nonumber \\
&\le c_3 |x-y|^{2\beta} \mu(|x-y|^2)
\int_0^{|x-y|^2}\left(\frac{\delta_D(x)}{\sqrt{t}\wedge 1} \wedge 1 \right)
\left(\frac{\delta_D(y)}{\sqrt{t}\wedge 1} \wedge 1 \right)t^{-\beta-d/2}\exp
\left(-c_2\frac{|x-y|^2}{t}\right)\, dt \\
&\quad
+c_1 \mu(|x-y|^2)
\int_{|x-y|^2}^{\infty} \left(\frac{\delta_D(x)}{\sqrt{t}\wedge 1} \wedge 1 \right)
\left(\frac{\delta_D(y)}{\sqrt{t}\wedge 1} \wedge 1 \right)t^{-d/2}\, dt.
\end{align*}
Together with  \cite[(3.16)--(3.17)]{KSV16} we obtain the upper bound.

\noindent
{\it Lower bound}:
Since $\mu$ is decreasing,  by
\cite[(3.15)]{KSV16},
\begin{align*}
&J^{Y^D}(x,y)\ge c_4\int_0^{|x-y|^2}
\left(\frac{\delta_D(x)}{\sqrt{t}\wedge 1} \wedge 1 \right)
\left(\frac{\delta_D(y)}{\sqrt{t}\wedge 1} \wedge 1 \right)
 t^{-d/2}\exp\left(
-\frac{c_5 |x-y|^2}{t}\right) \mu(t)\, dt \nonumber\\
&\ge c_4\mu(|x-y|^2)\int_0^{|x-y|^2}
\left(\frac{\delta_D(x)}{\sqrt{t}\wedge 1} \wedge 1 \right)
\left(\frac{\delta_D(y)}{\sqrt{t}\wedge 1} \wedge 1 \right)
t^{-d/2}\exp\left(
-\frac{c_5 |x-y|^2}{t}\right) \, dt\, .
\end{align*}
Combining this and \cite[(3.18)]{KSV16}
we arrive at
\begin{align*}
J^{Y^D}(x,y)&\ge c_6
\left(\frac{\delta_D(x)}{|x-y|\wedge 1} \wedge 1 \right)
\left(\frac{\delta_D(y)}{|x-y|\wedge 1} \wedge 1 \right)
\frac{\mu(|x-y|^2)}
{|x-y|^{d-2}}\, . \end{align*}
\qed


\section{Carleson estimate in $C^{1,1}$ domain}\label{s:ce}

Before we discuss the Carleson estimate in $C^{1,1}$ domains, 
we first present two preliminary results. 
Recall that $\zeta$ is the lifetime of $Y^D$.
\begin{lemma}\label{l:regularity}
Suppose that $D$ is an open set in $\R^d$. Let $x_0\in \R^d$, and $r_1<r_2$ be two positive numbers such that $ D \cap B(x_0, r_1) \neq \emptyset$.
Suppose $f$ is a non-negative function in $D$ that is
harmonic in $D\cap B(x_0, r_2)$ with respect to $Y^D$ and vanishes
continuously on $\partial D\cap B(x_0, r_2)$. Then $f$ is regular harmonic in
$D\cap B(x_0, r_1)$ with respect to $Y^D$, i.e.,
\begin{equation}\label{e:regularity}
f(x)=\E_x\left[ f(Y^D_{\tau_{D\cap B(x_0, r_1)}})\right], \qquad \hbox{
for all }x\in D\cap B(x_0, r_1)\, .
\end{equation}
\end{lemma}

\pf 
For $n\ge 1$, let  $B_n=\{y\in D\cap B(x_0, r_1):\, \delta_D(y)>1/n\} $.
Then for large $n$, $B_n$ is a non-empty open subset of $D\cap B(x_0, r_1)$
whose closure is contained in $D\cap B(x_0, r_2)$. Since $f$ is harmonic in
$D\cap B(x_0, r_2)$ with respect to $Y^D$, for $x\in D\cap B(x_0, r_1)$ and $n$ large
enough so that $x\in B_n$, we have
\begin{align*}
&f(x)=\E_x\left[f\big(Y^D_{\tau_{B_n}}\big)\right]\\
&=\E_x\left[
f\big(Y^D_{\tau_{B_n}}\big); \, \tau_{B_n}< \tau_{D\cap B(x_0, r_1)}\right]
+\E_x\left[f\big(Y^D_{\tau_{B_n}}\big);\, \tau_{B_n}=\tau_{D\cap B(x_0, r_1)}\right]\, .
\end{align*}
Hence
\begin{align}
&\left|f(x)-\E_x\left[f\big(Y^D_{\tau_{D \cap B(x_0, r_1)}}\big)
\right] \right|\nonumber\\
&\le \E_x\left[f\big(Y^D_{\tau_{B_n}}\big);\, \tau_{B_n}<
\tau_{D\cap B(x_0, r_1)}\right] +\E_x\left[f\big(Y^D_{\tau_{D\cap
B(x_0, r_1)}}\big);\, 
\tau_{B_n}< \tau_{D\cap B(x_0, r_1)}\right]\, .
\label{e:reg2}
\end{align}
It follows from \eqref{e:YD-killed-inside}   that 
$\cap_{n=1}^{\infty}\{\tau_{B_n} < \tau_{D\cap B(x_0, r_1)}\} =\emptyset$
 almost surely under each $\P_x$. 
 Thus,
$$
\lim_{n \to \infty} \E_x\left[f\big(Y^D_{\tau_{D\cap B(x_0, r_1)}}\big);\,
\tau_{B_n}< \tau_{D\cap B(x_0, r_1)}\right]=0\, .
$$
For the first term in \eqref{e:reg2}, note that
$\delta_D(Y^D_{\tau_{B_n}}) \le 1/n$ on 
$\{\tau_{B_n}<\tau_{D\cap B(x_0, r_1)}\}$. 
Since $f$ vanishes continuously on $(\partial D)\cap B(x_0, r_2)$, an easy compactness argument yields that there is $n_0\geq 1$ so
that $f$ is bounded in $(D\cap B(x_0, r_1))\setminus B_{n_0}$.
Hence by the
bounded convergence theorem we have
$$
\lim_{n \to  \infty}
\E_x\left[
f\big(Y^D_{\tau_{B_n}}\big);\, \tau_{B_n}< \tau_{D\cap
B(x_0, r_1)}\right]=0\,.
$$
This proves the lemma.
\qed
\begin{lemma}\label{l:detlax}
 For any $x,y \in \R^d$ and any open set $V$, 
 we have 
 $$ 
 \delta_V(y)\left(\frac{\delta_V(x)}{|x-y|}\wedge 1 \right) \le 2 \delta_V(x).
 $$
 \end{lemma}
 \pf 
 If $|x-y| \le \delta_V(x)$,
 $$
  \delta_V(y)\left(\frac{\delta_V(x)}{|x-y|}\wedge 1 \right) =
  \delta_V(y) \le   \delta_V(x)+|x-y| \le  2 \delta_V(x). 
$$
 If $|x-y| > \delta_V(x)$,
 $$
   \delta_V(y)
  \left(\frac{\delta_V(x)}{|x-y|}\wedge 1 \right) =
  \delta_V(x) \frac{\delta_V(y)}{|x-y|}\le 
   \delta_V(x) \frac{ \delta_V(x)+|x-y|}{|x-y|}\le   2 \delta_V(x). 
$$
\qed

In  the remainder of this section 
we will assume that $d \ge 2$ and that {\bf (A1)}--{\bf (A3)} hold.  If $d=2$, we assume that 
 {\bf (A4)} and {\bf (A5)} also hold. We assume that 
$D\subset \R^d$ is either a bounded $C^{1,1}$ domain, or
a $C^{1,1}$ domain with compact complement or a
domain consisting of all the points above the graph of a bounded globally 
$C^{1,1}$ function. 
As before, unbounded sets are considered only in dimensions $d\ge 3$, and then we assume that {\bf (A6)} holds. 
Let $(R, \Lambda)$ be the $C^{1, 1}$ characteristics
of $D$.
Without loss of generality we assume that $R \le 1$.
Recall $\kappa=(1+(1+ \Lambda)^2)^{-1/2}$,  $\rho_z (x)= x_d -  \varphi_z (\wt x)$ and $V_z(r)$ is  a $C^{1,1}$
domain  with characteristic $(rR/L, \Lambda L/r)$
such that $D_z( r/2, r/2) \subset V_z(r)  \subset  D_z( r, r) $ where  $L=L(R, \Lambda, d)>0$.
\begin{lemma}\label{lower bound}   There exists
a constant
$\delta_*=\delta_*(R, \Lambda)>0$
 such that for all
$Q \in \partial D$ and  $x\in D$ with $\rho_Q(x) < R/2$,
$$
\P_x\left(  \tau(x)=\zeta     \right)\ge \delta_*\, ,
$$
where $\tau(x):=\tau^{Y^D}_{D\cap B(x,
\rho_Q(x)/2)}=\inf\{t>0:\,
Y^D_t\notin D\cap B(x,
\rho_Q(x)/2)\}$.
\end{lemma}

\pf
By \cite[Theorem 4.5.4(1)]{FOT},
$$
\P_x\left(  \tau(x)=\zeta     \right)=\P_x(Y^D_{\zeta-}\in D\cap B(x,
\rho_Q(x)/2)
= \int_{D\cap B(x,
\rho_Q(x)/2)} G^{Y^D}(x,y) \kappa^{Y^D} (y)dy,
$$
where  $\kappa^{Y^D}$ 
is the density of the killing measure of $Y^D$ given by
$$
\kappa^{Y^D}(x)=
\int_0^\infty (1-P^D_t 1(x))\, \mu(t)dt\, ,\quad x\in D\, .
$$
Since $D$ is a $C^{1,1}$ domain, 
we have (see the proof of \cite[Lemma 5.7]{KSV16}) 
\begin{equation}\label{e:estimate-kappa-D}
\kappa^{Y^D} (y)\ge \kappa^{X^D}(y) \ge c_1 \phi (\delta_D(y)^{-2}), \qquad y \in D\cap B(x, \rho_Q(x)/2).
\end{equation}
Here \eqref{e:killing-functions-relation} is used in the first inequality. 
Thus, using \eqref{e:estimate-kappa-D} and Proposition \ref{p:gfe0},
\begin{align*}
&\P_x\left(  \tau(x)=\zeta     \right) = \int_{D\cap B(x,
\rho_Q(x)/2)} G^{Y^D}(x,y) \kappa^{Y^D} (y)dy\\
& \ge c_1 \int_{D\cap B(x,
\rho_Q(x)/2)} 
 \left(\frac{\delta_D(x)}{|x-y|}\wedge 1
\right)\left(\frac{\delta_D(y)}{|x-y|}\wedge 1
\right)\frac{\phi'(|x-y|^{-2}) \phi (\delta_D(y)^{-2})}
{|x-y|^{d+2}\phi(|x-y|^{-2})^2}
 dy\\
  & \ge c_2 \int_{ B(x,\rho_Q(x)/2)} 
\frac{\phi'(|x-y|^{-2})  \phi (\delta_D(y)^{-2})}
{|x-y|^{d+2}\phi(|x-y|^{-2})^2}
dy\\
  & \ge c_3 \phi (\rho_Q(x)^{-2})  \int_{ B(x,\rho_Q(x)/2)} 
\frac{\phi'(|x-y|^{-2})}
{|x-y|^{d+2}\phi(|x-y|^{-2})^2}
dy \\
&\ge c_4 \phi (\rho_Q(x)^{-2}) \int_{0}^{\rho_Q(x)/2}  (1/\phi(t ^{-2}))'dt \ge c_5.
\end{align*}
 \qed

For $a \ge 0$, let
\begin{align}  \label{e:w_1}
a_1=\begin{cases}
a & \text{ if $D$ is either a bounded $C^{1,1}$ domain, or a domain consisting}\\
& \text{\ of all the points above the graph of a bounded globally $C^{1,1}$ function},\\
a \wedge 1 &\text{ if  $D$ is a $C^{1,1}$
 domain with compact complement.}
\end{cases}
\end{align}
Note that for every $a \in (0,1]$ and $b>0$, we have 
$$
\frac{a}{(ab) \wedge 1} \le \frac{a}{(ab) \wedge a} \le   \frac{1}{b \wedge 1}.
$$
Thus for every $a \in (0,1]$ and $b, c>0$, 
\begin{align}
\label{e:aba+}
\frac{1}{(ab)_1} \le a^{-1} \frac{1}{b_1} \quad \text{ and } \quad  
\frac{c}{(ab)_1} \wedge 1 \le a^{-1} \frac{c}{b_1} \wedge 1 \le a^{-1} (\frac{c}{b_1} \wedge 1).
\end{align}
\begin{thm}[Carleson estimate]\label{carleson}
There exists a constant
$C=C(R, \Lambda)>0$ such that for every $Q\in
\partial D$, $0<r<R/2$, and every non-negative function
$f$ in $D$ that is harmonic in $D \cap B(Q, r)$ with respect to
$Y^D$ and vanishes continuously on $ \partial D \cap B(Q, r)$, we have
\begin{equation}\label{e:carleson}
f(x)\le C f(x_0) \qquad \hbox{for }  x\in D\cap B(Q,r/2),
\end{equation}
where $x_0\in D
\cap B(Q,r)$ with $\rho_Q(x_0)=r/2$.
\end{thm}

\pf
In this proof, the constants $\delta_*, \nu, \gamma, \beta_1, \eta$ and $c_i$'s are always
independent of $r$. Without loss of the generality,
we assume that diam($D$) $\le 1$ if $D$ is bounded.
Using the assumption that $D$ is $C^{1,1}$ and $r<R/2$, by
Theorem \ref{uhp}
and a standard chain argument, it suffices to prove (\ref{e:carleson}) for 
$x\in D\cap B(Q, \kappa r/(24))$
and $\wt x_0 = \wt 0$ in $CS_Q$.

Let 
$$
k(s)=\frac{\mu(s^2)}{s^{d-2}}.
$$ 
Then $k$ is decreasing and there exists $c>1$ such that
\begin{eqnarray}
k(s)&\le& c\, k(2s), \quad \forall s\in (0, 3)  \text{ if  $D$ is bounded, }\label{H:1}\\
k(s)&\le& c\, k(2s), \quad \forall  s\in (0, \infty)  \text{ if  $D$ is unbounded.}
\label{H:1a}
\end{eqnarray}
In fact, by \eqref{e:upper-estimate-mu} and 
\eqref{e:lower-estimate-mu}, we have 
\begin{align*}
&k(2s) =2^{-d+2} \mu(4s^2)
s^{-d+2}  \ge c_0 2^{-d} \phi'(4^{-1}s^{-2}) s^{-d-2}\\
&\ge c_0 2^{-d} \phi'(s^{-2}) s^{-d-2} 
\ge c_0 2^{-d}  (1-2e^{-1})   k(s).
  \end{align*}
If $D$ is unbounded,
then \eqref{H:1a} follows from  {\bf (A6)}.

Note that, as a consequence of \eqref{H:1}, 
there is a $\nu > 2$ such that 
\begin{equation}\label{H:1n}
k(as) \le c_1\, a^{-\nu+2} k(s), \qquad \forall s\in (0, 3)
\quad\text{and}\quad  a \in (0, 1).
\end{equation}
Choose $0<\gamma < \nu^{-1}$. For any $x\in D\cap B(Q,\kappa r/(12))$, define
$$
D_0(x)=D\cap B(x,2\rho_Q(x))\, ,\qquad B_1(x)=B(x,r^{1-\gamma}
\rho_Q(x)^{\gamma})\,
$$
and
$$
B_2=B(x_0,\kappa\rho_Q(x_0)/3)\, ,\qquad 
B_3=B(x_0, 2\kappa\rho_Q(x_0)/3).
$$
Since $x\in B(Q,\kappa r/(12))$, we have $\rho_Q(x)<r/(12)$. By 
the choice of $\gamma<1/2$, 
we have that $D_0(x)\subset B_1(x)$.
By Lemma \ref{lower bound}, there exists $\delta_*=\delta_*(R,
\Lambda)>0$ such that
\begin{equation}\label{e:c:1}
\P_x(\tau^{Y^D}_{D_0(x)}
=\zeta)
\ge\P_x(\tau^{Y^D}_{D\cap B(x,\rho_Q(x)/2)}
=\zeta) \ge \delta_*\, ,
\quad x\in D\cap B(Q,\kappa r/(12))\, .
\end{equation}
By Theorem \ref{uhp} and a chain argument, there exists
$\beta_1>0$ such that
\begin{equation}\label{e:c:2}
f(x)<(\rho_Q(x)/r)^{-\beta_1} f(x_0)\, ,
\quad x\in D\cap B(Q,\kappa r/(12)))\, .
\end{equation}
Since $f$ is regular harmonic in
$D_0(x)$ with respect to $Y^D$, by Lemma \ref{l:regularity}, 
for every $x\in D\cap B(Q,\kappa r/(12)))$,
\begin{align}
&f(x)=\E_x\big[f\big(Y^D(\tau_{D_0(x)})\big); Y^D(\tau_{D_0(x)})\in
B_1(x)\big]\nonumber\\
&\qquad+ \E_x\big[f\big(Y^D(\tau_{D_0(x)})\big);
Y^D(\tau_{D_0(x)})\notin B_1(x)\big]. \label{e:c:3}
\end{align}
We first show that there exists $\eta>0$ such that
for all $x\in D \cap B(Q, \kappa r/(12)) $ 
with $\rho_Q(x) < \eta r$,
\begin{equation}\label{e:c:4}
\E_x\big[f\big(Y^D(\tau_{D_0(x)})\big); Y^D(\tau_{D_0(x)})\notin
B_1(x)\big]\le f(x_0). 
\end{equation}
Let $\eta_0 :=2^{-2 \nu }$,
then,
since $\gamma < 1-\nu^{-1}$ (so $\nu > (1-\gamma)^{-1}$),  we have $\eta_0 < 4^{-(1-\gamma)^{-1}}$.
Thus
for $\rho_Q(x)< \eta_0 r$,
$$
2\rho_Q(x) \le r^{1-\gamma} \rho_Q(x)^{\gamma} - 2\rho_Q(x).
$$
Thus if $x\in D \cap B(Q, \kappa r/(12))$
with $\rho_Q(x) < \eta_0r$,  then
$|x-y|\le 2|z-y|$ for $z\in D_0(x)$, $y\notin B_1(x)$.  
Thus  by
\eqref{H:1}, \eqref{H:1a},
 Proposition \ref{p:J^{Y^D}(z,y)} and Theorem  \ref{t:Jfe:nn}, we have that  
for $y\in D\setminus B_1(x)$ and  $z \in D_0(x)$,
 \begin{align}
\label{e:Jlow1}
J^{Y^D}(z,y)dy & \asymp \left(\frac{\delta_D(z)}{|z-y|_1} \wedge 1 \right)
\left(\frac{\delta_D(y)}{|z-y|_1} \wedge 1 \right)k(|z-y|)\nn\\
& \le c_2 \left(\frac{\delta_D(z)}{|x-y|_1} \wedge 1 \right)
\left(\frac{\delta_D(y)}{|x-y|_1} \wedge 1 \right)k(|x-y|)\nn\\
& \le c_2 \frac{k(|x-y|)  }{|x-y|_1}
\left(\frac{\delta_D(y)}{|x-y|_1} \wedge 1 \right) .
\end{align}
By Proposition \ref{p:gfe0},
$$
G^{Y^D}(x,z)
\le c_3 \left(\frac{\delta_D(x)}{|x-z|}\wedge 1
\right)\left(\frac{\delta_D(z)}{|x-z|}\wedge 1
\right)\frac{\phi'(|x-z|^{-2})}
{|x-z|^{d+2}\phi(|x-z|^{-2})^2}, \quad z \in  D_0(x).
$$
Thus, using Lemma \ref{l:detlax}, we get
\begin{align}
\label{e:D3_1}
&\E_{x}
\int_0^{\tau_{ D_0(x)}}\delta_D(Y^D_t) dt 
=\int_{D_0(x)} G^{Y^D}_{D_0(x)}(x,z) \delta_D(z) dz\nn\\
&\le   \int_{D_0(x)} G^{Y^D}(x,z) \delta_D(z) dz\nn\\
&\le c_3  \int_{D_0(x)} \left(\frac{\delta_D(x)}{|x-z|}\wedge 1
\right)\frac{\phi'(|x-z|^{-2})}
{|x-z|^{d+2}\phi(|x-z|^{-2})^2} \delta_D(z) dz\nn\\
&\le 2c_3  \rho_Q(x) \int_{B(x,2\rho_Q(x))}\frac{\phi'(|x-z|^{-2})}
{|x-z|^{d+2}\phi(|x-z|^{-2})^2}  dz \le  c_4 \frac{ \rho_Q(x)}{ \phi(\rho_Q(x)^{-2})}.
\end{align}
Therefore, by \eqref{e:Jlow1}, \eqref{e:D3_1} and the fact that $D_0(x)\subset B_1(x)$,
\begin{align}\label{e:c:5}
&\E_x\big[f\big(Y^D(\tau_{D_0(x)})\big); Y^D(\tau_{D_0(x)}) \notin
B_1(x)\big]\nonumber\\
&=\E_x \int_0^{\tau_{D_0(x)}}   
\int_{D\setminus B_1(x)}
J^{Y^D}(Y^D_t,y)f(y)\, dy\, dt
\nonumber \\
&\le c_5 \E_x [\int_0^{\tau_{D_0(x)}}  \delta_D (Y^D_t)dt ]      
\int_{D\setminus B_1(x)}
    \frac{f(y) }{|x-y|_1}
\left(\frac{\delta_D(y)}{|x-y|_1} \wedge 1 \right)     k(|x-y|)   dy
\nonumber \\
&\le c_6 \frac{\rho_Q(x)}{ \phi(\rho_Q(x)^{-2})}
\left(
\int_{(D\setminus B_1(x))\cap B_3^c }
    \frac{f(y) }{|x-y|_1}
\left(\frac{\delta_D(y)}{|x-y|_1} \wedge 1 \right)     k(|x-y|)  dy \right.\nonumber \\
&\quad \quad \quad +\left.
\int_{(D\setminus B_1(x))\cap B_3}
 \frac{f(y) }{|x-y|_1}
\left(\frac{\delta_D(y)}{|x-y|_1} \wedge 1 \right)   k(|x-y|)   \, dy\right)\nonumber\\
&=:\, c_6 \frac{\rho_Q(x)}{ \phi(\rho_Q(x)^{-2})}
(I_1+I_2)\, .
\end{align}
On the other hand, for $z\in B_2$ and $y\notin B_3$, we have
$|z-y|\le |z-x_0|+|x_0-y|\le \kappa\rho_Q(x_0)/3+|x_0-y|\le 2|x_0-y|$ and
$|z-y|\le |z-x_0|+|x_0-y|\le 1+|x_0-y|$. 
By \eqref{H:1}, \eqref{H:1a}, 
Proposition \ref{p:J^{Y^D}(z,y)} and Theorem  \ref{t:Jfe:nn}, we have that  
for $(z, y) \in B_2 \times (D\setminus B_3) $, 
 \begin{align}
\label{e:Jlow1_1}
 J^{Y^D}(z,y)&\asymp
    \left(\frac{\delta_D(z)}{|z-y|_1} \wedge 1 \right)
\left(\frac{\delta_D(y)}{|z-y|_1} \wedge 1 \right)k(|z-y|) \nn\\
&\ge 
c_7 
    \left(\frac{\delta_D(z)}{|x_0-y|_1} \wedge 1 \right)
\left(\frac{\delta_D(y)}{|x_0-y|_1} \wedge 1 \right)k(|x_0-y|)   \nn\\
&\ge 
c_8 r
  \frac{f(y)}{|x_0-y|_1}     
\left(\frac{\delta_D(y)}{|x_0-y|_1} \wedge 1 \right)
k(|x_0-y|) ].
\end{align}
Thus, by Lemma \ref{l:tauB}
we have 
\begin{align}\label{e:c:6}
&f(x_0)\,\ge\, \E_{x_0}\left[f(Y^D(\tau_{B_2})); Y^D(\tau_{B_2})\notin B_3\right]\nonumber \\
&=    \E_{x_0} \int_0^{\tau_{B_2}}\left(
\int_{D\setminus B_3}
  J^{Y^D}(Y^D_t, y)f(y)\, dy
\right) dt\nonumber \\
&\ge 
c_8 r \E_{x_0} [ {\tau^{Y^D}_{B_2}}  ] 
\left(
\int_{D\setminus B_3}
\frac{f(y)}{|x_0-y|_1}     
\left(\frac{\delta_D(y)}{|x_0-y|_1} \wedge 1 \right) k(|x_0-y|) 
  dy\right)\nonumber \\
&\ge 
c_9\frac{ r}
{\phi( r^{-2})}
\int_{D\setminus B_3}
\frac{f(y)}{|x_0-y|_1}     
\left(\frac{\delta_D(y)}{|x_0-y|_1} \wedge 1 \right) k(|x_0-y|) \, dy\, .
\end{align}

Suppose now that $|y-x|\ge r^{1-\gamma}\rho_Q(x)^{\gamma}$ and 
$x\in B(Q, \kappa r/(12))$. Then
$$
|y-x_0|\le |y-x|+r\le
|y-x|+r^{\gamma}\rho_Q(x)^{-\gamma}|y-x|\le 2r^{\gamma}\rho_Q(x)^{-\gamma}|y-x|.
$$
Thus, using \eqref{e:aba+} and \eqref{H:1n}, we
get  for $|x-y| \le 2$,
\begin{align}\label{e:gf1}  
&
\frac{k(|x-y|)}{|x-y|_1}     
\left(\frac{\delta_D(y)}{|x-y|_1} \wedge 1 \right)\nn\\
&\le
\frac{k(2^{-1}  (\rho_Q(x)/r)^{\gamma} |x_0-y|)}{(2^{-1}  (\rho_Q(x)/r)^{\gamma} |x_0-y|)_1}     
\left(\frac{\delta_D(y)}{(2^{-1}  (\rho_Q(x)/r)^{\gamma} |x_0-y|)_1} \wedge 1 \right)\nn \\
& \le  
c_1 2^\nu (\rho_Q(x)/r)^{-\nu \gamma}  \frac{k(|x_0-y|)}{|x_0-y|_1}     
\left(\frac{\delta_D(y)}{|x_0-y|_1} \wedge 1 \right).
\end{align}
Now, using 
\eqref{H:1} and \eqref{H:1a} 
together with $|y-x_0|\le
|y-x|+1/2 \le 2 |y-x|$ if $|y-x| \ge 2$
and \eqref{e:c:6}--\eqref{e:gf1},
\begin{align}\label{e:c:7}
&I_1\le c_{10} 
\int_{D\cap\{2>|y-x|>r^{1-\gamma}\rho_Q(x)^{\gamma}\}\cap B_3^c }(\rho_Q(x)/r)^{-\nu \gamma}
\frac{k(|x_0-y|)}{|x_0-y|_1}      \left(\frac{\delta_D(y)}{|x_0-y|_1} \wedge 1 \right) f(y) \, dy\nonumber\\
&\quad +c_{10}
\int_{D\cap \{|y-x|\ge 2\} \cap B_3^c   } \frac{k(|x_0-y|)}{|x_0-y|_1}     
\left(\frac{\delta_D(y)}{|x_0-y|_1} \wedge 1 \right)f(y)\, dy\nonumber\\
&\le  c_{11} \left( (\rho_Q(x)/r)^{- \nu \gamma}+1\right)\int_{D\setminus B_3}
\frac{k(|x_0-y|)}{|x_0-y|_1}      \left(\frac{\delta_D(y)}{|x_0-y|_1} \wedge 1 \right)\, f(y)\, dy\nonumber\\
&\le 
c_9^{-1}
 c_{11}\frac{\phi(r^{-2})} {r}    \left( (\rho_Q(x)/r)^{-
\nu \gamma}+1\right)f(x_0)\nonumber\\
&\le 2
c_9^{-1} c_{11} (\rho_Q(x)/r)^{-
\nu \gamma}
\frac{\phi(r^{-2})} {r} f(x_0) \, ,
\end{align}
where the second to last inequality is due to \eqref{e:c:6}.
 
If $y\in B_3(x)$, then $\delta_D(y)\le c_{12}r$ and 
$|y-x|\ge |x_0-Q|-|x-Q|-|y-x_0|>\kappa\rho_Q(x_0)/6$.
By Theorem \ref{uhp}, there exists $c_{13}>0$
such that $f(y)\le c_{13}f(x_0)$ for all $y\in B_3(x)$. 
Thus by \eqref{e:upper-estimate-mu}, 
\begin{align}\label{e:c:8}
&I_2 \le c_{14} f(x_0) r
\int_{(D\setminus B_1(x))\cap B_3}
 \frac{k(|x-y|)  }{(|x-y|_1)^2}       \, dy\nonumber\\
&\le  
  c_{14} f(x_0)r
\int_{|y-x|>\kappa\rho_Q(x_0)/6} 
 \frac{k(|x-y|)  }{(|x-y|_1)^2}       \, dy\nonumber\\
&\le   c_{14} f(x_0)r\left(
\int_{1>|z|>\kappa\rho_Q(x_0)/6}
|z|^{-2} k(|z|)\, dz +
\int_{
1 \le |z|} k(|z|)\, dz \right)\nonumber\\
 &\le   c_{15} f(x_0) r\left(
\int_{1>|z|>\kappa\rho_Q(x_0)/6} 
\frac{\phi'(|z|^{-2})}
{|z|^{d+4}}\, dz +
1 \right)\nn\\
 &\le   c_{16} f(x_0) r\left(
 \int_{\kappa r/(12)}^1
 \frac{\phi'(s^{-2})}
{s^{5}}\, ds +
1 \right)\nn\\
 &\le   c_{16} f(x_0) r\left(
(12)^2 (\kappa r)^{-2}  \int_{\kappa r/(12)}^\infty 
\frac{\phi'(s^{-2})}
{s^{3}}\, ds +
1 \right)\,\le  \,  c_{17} f(x_0) 
r^{-1} \phi(r^{-2}).
\end{align}
Combining \eqref{e:c:5}, \eqref{e:c:7} and \eqref{e:c:8}, we obtain
\begin{align}\label{e:c:9}
&\E_x[f(Y^D(\tau_{D_0(x)}));\, Y^D(\tau_{D_0(x)})\notin B_1(x)]\nonumber\\
&\le c_{18} f(x_0)\Big( \frac{\rho_Q(x)}{ \phi(\rho_Q(x)^{-2})}(\rho_Q(x)/r)^{-
\nu \gamma}
\frac{\phi(r^{-2})} { r}
+
\frac{\rho_Q(x)}{ \phi(\rho_Q(x)^{-2})}  r^{-1}    \phi(r^{-2}) \Big)\nonumber \\
&=c_{18} f(x_0)\Big( \frac{\phi(r^{-2})}{ \phi(\rho_Q(x)^{-2})}(\rho_Q(x)/r)^{1-
\nu \gamma}
 +
\frac{\phi(r^{-2})}{ \phi(\rho_Q(x)^{-2})}(\rho_Q(x)/r) \Big)\nonumber \\
&\le c_{18} f(x_0)\left((\rho_Q(x)/r)^{1-
\gamma
\nu}
  +  (\rho_Q(x)/r)   \right).
\end{align}
  Since
$1-\gamma
\nu>0$, choose now $\eta\in (0, \eta_0)$ so that
$$
c_{18}\,\left(\eta^{1-\gamma
\nu} +\eta
\right)\,\le\, 1\, .
$$
Then  for $x\in  D \cap B(Q, \kappa r/(12))$
with $\rho_Q(x) < \eta r$, we
have by \eqref{e:c:9},
\begin{eqnarray*}
\E_x\left[f(Y^D(\tau_{D_0(x)}));\, Y^D(\tau_{D_0(x)})\notin
B_1(x)\right] &\le & c_{18}\,
f(x_0)\left(\eta^{1-\gamma
\nu}+\eta \right) \le
f(x_0)\, .
\end{eqnarray*}
This completes the proof of \eqref{e:c:4}.

We now prove the Carleson estimate \eqref{e:carleson} for 
$x\in D\cap B(Q, \kappa r/(24))$ 
by a method of contradiction.
Without loss of generality, we may assume that $f(x_0)=1$.
Suppose that there exists $x_1\in D\cap B(Q, \kappa r/(24))$ such
that $f(x_1)\ge K>\eta^{-\beta_1}\vee (1+\delta_*^{-1})$, where $K$ is a
constant to be specified later. By \eqref{e:c:2} and the assumption
$f(x_1)\ge K>\eta^{-\beta_1}$, we have
$(\rho_Q(x_1)/r)^{-\beta_1}>f(x_1)\ge K> \eta^{-\beta_1}$, and hence
$\rho_Q(x_1)<\eta r$.
By (\ref{e:c:3}) and  (\ref{e:c:4}),
$$
K\le f(x_1)\le \E_{x_1}\left[f(Y^D(\tau_{D_0(x_1)}));
Y^D(\tau_{D_0(x_1)}) \in B_1(x_1) \right]+1\, ,
$$
and hence
$$
\E_{x_1}\left[f(Y^D_{\tau_{D_0(x_1)}}); Y^D_{\tau_{D_0(x_1)}} \in
B_1(x_1)\right] \ge f(x_1)-1 > \frac{1}{1+\delta_*}\, f(x_1)\, .
$$
In the last inequality of the display above we used the assumption
that  $f(x_1)\ge K>1+\delta_*^{-1}$. 
If $K \ge (24/\kappa)^{\beta_1/\gamma}$, then
$(\rho_Q(x_1)/r)^\gamma <\kappa/(24)$.
Thus $\overline{ B_1(x_1)}\subset B(Q, \kappa r/(12))$. 
We now get from \eqref{e:c:1} that
\begin{align*}
&\E_{x_1}[f(Y^D(\tau_{D_0(x_1)})), Y^D(\tau_{D_0(x_1)})\in B_1(x_1)]\\
&=\E_{x_1}[
f(Y^D(\tau_{D_0(x_1)})), Y^D(\tau_{D_0(x_1)})\in B_1(x_1)\cap D]\\
& \le \P_x(Y^D(\tau_{D_0(x_1)})\in
D) \, \sup_{B_1(x_1)}f  \le (1-\delta_*) \, \sup_{B_1(x_1)}f \, .
\end{align*}
Therefore, 
$\sup_{B_1(x_1)}f> f(x_1)/(1-\delta_*^2)$, i.e.,
there exists $x_2\in D\cap B(Q, \kappa r/(12))$ such that
$$
|x_1-x_2|\le r^{1-\gamma}\rho_Q(x_1)^{\gamma} \quad \hbox{ and }
\quad
f(x_2)>\frac{1}{1-\delta_*^2}\, f(x_1)\ge \frac{1}{1-\delta_*^2}\, K\, .
$$
Similarly, if $x_k\in D\cap B(Q,  \kappa r/(12))$ 
with $f(x_k)\geq
K/(1-\delta_*^2)^{k-1}$ for $k\ge 2$, then there exists $x_{k+1}\in D$
such that
\begin{equation}\label{e:c:10}
|x_k-x_{k+1}|\le r^{1-\gamma}\rho_Q(x_k)^{\gamma}  \quad \hbox{ and
} \quad f(x_{k+1}) > \frac{1}{1-\delta_*^2}\, f(x_k)>
\frac{1}{(1-\delta_*^2)^k}\, K\, .
\end{equation}
From (\ref{e:c:2}) and (\ref{e:c:10}) it follows that
$\rho_Q(x_{k})/r \le (1-\delta_*^2)^{(k-1)/\beta_1}K^{-1/\beta_1}$, for
every $k\ge 1$. 
Therefore by this and \eqref{e:c:10}, 
\begin{align*}
&|x_k-Q|\,\le\,|x_1-Q|
+\sum_{j=1}^{k-1}|x_{j+1}-x_j|\,\le\, \frac{\kappa r}{24} +
 \sum_{j=1}^{\infty} r^{1-\gamma}\rho_Q(x_j)^{\gamma}\\
&\le \frac{\kappa r}{24}
+r^{1-\gamma}\sum_{j=1}^{\infty}(1-\delta_*^2)^{(j-1)
\gamma/\beta_1}K^{-\gamma/\beta_1}r^{\gamma}\,=
\frac{\kappa r}{24}
+ r K^{-\gamma/\beta_1}\,
\frac{1}{1-(1-\delta_*^2)^{\gamma/\beta_1}}.
\end{align*}
Choose
$$
K=\eta^{-\beta_1}
\vee (1+\delta_*^{-1})\vee 
[(24/\kappa)^{\beta_1/\gamma}(1-(1-\delta_*^2)^{\gamma/\beta_1})^{-\beta_1/\gamma}]
$$
so that $K^{-\gamma/\beta_1}\, (1-(1-\delta_*^2)^{\gamma/\beta_1})^{-1}\le
\kappa/(24)$. Hence $x_k\in D\cap B(Q, \kappa r/(12))$ 
for every $k\ge 1$. 
Since $\lim_{k\to \infty}f(x_k)=\infty$ by \eqref{e:c:10},
 this contradicts the fact that
$f$ is bounded on $B(Q,r/2)$. This contradiction shows that $f(x)< K$
for every $x\in D\cap B(Q, \kappa r/(24))$.
This completes the proof of
the theorem.
 \qed
 
\section{Boundary Harnack principle in $C^{1,1}$ domain}\label{s:bhi}
 
We continue assuming that $d\ge 2$, {\bf (A1)}--{\bf (A3)} hold true, that 
  {\bf (A4)} and {\bf (A5)} also hold if $d=2$,
and that $D\subset \R^d$ is either a bounded $C^{1,1}$ domain, or
a $C^{1,1}$ domain with compact complement or a
domain consisting of all the points above the graph of a bounded  globally $C^{1,1}$
function. 
Unbounded sets are considered only when $d\ge 3$, and then we  assume that {\bf (A6)} holds.
Let $(R, \Lambda)$ be the $C^{1, 1}$ characteristics of $D$.
Without loss of generality we assume that $R \le 1$.
Recall $\kappa=(1+(1+ \Lambda)^2)^{-1/2} $ and $V_z(r)$ is  a $C^{1,1}$
domain  with characteristic $(rR/L, \Lambda L/r)$
such that $D_z( r/2, r/2) \subset V_z(r)  \subset  D_z( r, r) $ where  $L=L(R, \Lambda, d)>0$.
 \begin{lemma}\label{L:ggggg}
There exists a constant  $C=C(R,  \Lambda)>0$ such that for every $r \le  \kappa^{-1} R/2$, $Q \in D$ 
 and $x \in D_Q(2^{-3}r, 2^{-3}r)$,
 \begin{align}
\label{e:D2_43}
&\P_x\Big( Y^D(\tau_{ V_Q(r)}) \in D_Q (  2 r , r)\setminus D_Q (  3 r/2 , r)\Big) 
\ge 
C r^{-3} \phi'( r^{-2} )     \delta_{D}(x)/ \phi( r^{-2} ).
\end{align}
 \end{lemma}
 \pf  Without loss of generality, we assume $Q=0$.  Note that 
$V_0(r) \subset  D_0( r, r) \subset D\cap B(0,  r /\kappa)$. 
For $a, b>0$, we define the cone $\sC(x,a ,b)$ above $x$ by
$$
\sC(x,a ,b):=\{y=(\wt y,y_d)\in B(x,a) \mbox{ in }  CS : y_d>x_d ,|\wt x-\wt y|<b (y_d -x_d)\}.
$$
Since $x \in D_0(2^{-3}r, 2^{-3}r)$, 
there exists $\eta \in (0, (2(1+ \Lambda))^{-2})$ such that \begin{align}
\label{e:Ceta}
\sC(x,2^{-4} r ,\eta) \subset D_0(2^{-2}r, 2^{-2}r).
 \end{align}

Since $G^{Y^D}_{V_0(r)}\ge G^{Y^{V_0(r)}}$ by \eqref{e:UB-smaller-UDB}, 
we have by Proposition \ref{p:gfe},
\begin{align}
\label{e:D2_41}
& \E_{x}
\int_0^{\tau_{ V_0(r)}}\delta_D(Y^D_t) dt =\int_{V_0(r)} G^{Y^D}_{V_0(r)}(x,z) \delta_D(z) dz\nn\\
&\ge   \int_{V_0(r)} G^{Y^{V_0(r)}}(x,z) \delta_D(z) dz\nn\\
&\ge  c_1\int_{V_0(r) }  \left(\frac{
\delta_{V_0(r)}(x)}{|x-z|}\wedge 1
\right)\left(\frac{\delta_{V_0(r)}(z)}{|x-z|}\wedge 1
\right)\frac{\phi'(|x-z|^{-2})}
{|x-z|^{d+2}\phi(|x-z|^{-2})^2} \delta_D(z) dz.
\end{align}
It is easy to see that there exists $c_2\in (0,1]$ such that 
\begin{equation}\label{e:D2_421}
c_2|x-z|\le \delta_{V_0(r)}(z)=
\delta_{D} (z),\qquad z \in 
\overline{\sC(x,2^{-5} r, 2^{-1}\eta)}.
\end{equation}
Thus  by \eqref{e:Ceta} and \eqref{e:D2_421},
 \begin{align}
\label{e:D2_42}
&\int_{V_0(r) }  \left(\frac{
\delta_{V_0(r)}(x)}{|x-z|}\wedge 1
\right)\left(\frac{\delta_{V_0(r)}(z)}{|x-z|}\wedge 1
\right)\frac{\phi'(|x-z|^{-2})}
{|x-z|^{d+2}\phi(|x-z|^{-2})^2} \delta_D(z) dz\nn\\
& \ge c_2  \int_{\sC(x,2^{-5} r, 2^{-1}\eta) }  
\left(\frac{\delta_{D}(x)}{|x-z|}\wedge 1 \right)\frac{\phi'(|x-z|^{-2})}
{|x-z|^{d+2}\phi(|x-z|^{-2})^2}
\delta_{D}(z)dz.
\end{align}

We claim that for $z \in \sC(x,2^{-5} r, 2^{-1}\eta) $,
 \begin{align}
\label{e:D2_420}
\left(\frac{\delta_{D}(x)}{|x-z|}\wedge 1
\right)
\delta_{D}(z)  \ge c_3 \delta_{D}(x).
\end{align}
If $z \in \sC(x,2^{-5}r, 2^{-1}\eta)  \setminus  B(x,   \delta_D(x)/2)$,
then $|x-z| \ge \delta_D(x)/2$, 
so by \eqref{e:D2_421},
 \begin{align*}
\left(\frac{\delta_{D}(x)}{|x-z|}\wedge 1
\right)
\delta_{D}(z) \ge 
c_4\left(\frac{\delta_{D}(x)}{|x-z|}\wedge 1
\right)|x-z|\ge c_5 \delta_{D}(x).
\end{align*}
If $z \in \sC(x,2^{-5}r , 2^{-1}\eta)\cap  B(x,   \delta_D(x)/2)$ 
then $|x-z| < \delta_D(x)/2$ and 
$ 
\delta_{D}(z) \ge \delta_{D}(x) -|x-z| > \delta_{D}(x)/2$. Thus,
 \begin{align*}
\left(\frac{\delta_{D}(x)}{|x-z|}\wedge 1
\right)
\delta_{D}(z)  \ge 
\delta_{D}(z) \ge
\frac{1}2 \delta_{D}(x).
\end{align*}
We have proved \eqref{e:D2_420}.

Combining \eqref{e:D2_41},   \eqref{e:D2_42} and \eqref{e:D2_420}, we get 
 \begin{align}
\label{e:D2_422}
&\E_{x}
\int_0^{\tau_{ V_0(r)}}\delta_D(Y^D_t) dt \nn\\
& \ge c_1c_2c_3 \delta_{D}(x) \int_{\sC(x,2^{-5}r, 2^{-1}\eta)  }  
\frac{\phi'(|x-z|^{-2})}
{|x-z|^{d+2}\phi(|x-z|^{-2})^2} dz\nn\\
 &\ge  c_6 \delta_{D}(x) \int_{0}^{2^{-5}r}
 \frac{\phi'(s^{-2})}
{s^{3}\phi(s^{-2})^2 } ds  
\ge  c_{7}  \delta_{D}(x)/ \phi( r^{-2} ).
\end{align}
We have used Lemma \ref{l:pbf}(a)   in the last inequality.

By Proposition  \ref{p:J^{Y^D}(z,y)}, \eqref{H:1} and Lemma \ref{l:pbf},
 we have that  for $z \in V_0(r)$, 
 \begin{align}
\label{e:Jlow22n}
&\int_{D_0 (  2 r , r)\setminus D_0 (  3 r/2 , r)}     J^{Y^D}(z,y)dy\nn\\ 
&\asymp
\int_{D_0 (  2 r , r)\setminus D_0 (  3 r/2 , r)}  \left(\frac{\delta_D(z)}{|z-y|} \wedge 1 \right)
\left(\frac{\delta_D(y)}{|z-y|} \wedge 1 \right)
k(|y-z|) dy\nn\\
 &\asymp
\int_{D_0 (  2 r , r)\setminus D_0 (  3 r/2 , r)}     \left(\frac{\delta_D(z)}{|y|} \wedge 1 \right)
\left(\frac{\delta_D(y)}{|y|} \wedge 1 \right)
k(|y|) dy\nn\\
&\asymp
\int_{D_0 (  2 r , r)\setminus D_0 (  3 r/2 , r)}     \left(\frac{\delta_D(z)}{|y|} \wedge 1 \right)
\left(\frac{\delta_D(y)}{|y|} \wedge 1 \right)
\frac{\phi'(|y|^{-2})}
{|y|^{d+2}}  dy \nn\\
&\asymp \delta_D(z)
\int_{D_0 (  2 r , r)\setminus D_0 (  3 r/2 , r)}    \frac{\delta_D(y)}{|y|^2}\frac{\phi'(|y|^{-2})}
{|y|^{d+2}}dy     \nn\\
&\asymp  \delta_D(z) r
\int_{D_0 (  2 r , r)\setminus D_0 (  3 r/2 , r)}    \frac{1}{|y|^2}\frac{\phi'(|y|^{-2})}
{|y|^{d+2}}dy    \nn\\
&\ge c_8 \delta_D(z) r
|D_0 (  2 r , r)\setminus D_0 (  3 r/2 , r)|  \frac{1}{r^2}\frac{\phi'(r^{-2})}
{r^{d+2}}
 \ge  c_9  \delta_D(z) r^{-3} \phi'( r^{-2} ).
\end{align}
We have used Lemma \ref{l:pbf}(a) in the first inequality of the last line.

Thus \eqref{e:Jlow22n} and \eqref{e:D2_422} imply
 \begin{align*}
&\P_x\Big( Y^D(\tau_{ V_0(r)}) \in D_0 (  2 r , r)\setminus D_0 (  3 r/2 , r)\Big) =\E_{x}
\int_0^{\tau_{ V_0(r)}}
\int_{D_0 (  2 r , r)\setminus D_0 (  3 r/2 , r)}     J^{Y^D}(Y^D_t,y)dydt\nn\\
&\ge c_{9}  r^{-3} \phi'( r^{-2} )  \E_{x}
\int_0^{\tau_{ V_0(r)}}\delta_D(Y^D_t) dt   
\ge c_{10}  r^{-3} \phi'( r^{-2} )     \delta_{D}(x)/ \phi( r^{-2} ).
\end{align*}
 \qed
\begin{lemma}\label{L:2}
 If the constant $\delta$ in {\bf (A3)} satisfies $0<\delta \le 1/2$, then we assume 
that  {\bf (A7)} holds.
There exists  $C=C(R,  \Lambda)>0$ such that for every $r \le  \kappa^{-1} R/2$, $Q \in \partial D$ and 
$x \in D_Q ( 2^{-3} r , 2^{-4} r )$,
  \bee\label{e:L:2}
\P_{x}\left(Y^D( \tau_{ V_Q (r)}) \in
D_Q (  2 r , 2r)\right) \,\le\, C \,  \, 
 \frac{ \delta_D(x)   \phi'(r^{-2})  }
{r^{3}\phi(r^{-2})}.   
 \eee
\end{lemma}

\pf 
Without loss of generality, we assume $Q=0$ and fix $r \le  \kappa^{-1} R/2$.  
Note that 
$V_0(r) \subset  D_0( r, r) \subset D\cap B(0,  r /\kappa)$.
By using the L\'evy system formula and \eqref{e:exit-time-YD} we get that for any Borel subset $U \subset D_0 (  2 r , 2r)\setminus V_0(r)$,
\begin{align}
\label{e:ub1}
\P_{x}\left(Y^D( \tau_{ V_0 (r)}) \in
U \right)&=
\E_{x}
\int_0^{\tau_{ V_0(r)}}
\int_{U}    J^{Y^D}(Y^D_t,z)dzdt\nn\\
&=
\int_{ V_0(r)} G^{Y^D}_{V_0(r)}(x,y)
\int_{U}    J^{Y^D}(y,z)dzdy\nn\\
& \le 
\int_{U}  \int_{ V_0(r)} G^{Y^D}(x,y)
   J^{Y^D}(y,z)dy dz.
\end{align}
Recall that $g$ and $j$ are defined in \eqref{e:defofgandj}.
By Lemma \ref{l:pbf}, Propositions \ref{p:gfe0} and \ref{p:J^{Y^D}(z,y)}, 
for $z \in D_0 (  2 r , 2r)\setminus V_0(r)$,
\begin{align}
\label{e:ub2}
&\int_{ V_0(r)} G^{Y^D}(x,y)
   J^{Y^D}(y,z)dy \nn\\
   &\le 
 c_0  \int_{ V_0(r)}
   \left(\frac{\delta_D(x)}{|x-y|}\wedge 1
\right)\left(\frac{\delta_D(y)}{|x-y|}\wedge 1
\right)\frac{\phi'(|x-y|^{-2})}
{|x-y|^{d+2}\phi(|x-y|^{-2})^2}
    \nn\\
&\quad \times \left(\frac{\delta_D(y)}{|z-y|}\wedge 1
\right)
\left(\frac{\delta_D(z)}{|z-y|}\wedge 1
\right)
\frac{\phi'(|y-z|^{-2})}
{|y-z|^{d+2}}dy\nn\\
 & \le c_1
 \frac{\delta_D(x)}{|x-z|^2}g(|x-z|)
\int_{ V_0(r) \cap \{  |x-z|  \le 2|x-y| \}}
  {\delta_D(y)}   
\left(\frac{\delta_D(z)}{|y-z|}\wedge 1\right)
j(|y-z|)dy    \nn\\
&\quad+ c_1
   \int_{ V_0(r) \cap \{  |x-z|  > 2|x-y| \}}
 \left(\frac{\delta_D(x)}{|x-y|}\wedge 1\right)\left(\frac{\delta_D(y)}{|y-z|}\wedge 1\right) 
\left(\frac{\delta_D(z)}{|y-z|}\wedge 1\right)
 g(|x-y|) j(|y-z|)dy \nn\\
 &=: c_1(I+II).
\end{align}

  By Lemma \ref{l:detlax}, for $z \in D_0 (  2 r , 2r)\setminus V_0(r)$, 
\begin{align}
\label{e:ub3}
&I \le  2 \frac{\delta_D(x)\delta_D(z)}{|x-z|^2}g(|x-z|)
\int_{ V_0(r) \cap \{  |x-z|  \le 2|x-y| \}}j(|y-z|)dy.
\end{align}
Since 
$$
\int_{ V_0(r) \cap \{  |x-z|  \le 2|x-y| \}}j(|y-z|)dy \le c_2 
\int_{\delta_{V_0(r)}(z)}^{\delta_{V_0(r)}(z) +3r }
s^{d-1}\frac{\phi'(s^{-2})}
{s^{d+2}} ds
 \le c_3 \phi(\delta_{V_0(r)}(z)^{-2}), 
$$
combining this with \eqref{e:ub3}
we get
\begin{align}
\label{e:ub5}
&I \le c_4 \phi(\delta_{V_0(r)}(z)^{-2})
\frac{\delta_D(x)\delta_D(z) \phi'(|x-z|^{-2})}
{|x-z|^{d+4}\phi(|x-z|^{-2})^2}.
\end{align}

On the other hand, when $|x-z| \ge 2|x-y|$, we have 
\begin{align}\label{e:neww0}
|y-z| \ge |x-z| - |x-y| \ge \frac12 |x-z| \ge |x-y|,
\end{align}
implying
\begin{equation}\label{e:neww0-0}
\frac23|y-z| \le |x-z| \le 2|y-z|.
\end{equation}
Moreover, by Lemma \ref{l:detlax},
\begin{align}\label{e:neww00}
  & \left(\frac{\delta_D(x)}{|x-y|}\wedge 1
\right)\left(\frac{\delta_D(y)}{|y-z|}\wedge 1
\right)\left(\frac{\delta_D(z)}{|y-z|}\wedge 1
\right)\nn\\
&\le \frac{\delta_D(z)}{|y-z|^2}
\delta_D(y)    \left(\frac{\delta_D(x)}{|x-y|}\wedge 1
\right)  
\le 2\frac{\delta_D(x)\delta_D(z)}{|y-z|^2}.
\end{align}
Thus using \eqref{e:neww00} first, 
and then Lemma \ref{l:pbf}(c) with \eqref{e:neww0-0}, 
\begin{align}
\label{e:ub10}
&II \le
c_5 \delta_D(x)\delta_D(z)
 \int_{ V_0(r) \cap \{  |x-z|  > 2|x-y| \}}
g(|x-y|)
\frac{\phi'(|y-z|^{-2})}
{|y-z|^{d+4}} dy  \nn\\
&=
c_5 \delta_D(x)\delta_D(z)
 \int_{ V_0(r) \cap \{  |x-z|  > 2|x-y| \}}
g(|x-y|)
\frac{\phi'(|y-z|^{-2})}
{|y-z|^{d+4}\phi(|y-z|^{-2})^2}\phi(|y-z|^{-2})^2 dy  \nn\\
&\le c_6
\delta_D(x)\delta_D(z)\frac{\phi'(|x-z|^{-2})}
{|x-z|^{d+4}\phi(|x-z|^{-2})^2}  
 \int_{ V_0(r) \cap \{  |x-z|  > 2|x-y| \}}
g(|x-y|)
\phi(|y-z|^{-2})^2
dy.
\end{align}

Let $a:=|x-z|$. 
By the triangle inequality,
\begin{align}
\label{e:ub11}
&\int_{ V_0(r) \cap \{  |x-z|  > 2|x-y| \}}
g(|x-y|)
\phi(|y-z|^{-2})^2
dy\nn\\
&\le \int_{ V_0(r) \cap \{  |x-z|  > 2|x-y| \}}
g(|x-y|)
\phi(||x-z|-|x-y||^{-2}\wedge|x-y|^{-2})^2
dy
\nn\\
&\le c_7 \int_0^{\kappa^{-1} r}
\frac{\phi'(r^{-2})}
{r^{3}\phi(r^{-2})^2}
\phi(|a-r|^{-2}\wedge r^{-2})^2
dy.
\end{align}
Note that
\begin{align}
\label{e:ub12}
&  \int_0^{\kappa^{-1} r}
\frac{\phi'(r^{-2})}
{r^{3}\phi(r^{-2})^2}
\phi(|a-r|^{-2}\wedge r^{-2})^2
dy\nn\\ &\le 
\int_0^{a/2} 
\frac{\phi'(r^{-2})}
{r^{3}\phi(r^{-2})^2}
\phi(|a-r|^{-2})^2
dy+
\int_{a/2}^\infty \frac{\phi'(r^{-2})}
{r^{3}\phi(r^{-2})^2}
\phi(r^{-2})^2
dy\nn\\ &\le 
\phi(4a^{-2})^2
\int_0^{a/2} 
\frac{\phi'(r^{-2})}
{r^{3}\phi(r^{-2})^2}
dy+
\int_{a/2}^\infty \frac{\phi'(r^{-2})}
{r^{3}}
dy\nn\\ &= 2^{-1}
\phi(4a^{-2})^2
\int_0^{a/2} 
\left(\frac{1}{\phi(r^{-2})}\right)'
dy+2^{-1}
\int_{a/2}^\infty   (-\phi(r^{-2}))'
dy\nn\\ &= 
 2^{-1} \phi(4a^{-2})^2
\frac{1}{\phi(4a^{-2})}+
 2^{-1} \phi(4a^{-2}) \le c_8  \phi(|x-z|^{-2}).
\end{align}
Therefore
\begin{align}
\label{e:ub13}
&II \le  c_9  \delta_D(x)\delta_D(z)\frac{\phi'(|x-z|^{-2})}
{|x-z|^{d+4}\phi(|x-z|^{-2})^2}     \phi(|x-z|^{-2}) \nn\\
&\le c_{10} \phi(\delta_{V_0(r)}(z)^{-2})
\frac{\delta_D(x)\delta_D(z) \phi'(|x-z|^{-2})}
{|x-z|^{d+4}\phi(|x-z|^{-2})^2}.
\end{align}

Now for $x \in D_0 (2^{-3}r,  2^{-4}r)$,   
 we have 
$c_{11} r <|x-z| \le (1+\kappa^{-1})r$ for 
$z \in D_0 (  2 r , 2r)\setminus V_0(r)$.
Thus putting \eqref{e:ub1}, \eqref{e:ub2}, \eqref{e:ub5} and \eqref{e:ub13} together, 
and using Lemma \ref{l:pbf}, we see that for any Borel subset $U \subset D_0 (  2 r , 2r)\setminus V_0(r)$
\begin{align}
\label{e:ub14}
&\P_{x}\left(Y^D( \tau_{ V_0 (r)}) \in
U\right) \nn\\
&\le c_{12} \delta_D(x) \int_{U} 
\phi(\delta_{V_0(r)}(z)^{-2})
\frac{\delta_D(z) \phi'(|x-z|^{-2})}
{|x-z|^{d+4}\phi(|x-z|^{-2})^2}dz \nn \\
&\le c_{13}
\frac{ \delta_D(x)   \phi'(r^{-2})  }
{r^{d+4}\phi(r^{-2})^2}
\int_{U} 
\phi(\delta_{V_0(r)}(z)^{-2})
{\delta_D(z) }dz.
\end{align}

We now deal with two cases separately.

\noindent
{\it Case 1: The constant $\delta$ in {\bf (A3)} satisfies $1 \ge \delta > 1/2$.}
By the co-area formula, 
\begin{align}
\label{e:ub15}
& \int_{D_0 (  2 r , 2r)\setminus V_0(r)} 
\phi(\delta_{V_0(r)}(z)^{-2})
{\delta_D(z) }dz
 \le c_{14} \int_{0}^{2\kappa^{-1} r}\int_{0}^{2\kappa^{-1} r} s\phi(t ^{-2}) r^{d-2}dsdt\nn\\
 & =
 c_{15} r^d \int_{0}^{2\kappa^{-1} r}  \phi(t ^{-2})dt.
\end{align}
Since 
$1 \ge \delta > 1/2$
we choose $\eps \in (0, (\delta-1/2)/2)$ so that $2(1-\delta +\eps) <1$.
Using \cite[Lemma 3.2]{KM} with this $\eps>0$, 
we have 
\begin{align}
\label{e:ub15-1}
&\int_{0}^{2\kappa^{-1} r}  \phi(t ^{-2}) dt  \le \int_{0}^{ r}  \phi(t ^{-2}) dt+ 
 \phi(r^{-2}) \int_{r}^{2\kappa^{-1} r} dt
\nn\\  \le &\phi(r^{-2}) \int_{0}^{ r}  \frac{\phi(t ^{-2})}{\phi(r^{-2})} dt  + \phi(r^{-2})(2\kappa^{-1} -1)r   \nn\\
\le& c_{16}\phi(r^{-2})  \left(\int_{0}^{ r}  \left( \frac{t ^{-2}}{r^{-2}} \right)^{1-\delta+\eps} dt + r \right) \nn\\
=&c_{16} \phi(r^{-2})  \left( r ^{2(1-\delta+\eps)} \int_{0}^{ r}  t^{-2+2\delta-2\eps}dt + r \right) 
\le c_{17}r \phi(r^{-2}) .
\end{align}
Combining \eqref{e:ub14}--\eqref{e:ub15-1},  we conclude that 
\begin{align*}
&\P_{x}\left(Y^D( \tau_{ V_0 (r)}) \in
D_0 (  2 r , 2r)\right) \le  c_{11} c_{14} c_{16} \frac{ \delta_D(x)   \phi'(r^{-2})  }
{r^{d+4}\phi(r^{-2})^2}    r^d  r \phi(r ^{-2})\\ 
&= c_{13} c_{15} c_{17} \frac{ \delta_D(x)   \phi'(r^{-2})  }
{r^{3}\phi(r^{-2})}.   
\end{align*}

\noindent
{\it Case 2: The constant $\delta$ in {\bf (A3)} satisfies $0<\delta \le 1/2$ and  {\bf (A7)} holds.}
By \eqref{e:ub14} we have that for $x \in D_0 (2 ^{-3}r,  2^{-4}r)$,
\begin{align}
\label{e:D2_43nu}
&\P_x\Big( Y^D(\tau_{ V_0(r)}) \in D_0 (  2 r , r)\setminus D_0 (  3 r/2 , r)\Big) 
\nn\\
\le&  c_{13}
\frac{ \delta_D(x)   \phi'(r^{-2})  }
{r^{d+4}\phi(r^{-2})^2}
\int_{D_0 (  2 r , r)\setminus D_0 (  3 r/2 , r)} 
\phi(\delta_{ V_0(r)}(z)^{-2})
{\delta_D(z) }dz\nn\\
\le&  c_{13}
\frac{ \delta_D(x)   \phi'(r^{-2})  }
{r^{d+4}\phi(r^{-2})^2}
|D_0 (  2 r , r)\setminus D_0 (  3 r/2 , r)|
\phi(4r^{-2}) 2r/\kappa \nn\\
\le &
c_{18}  r^{-3} \phi'( r^{-2} )     \delta_{D}(x)/ \phi( r^{-2} ). 
\end{align}

Note that, using Lemma \ref{l:detlax} and Proposition \ref{p:gfe0}, we get 
that for $w \in D_0(r/2, r/2)$ and any $\tilde{r}\in (0,r/2)$,
\begin{align}
\label{e:D2_4n1}
& \int_{D \cap B(w,  \tilde{r})} G^{Y^D}(w,z)\delta_D(z)dz\nn\\
&\le
c_{19}  \int_{B(w,  \tilde{r})}  \delta_D(z)\left(\frac{\delta_D(w)}{|w-z|}\wedge 1
\right)\frac{\phi'(|w-z|^{-2})}
{|w-z|^{d+2}\phi(|w-z|^{-2})^2}  dz\nn\\
&\le
c_{20}  \delta_D(w) \int_{B(w,  \tilde{r})} \frac{\phi'(|w-z|^{-2})}
{|w-z|^{d+2}\phi(|w-z|^{-2})^2}  dz\le   
c_{21}   \delta_D(w)/ \phi( \tilde{r}^{-2}).
\end{align}
Recall that $a_1$ is defined in \eqref{e:w_1}.
Using \eqref{e:D2_4n1}, Proposition \ref{p:J^{Y^D}(z,y)}, Theorem  \ref{t:Jfe:nn} and  {\bf (A7)}, we have that for $w \in D_0(r/2, r/2)$ and  $0 <4R_1 \le R_2 <r $,
\begin{align}\label{e:ggggg1}
&\P_w \left(Y^D_{\tau_{D \cap B(w, R_1)}} \in D\setminus B(w, R_2)\right)\nn\\
\le& 
\int_{D \cap B(w, R_1)} G_{D \cap B(w, R_1)}^{Y^D}(w,z) \int_{D\setminus B(w, R_2)} J^{Y^D}(y,z) dydz\nn
\\
\le&c_{22}
\int_{D \cap B(w, R_1)} G^{Y^D}(w,z)\delta_D(z)dz \int_{\R^d\setminus B(0, R_2)}
\frac{\phi'(|y|^{-2})}
{|y|_1 |y|^{d+2}}
dy
\nn
\\
\le & c_{23}\frac{\delta_D(w)}{ \phi(R_1^{-2})   } \left( \int_{R_2}^{1} s^{-4}\phi'(s^{-2}) ds + \int_{1}^{\infty} s^{-3}\phi'(s^{-2}) ds\right) \le c_{24}\frac{\delta_D(w)}{ \phi(R_1^{-2})   } R_2^{-1}\phi(R_2^{-2}). 
\end{align}

Let 
$$H_2=\{Y^D( \tau_{ V_0(r)}) \in
D_0 (  2 r , 2r)  \}, \quad 
H_1=\{Y^D( \tau_{ V_0(r)}) \in
D_0 (  2 r , r)\setminus D_0 (  3 r/2 , r)\}. 
$$
We claim that $\P_x( H_2)\leq c_{25} \P_x( H_1)$ 
for all $r \le  \kappa^{-1} R/2$ and $x \in D_0 ( 2^{-3} r , 2^{-4} r )$. 
Combining this claim with \eqref{e:D2_43nu}
we arrive at the conclusion of the lemma:
$$
\P_x\left( Y^D( \tau_{ V_0(r)}) \in D_0 (  2 r , 2r) \right)
\leq c_{25} \P_x( H_1) \le 
c_{26} \frac{ \delta_D(x)   \phi'(r^{-2})  }
{r^{3}\phi(r^{-2})}
\quad x \in D_0 ( 2^{-3} r , 2^{-4} r )\,   .
$$

Now we give the proof of the claim, which is inspired by the proof of \cite[Lemma 5.3]{G}.
Note that, by \cite[Lemma 1.3]{KM14} and \cite[Lemma 3.2]{KM},  we have $\phi'(r) r 
\asymp \phi(r)$ for $r \ge 1$. 
Thus, by Lemma \ref{L:ggggg}, we have that for $w \in D_0(r/8, r/8)$,
\begin{align}\label{e:D2_43n2}
&\P_w ( H_1 ) 
\ge c_{27}  r^{-3} \phi'( r^{-2} )     \delta_{D}(w)/ \phi( r^{-2} ) \asymp \delta_{D}(w)/ r.
\end{align}

We define, for $i\ge 1$,
$$
J_i=D_0(2^{-i-2}r,s_i)\setminus  D_0(2^{-i-3}r,s_i),\ \ \ \ s_i=\frac{1}{4}\left(\frac{1}{2}-\frac{1}{50}\sum_{j=1}^i\frac{1}{j^2}\right)r,
$$
and $s_0=s_1$.  
Note that $r/(10)<s_i <r/8$.
For $i\geq 1$,  set
\begin{align}\label{q1q}
d_i=d_i(r)=\sup_{z\in J_i}\P_z( H_2)/\P_z( H_1),\ \ \
\widetilde{J}_i=D_0(2^{-i-2}r,s_{i-1}),\ \ \ \
\tau_i=\tau_{\widetilde{J}_i}.
\end{align}
It follows from \eqref{e:D2_43n2} that
$\sup_{r \le \kappa^{-1} R/2} d_i(r)$ is finite for all $i \ge 1$.
Using the strong Markov property  and the fact that $\tau_i\leq \tau_{V_0(r)}$, we get that
for $z\in J_i$ and $i\geq 2$,
\begin{align*}  
&\P_z\left(H_2, \ Y^{D}_{\tau_i}\in \cup_{k=1}^{i-1}J_k\right)=\sum_{k=1}^{i-1}  \P_z\left(Y^{D}_{\tau_{V_0(r)}}\in D_0 (  2 r , 2r), \ Y^{D}_{\tau_i}\in J_k\right)\\
&= \sum_{k=1}^{i-1}
 \E_z\left[ \P_{Y^{D}_{\tau_i}}(H_2) :  \ Y^{D}_{\tau_i}\in J_k \right]
\leq  \sum_{k=1}^{i-1}  d_k \E_z\left[ \P_{Y^{D}_{\tau_i}}(H_1) :  \ Y^{D}_{\tau_i}\in J_k\right] \nonumber\\
&  \le   \left ( \sup_{1\leq k\leq i-1}  d_k \right)\sum_{k=1}^{i-1}  \E_z\left[ \P_{Y^{D}_{\tau_i}}(H_1) :  \ Y^{D}_{\tau_i}\in J_k\right]   = 
  \left ( \sup_{1\leq k\leq i-1}  d_k \right) \P_{z}\left(H_1, \ Y^{D}_{\tau_i}\in \cup_{k=1}^{i-1}J_k\right) \\
&\le 
  \left ( \sup_{1\leq k\leq i-1}  d_k \right) \P_{z}(H_1).
\end{align*}
Thus, for $z\in J_i$ and $i\geq 2$,
\begin{align}\label{q12q}  
\P_z( H_2) = 
& \P_z\left(H_2, \ Y^{D}_{\tau_i}\in \cup_{k=1}^{i-1}J_k\right)+\P_z\left( Y^{D}_{\tau_i}\in D_0 (  2 r , 2r) \setminus
\cup_{k=1}^{i-1}J_k\right)\nonumber\\
\leq  &  \left ( \sup_{1\leq k\leq i-1}  d_k \right) \P_{z}(H_1) +\P_z\left( Y^{D}_{\tau_i}\in
D_0 (  2 r , 2r) \setminus \cup_{k=1}^{i-1}J_k\right).
\end{align}

For  $i\ge 2$, define $\sigma_{i, 0}=0, \sigma_{i, 1}=\inf\{t>0: |Y^{D}_t-Y^{D}_0|\geq 2^{-i-2}r\} $ and
$\sigma_{i, m+1}=\sigma_{i,1}\circ\theta_{\sigma_{i, m}}$
for $m\geq 1$. By Lemma \ref{lower bound},  we have that there exists $k_1 \in (0,1)$ such that 
\begin{align}\label{e:rs}
\P_{w}\big(Y^{D}_{\sigma_{i, 1}}\in \widetilde{J}_i\big) \le 1-\P_{w}( \sigma_{i, 1}=\zeta ) \le 1-\P_{w}\big(\tau^{Y^D}_{D\cap B(w,\rho_Q(w)/2)}=\zeta \big)  <k_1,\ \ \ w\in  \widetilde{J}_i.
\end{align}
For the purpose 
of  further estimates, we now choose a 
positive integer $l$ 
such that $k_1^l\le 2^{\gamma-2}$, where $\gamma$ is the constant from \textbf{(A7)}. Next we choose $i_0 \ge 2$ large enough so that 
$2^{-i}<1/(200 l i^3)$ for all $i\ge i_0$,  
and that
$$
200^3 c_{24} \sigma_2^{-1} l^{4-2\gamma}\le 
2^{(i+2)(1-\gamma)-1}\, , \quad i\ge i_0.
$$
Here $c_{24}$ is the constant from \eqref{e:ggggg1} and $\sigma_2$ the constant from \textbf{(A7)}. 
With $l$ and $i_0$ so chosen, we have that
\begin{equation}\label{e:precise-inequality}
k_1^{li}+200^3 c_{24} \sigma_2^{-1} 
\frac{l^{4-2\gamma}i^{4+6(1-\gamma)}}
{2^{2(i+2)(1-\gamma)+i}}\le 
\frac{i^{4+6(1-\gamma)}}{2^{(i+1)(1-\gamma)+i}}, \quad i\ge i_0.
\end{equation}
Now we assume $i\ge i_0$. Using \eqref{e:rs} and the strong Markov property we have  that  for   $z\in J_i$,
\begin{align}\label{q132q}    
&\P_z( \tau_{i}>\sigma_{i, li})\leq 
\P_z\big(Y^{D}_{\sigma_{i, k}}\in \widetilde{J}_i, 1\leq k\leq
li\big )\nonumber\\
&=
 \E_z \left[ \P_{Y^{D}_{\sigma_{i, li-1}}} (Y^{D}_{\sigma_{i, 1}}\in  \widetilde{J}_i) : Y^{D}_{\sigma_{i, li-1}} \in  \widetilde{J}_i,  Y^{D}_{\sigma_{k}}\in \widetilde{J}_i, 1\leq k\leq
li-2
  \right]
 \nonumber\\      
& \leq \P_z\big(Y^{D}_{\sigma_{i, k}}\in \widetilde{J}_i, 1\leq k\leq
li-1 \big )k_1\leq k_1^{li}.
\end{align}
Note that if $z\in J_i$ and $y\in D_0 (  2 r , 2r) \setminus[ \widetilde{J}_i \cup(\cup_{k=1}^{i-1}J_k)]$, 
then $|y-z|\ge (s_{i-1}-s_i) \wedge (2^{-3}-2^{-i-2}) r = r/(200 i^2)$. 
Furthermore, since  $2^{-i-2} r< r/(200 i^2)$, $\tau_i$ must be one of  the
$\sigma_{i, k}$'s, $k\le li$. 
Hence, on
$\{Y^{D}_{\tau_{i}}\in D_0 (  2 r , 2r) \setminus \cup_{k=1}^{i-1}J_k,\ \ \tau_{i}\leq
\sigma_{i, li}\}$ with $Y^{D}_0=z\in J_i$,  there exists $k$, $1\le k\le li$, such that 
$|Y^D_{\sigma_{i, k}}-Y^D_0|=|Y^D_{\tau_i}-Y^D_0|> r/(200 i^2)$.
Thus  for some $1\leq k\leq li$,
$$
\sum_{j=1}^k\big|Y^D_{\sigma_{i, j}}-Y^D_{\sigma_{i, j-1}}\big|>
\frac{r}{200i^2}\, .
$$
which implies for some $1\leq k'\leq k\le  li$, 
$$ 
\big|Y^{D}_{\sigma_{i, k'}}-Y^{D}_{\sigma_{i, k'-1}}\big|\geq
\frac1k \frac{r}{200i^2}\ge \frac{1}{ li} \frac{r}{200i^2}.
$$
Thus, we have 
\begin{align}& \{Y^{D}_{\tau_{i}}\in D_0 (  2 r , 2r) \setminus
\cup_{k=1}^{i-1}J_k,\ \ \tau_{i}\leq \sigma_{i, li}\}\nonumber\\
 \subset &   \cup_{k=1}^{li}\{|
Y^D_{\sigma_{i, k}}- Y^D_{\sigma_{i, k-1}}|\geq r/(200li^3),
Y^D_{\sigma_{i, k-1}}\in \widetilde{J}_{i}\},\ \ \ 
i\geq i_0\, .
\end{align}
  Using  the strong Markov property  and then using \eqref{e:ggggg1} (noting that  
$4 \cdot 2^{-i-2}<1/(200 l i^3)$ for all $i\ge i_0$)
 we have 
 \begin{align}\label{dsa}& \P_z \left(Y^{D}_{\tau_{i}}\in D_0 (  2 r , 2r ) \setminus
\cup_{k=1}^{i-1}J_k,\ \ \tau_{i}\leq \sigma_{i, li} \right)\nonumber\\
 \leq  &\sum_{k=1}^{li} \P_z\left(|
Y^D_{\sigma_{i, k}}- Y^D_{\sigma_{i, k-1}}|\geq r/(200li^3),
Y^D_{\sigma_{i, k-1}}\in \widetilde{J}_{i}
\right ) \nonumber\\
 \leq  &li \sup_{z\in \widetilde{J}_{i}} \P_z\left(|
Y^D_{\sigma_{i, 1}} - z |\geq r/(200li^3) \right)
\nonumber\\
  \leq  &c_{24}li  \frac{2^{-i}r}{ \phi(2^{2(i+2)} r^{-2})   } 
 \frac{200 li^3}{r} \phi((200 li^3)^2/r^2)\nn\\
= & c_{24}li   \frac{2^{-i}}{ \phi(2^{2(i+2)} r^{-2})   } 
(200 li^3)\phi((200 li^3)^2/r^2).
\end{align}
By (\ref{q132q}), (\ref{dsa}),  {\bf (A7)}, \eqref{e:wlsc-substitute} and \eqref{e:precise-inequality},  for $z\in J_i$, $i\ge i_0$, we have
\begin{align}\label{as}
&\P_z\left( Y^{D}_{\tau_{i}}\in D_0 (  2 r , 2r) \setminus \cup_{k=1}^{i-1}J_k \right) \leq
k_1^{li}+c_{24} li  
\frac{2^{-i}}{ \phi(2^{2(i+1)} r^{-2})   } (200li^3)\phi((200li^3)^2/r^2)\nn\\
& \le k_1^{li}+200^3 c_{24} \sigma^{-1}_2l^{4-2\gamma} 
\frac{ i^{4+6(1-\gamma )} }{ 2^{2(i+1)(1-\gamma )+i}}   \le    \frac{ i^{4+6(1-\gamma )} }{ 2^{(i+1)(1-\gamma )+i}}
\le \frac{ i^{4+6(1-\gamma )} }{ 2^{(2-\gamma )i}}. 
\end{align} 
Combining the above with  \eqref{e:D2_43n2}, we get for 
$z\in J_i, i\geq i_0$,
\begin{align}\label{as1}
\frac{\P_z( Y^{D}_{\tau_i}\in
D_0 (  2 r , 2r) \setminus \cup_{k=1}^{i-1}J_k)}{\P_z(H_1)} \le  
   c_{28}   \frac{ i^{4+6(1-\gamma )} }{ 2^{(1-\gamma )i}}.
\end{align} 
By this, (\ref{q12q}) and (\ref{as}), for $z\in J_i$ 
and $i\geq i_0$,
\begin{align}  
&\frac{\P_z( H_2)}{\P_z(H_1)} \leq  \sup_{1\leq k\leq i-1}  d_k   +\frac{\P_z( Y^{D}_{\tau_i}\in
D_0 (  2 r , 2r) \setminus \cup_{k=1}^{i-1}J_k)}{\P_z(H_1)} \nn\\
\leq& \sup_{1\leq k\leq i-1}d_k+  c_{28}   \frac{ i^{4+6(1-\gamma )} }{ 2^{(1-\gamma )i}}
.\nonumber
\end{align}
This implies that
\begin{align}  
d_i& \leq  \sup_{1\leq k\leq i_0-1} d_k+ 
c_{28} \sum_{k=1}^i\frac{ k^{4+6(1-\gamma )} }{ 2^{(1-\gamma )k}}  \leq \sup_{1\leq k\leq i_0-1\atop r \le \kappa^{-1} R/2} d_k(r)+   c_{28} \sum_{k=1}^{\infty} \frac{ k^{4+6(1-\gamma )} }{ 2^{(1-\gamma )k}}=:c_{29} <\infty. \nonumber
\end{align}
Thus the claim above is valid, since $D_0 ( 2^{-3} r , 2^{-4} r )\subset \cup_{k=1}^\infty J_k$. The proof is now complete.
\qed

\noindent
\textbf{Proof of Theorem \ref{t:main}}.
 In this proof, the constants $\eta$ and $c_i$
are always independent of $r$.

Note that,  since $D$ is a $C^{1,1}$ domain and $r<R$, by Theorem \ref{uhp} 
and a standard chain argument, it
suffices to prove \eqref{e:bhp_m} for $x,y \in D \cap B(Q,2^{-7} \kappa r)$. Throughout the remainder of the proof 
we assume that  $x \in D \cap B(Q,2^{-7} \kappa r)$.

Let $Q_x$ be
the point in $ \partial D$ 
so that $|x-Q_x|=\delta_{D}(x)$ and
let $x_0:=Q_x+\frac{r}{8}(x-Q_x)/|x-Q_x|$. We choose a
$C^{1,1}$ function $\varphi: \bR^{d-1}\to \bR$ satisfying $\varphi
(\wt 0)= 0$, $\nabla\varphi (\wt 0)=(0, \dots, 0)$, $\| \nabla \varphi
 \|_\infty \leq \Lambda$, $| \nabla \varphi (\wt y)-\nabla \varphi (\wt z)|
\leq \Lambda |\wt y-\wt z|$, and an orthonormal coordinate system $CS$ with
its origin at $Q_x$ such that
$$
B(Q_x, R)\cap D=\{ y=(\wt y, y_d) \in B(0, R) \mbox{ in } CS: y_d >
\varphi (\wt y) \}.
$$
In the coordinate system $CS$ we have $\wt x = \wt 0$ and $x_0=(\wt
0, r/8)$. For any $b_1, b_2>0$, we define
$$
D(b_1, b_2):=\left\{ y=(\wt y, y_d) \mbox{ in } CS: 0<y_d-\varphi(\wt
y)< 2^{-2}\kappa r b_1, \ |\wt y| <  2^{-2}\kappa r b_2\right\}.
$$

By \eqref{e:Uzr}, we have 
 that
$D(2, 2)\subset D\cap B(Q_x, r/2)\subset D\cap B(Q, r)$.
Thus, since
$f$ is
harmonic in $D\cap B(Q,r)$ and vanishes continuously in
$\partial D\cap B(Q, r)$, by Lemma \ref{l:regularity}, $f$ is regular
harmonic in 
$D(2, 2)$ and vanishes continuously in
$\partial D\cap D(2, 2)$. 

Recall that $V(1):=V_{Q_x}(2^{-2}\kappa r)$ is a $C^{1,1}$
domain  with $C^{1,1}$ characteristics $(rR/L, \Lambda L/r)$
such that $D( 1/2, 1/2) \subset V(1)  \subset  D( 1, 1) $,
where $L=L(R, \Lambda, d)>0$.
By the  Harnack inequality and Lemma \ref{L:ggggg}, we have
\begin{align}
&f(x) =  \E_x\left[f\big(Y^D(\tau_{ V(1)})\big)\right] \ge
\E_x\left[f\big(Y^D(\tau_{ V(1)})\big); Y^D_{
\tau_{ V(1)}} \in  
D(2,1)\setminus D(3/2,1)\right] \nn\\
&\ge c_{12} f(x_0) \P_x\Big( Y^D(\tau_{ V(1)}) \in 
D(2,1)\setminus D(3/2,1)\Big)
\ge 
c_{11} c_{12} f(x_0) \frac{ \delta_D(x)   \phi'(r^{-2})  }
{r^{3}\phi(r^{-2})}.
\label{e:BHP2}
\end{align}

Recall that $a_1$ is defined in \eqref{e:w_1}.
By Lemma \ref{l:pbf}, \eqref{e:upper-estimate-mu}, \eqref{e:lower-estimate-mu}, 
\eqref{H:1}, \eqref{H:1a},
Proposition \ref{p:J^{Y^D}(z,y)} and Theorem  \ref{t:Jfe:nn}, we have that  for $z \in V(1)$, 
 \begin{align}
\label{e:Jlow2}
&\int_{D\setminus D(2,2)} f(y)     J^{Y^D}(z,y)dy\nn\\ &\asymp
\int_{D\setminus D(2,2)} f(y)     \left(\frac{\delta_D(z)}{|z-y|_1} \wedge 1 \right)
\left(\frac{\delta_D(y)}{|z-y|_1} \wedge 1 \right)\frac{\mu(|z-y|^2)}
{|z-y|^{d-2}}    dy\nn\\
&\asymp
\int_{D\setminus D(2,2)} f(y)     \left(\frac{\delta_D(z)}{|y|_1} \wedge 1 \right)
\left(\frac{\delta_D(y)}{|y|_1} \wedge 1 \right)\frac{\mu(|y|^2)}
{|y|^{d-2}}    dy \nn\\
&\asymp \delta_D(z)
\int_{D\setminus D(2,2)} f(y)   \frac{1}{|y|_1} \left(\frac{\delta_D(y)}{|y|_1} \wedge 1 \right)  \frac{\mu(|y|^2)}
{|y|^{d-2}}    dy.
\end{align}

Since 
\begin{align*}
&\E_{x}\left[f \left(Y^D(\tau_{ V(1)})\right); \,
Y^D(\tau_{ V(1)})
\notin  D(2,2)\right]\,=\, \E_{x}
\int_0^{\tau_{ V(1)}}
\int_{D\setminus
D(2,2)}  f(y)     J^{Y^D}(Y^D_t, y)dydt, \end{align*} 
by  \eqref{e:Jlow2}, we have 
 \begin{align}
\label{e:D2_2}
&\E_{x}\left[f \left(Y^D(\tau_{ V(1)})\right); \,
Y^D(\tau_{ V(1)})
\notin  D(2,2)\right]\nn \\&\asymp \E_{x}
\int_0^{\tau_{ V(1)}}\delta_D(Y^D_t) dt  
\int_{D\setminus D(2,2)} f(y)    
 \frac{1}{|y|_1} \left(\frac{\delta_D(y)}{|y|_1} \wedge 1 \right)   \frac{\mu(|y|^2)}
{|y|^{d-2}}    dy.
\end{align}

By Proposition \ref{p:gfe0},  for $x,z \in  V(1)$,
$$
G^{Y^D}(x,z)
\le 
c_{13} \left(\frac{\delta_D(x)}{|x-z|}\wedge 1
\right)\left(\frac{\delta_D(z)}{|x-z|}\wedge 1
\right)\frac{\phi'(|x-z|^{-2})}
{|x-z|^{d+2}\phi(|x-z|^{-2})^2}.
$$
Thus, using Lemma \ref{l:detlax},
\begin{align}
\label{e:D2_4}
&\E_{x}
\int_0^{\tau_{ V(1)}}\delta_D(Y^D_t) dt 
=\int_{V(1)}  \delta_D(z) G^{Y^D}_{V(1)}(x,z) dz\,\le\,   \int_{V(1)} \delta_D(z) G^{Y^D}(x,z) dz\nn\\
&\le
c_{13}  \int_{V(1)}  \delta_D(z)\left(\frac{\delta_D(x)}{|x-z|}\wedge 1
\right)\left(\frac{\delta_D(z)}{|x-z|}\wedge 1
\right)\frac{\phi'(|x-z|^{-2})}
{|x-z|^{d+2}\phi(|x-z|^{-2})^2}  dz\nn\\
&\le 2 
c_{13}  \delta_D(x) \int_{V(1)} \left(\frac{\delta_D(z)}{|x-z|}\wedge 1
\right)\frac{\phi'(|x-z|^{-2})}
{|x-z|^{d+2}\phi(|x-z|^{-2})^2}  dz\nn\\
&\le  2 
c_{13}\delta_D(x) \int_{B(x, r/2)}\frac{\phi'(|x-z|^{-2})}
{|x-z|^{d+2}\phi(|x-z|^{-2})^2}  dz \le  
c_{14}   \delta_D(x)/ \phi(r^{-2}).
\end{align}
Combining \eqref{e:D2_422}, \eqref{e:D2_2} and  \eqref{e:D2_4}, 
we obtain 
\begin{align}
\label{e:D2_2n}
&\E_{x}\left[f \left(Y^D(\tau_{ V(1)})\right); \,
Y^D(\tau_{ V(1)})
\notin  D(2,2)\right]\nn \\&\asymp \frac{\delta_D(x)}{\phi(r^{-2})} 
\int_{D\setminus D(2,2)} f(y)
 \frac{1}{|y|_1} \left(\frac{\delta_D(y)}{|y|_1} \wedge 1 \right)   \frac{\mu(|y|^2)}
{|y|^{d-2}}    dy.
\end{align}

On the other hand, by Theorem \ref{carleson} and Lemma \ref{L:2}, we have
 \begin{align}
\label{e:D2_6}
&\E_x\left[f\left(Y^D(\tau_{ V(1)})\right);\,  Y^D(\tau_{
V(1)}) \in D(2,2)\right] \nn\\
& \le\, 
c_{15} \, f(x_0) \P_x\left(
Y^D(\tau_{V(1)}) \in  D(2,2)\right)\nn\\
&\le\, 
c_{16} \, f(x_0)
\frac{ \delta_D(x)   \phi'(r^{-2})  }
{r^{3}\phi(r^{-2})}.
\end{align}
Combining \eqref{e:BHP2},  \eqref{e:D2_2n},  \eqref{e:D2_6}, we get
 \begin{align}
&f(x) =  \E_x\left[f (Y^D(\tau_{ V(1)})); \,
Y^D(\tau_{ V(1)}) \in D(2,2)\right]+\,\E_x\left[ f
(Y^D(\tau_{ V(1)})); \,
Y^D(\tau_{ V(1)}) \notin  D(2,2)\right] \nn  \\
&\le  
c_{17} \delta_D(x)  \left( \frac{  \phi'(r^{-2})  }
{r^{3}\phi(r^{-2})}  f(x_0)   +\frac{1}{\phi(r^{-2})}  
\int_{D\setminus D(2,2)}   
\frac{f(w) }{|w|_1} \left(\frac{\delta_D(w)}{|w|_1} \wedge 1 \right)   \frac{\mu(|w|^2)}
{|w|^{d-2}}    dw   \right) \label{e:BHP1}
\end{align}
and 
 \begin{align}
&f(x) = \frac{1}{2}f(x) +\frac{1}{2}f(x)\nn\\
& \ge \frac{1}{2}\E_x\left[f\big(Y^D(\tau_{ V(1)})\big); Y^D_{
\tau_{ V(1)}} \in  
D(2,1)\setminus D(3/2,1)\right] \nonumber\\
&\quad+ \frac{1}{2} \E_x\left[ f
(Y^D(\tau_{ V(1)})); \,
Y^D(\tau_{ V(1)}) \notin  D(2,2)\right] \nn  \\
&\ge  
c_{18}  \delta_D(x)  \left( \frac{  \phi'(r^{-2})  }
{r^{3}\phi(r^{-2})}  f(x_0)   +\frac{1}{\phi(r^{-2})}
\int_{D\setminus D(2,2)}  
\frac{f(w) }{|w|_1} \left(\frac{\delta_D(w)}{|w|_1} \wedge 1 \right)   \frac{\mu(|w|^2)}
{|w|^{d-2}}    dw   \right). \label{e:BHP1n}
\end{align}
Therefore, 
 \begin{align}
f(x)  \asymp   \delta_D(x)  \left( \frac{  \phi'(r^{-2})  }
{r^{3}\phi(r^{-2})}  f(x_0)   +\frac{1}{\phi(r^{-2})}
\int_{D\setminus D(2,2)}   
 \frac{f(w) }{|w|_1} \left(\frac{\delta_D(w)}{|w|_1} \wedge 1 \right)   \frac{\mu(|w|^2)}
{|w|^{d-2}}    dw   \right). \label{e:BHP-yes}
\end{align}
For any $y\in D \cap B(Q,2^{-7} \kappa_0 r)$, we have the same estimate with $f(y_0)$ instead of $f(x_0)$ where $y_0=Q_y+\frac{r}{8}(y-Q_y)/|y-Q_y|$ and $Q_y\in \partial D$ with $|y-Q_y|=\delta_{D}(y)$. 
Since $D$ is $C^{1,1}$, using Theorem \ref{uhp}, 
$f(y_0)\asymp f(x_0)$. 
Therefore it follows 
from \eqref{e:BHP-yes} that for all $x,y \in D \cap B(Q,2^{-7} \kappa_0 r)$,
$$
\frac{f(x)}{f(y)}\,\le \,
c_{19}\,\frac{\delta_D(x)}{\delta_D(y)},
$$
which  proves the theorem. \qed

\section{Boundary Harnack principle in the interior of $D$}\label{s:bhi2}

In this section we assume that $d\ge 2$, 
$D$ is a domain in $\R^d$ and that \textbf{(A1)}-\textbf{(A2)} hold.

Recall that, for an open set  $U\subset D$, $X^U$ (respectively, $Y^{D,U}$)
is the process $X$ (respectively, $Y^D$)  killed upon exiting $U$.
One of the goals of this section is to show that when $U$ is relatively compact subset of $D$, the process $Y^{D,U}$ can be thought of as a non-local Feynman-Kac transform of $X^U$. 
Moreover, if $U$ is a certain $C^{1,1}$ domain,
the conditional gauge function related to this transform is bounded 
between two  positive constants
 which will imply that the Green functions of $X^U$ and $Y^{D,U}$ are comparable. We will prove a uniform version of this result in the sense that the comparability constants are independent of the set $U$ as long as its diameter is small and not larger than a multiple of its distance to the boundary.

Let $(\EE^{X^U}, \DD(\EE^{X^U}))$ be the Dirichlet form of $X^U$. Then, cf.~\cite[Section 13.4]{SSV},
\begin{eqnarray}
\EE^{X^U}(f,f)&=&\int_0^{\infty}\int_U f(x)(f(x)-P_s f(x))\, dx\,  \mu(s)ds\, , \label{e:df-of-XU}\\
\DD(\EE^{X^U})&=&\{f\in L^2(U, dx):\, \EE^{X^U}(f,f)<\infty\}\, .\label{e:dom-df-of-XU}
\end{eqnarray} 
Note that in computing $P_s f$, we extend the function $f$ to $\R^d$ by setting $f(x)=0$ for $x\in \R^d\setminus U$. 
Furthermore, for $f\in \DD(\EE^{X^U})$,
\begin{equation}\label{e:df-of-XU-alt}
\EE^{X^U}(f,f)=\frac12 \int_U \int_U (f(x)-f(y))^2 J^X(x,y)dy dx+\int_U f(x)^2 \kappa^X_U(x)dx,
\end{equation}
where $J^X$ is defined in \eqref{e:JX} and
$$
\kappa^X_U(x)=\int_{\R^d\setminus U}J^X(x,y)dy\, .
$$
The Dirichlet form $(\EE^{Y^D}, \DD(\EE^{Y^D}))$ of $Y^D$ is given by
$$
\EE^{Y^D}(f,f)=\int_0^{\infty}\int_D f(x)(f(x)-P_s^D f(x))\, dx\,  \mu(s) ds
$$
and $\DD(\EE^{Y^D})=\{f\in L^2(D,dx):\, \EE^{Y^D}(f,f)<\infty\}$. Moreover, for $f\in \DD(\EE^{Y^D})$, 
$$
\EE^{Y^D}(f,f)=\frac12 \int_D \int_D (f(x)-f(y))^2 J^{Y^D}(x,y)dydx+\int_D f(x)^2 \kappa^{Y^D}(x) dx,
$$
where $J^{Y^D}$ is defined in \eqref{e:JY} and
$$
\kappa^{Y^D}(x)=\int_0^{\infty}(1-P_t^D 1(x))\mu(t)\, dt\, .
$$
Hence, it follows that the Dirichlet form $(\EE^{Y^{D,U}}, \DD(\EE^{Y^{D,U}}))$ of $Y^{D,U}$ is equal to
\begin{eqnarray}
\EE^{Y^{D,U}}(f,f)&=&\int_0^{\infty}\int_U f(x)(f(x)-P^D_s f(x))\, dx\,  \mu(s)ds\, , \label{e:df-of-YDU}\\
\DD(\EE^{Y^{D,U}})&=&\{f\in L^2(U, dx):\, \EE^{Y^{D,U}}(f,f)<\infty\}\, .\label{e:dom-df-of-YDU}
\end{eqnarray}
Moreover, for $f\in \DD(\EE^{Y^{D,U}})$,
\begin{eqnarray}\label{e:df-of-YDU-alt}
\EE^{Y^{D,U}}(f,f)&=&\frac12 \int_U\int_U (f(x)-f(y))^2 J^{Y^D}(x,y)dy dx +\int_U f(x)^2 \left(\int_{D\setminus U}J^{Y^D}(x,y)dy\right)dx \nonumber \\
& &+\int_D f(x)^2 \kappa^{Y^D}(x) dx \nonumber\\
&=&\frac12\int_U\int_U (f(x)-f(y))^2 J^{Y^D}(x,y)dy dx  + \int_D f(x)^2 \kappa^{Y^D}_U(x) dx,
\end{eqnarray}
where
$$
 \kappa^{Y^D}_U(x)=\kappa^{Y^D}(x)+\int_{D\setminus U} J^{Y^D}(x,y)dy\, .
$$

We first need the following simple result.
\begin{lemma}\label{l:J-difference}
For $x,y\in D$,
\begin{equation}\label{e:J-difference}
J^X(x,y)-J^{Y^D}(x,y)\le j^X(\delta_D(y))\, .
\end{equation}
\end{lemma}
\pf
By \eqref{e:pD},  \eqref{e:JY} and \eqref{e:JX}, we have that 
\begin{eqnarray*}
J^X(x,y)-J^{Y^D}(x,y)&=&\int_0^{\infty} \E_x[p(t-\tau_D, W_{\tau_D}, y), \tau_D<t]\mu(t)\, dt\\
&=&\E_x\left[\int_{\tau_D}^{\infty}p(t-\tau_D, W_{\tau_D}, y)\mu(t)dt\right]\, .
\end{eqnarray*}
Since for every $s<t$ and $z\in \partial D$,
\begin{eqnarray*}
\int_s^{\infty}p(t-s,z,y)\mu(t)dt&=&\int_s^{\infty}(4\pi (t-s))^{-d/2}e^{-\frac{|z-y|^2}{4(t-s)}}\mu(t)dt\\
&\le &\int_s^{\infty}(4\pi (t-s))^{-d/2}e^{-\frac{\delta_D(y)^2}{4(t-s)}}\mu(t)dt\\
&=&\int_0^{\infty}(4\pi u)^{-d/2}e^{-\frac{\delta_D(y)^2}{4u}}\mu(u+s)du\\
&\le &\int_0^{\infty}(4\pi u)^{-d/2}e^{-\frac{\delta_D(y)^2}{4u}}\mu(u)du=j^X(\delta_D(y))\, ,
\end{eqnarray*} 
the claim of the lemma follows. \qed
\begin{lemma}\label{l:domains-equal}
Let $U$ be a relatively compact  open subset of $D$.
Then $\DD(\EE^{X^U})=\DD(\EE^{Y^{D,U}})$.
\end{lemma}
\pf Let $f\in L^2(U,dx)$. Using Fubini's theorem in the third line, Lemma 
\ref{l:J-difference} 
 and  $L^2(U,dx)\subset L^1(U,dx)$ in the last line, we get
\begin{eqnarray*}
\lefteqn{\int_0^{\infty}\int_U \big|f(x)(P_s^D f(x)-P_s f(x))\big| \, dx\, \mu(s)ds}\\
&\le & \int_0^{\infty}\int_U |f(x)| \int_U (p(s,x,y)-p^D(s,x,y))|f(y)|\, dy\, \mu(s)ds\\
&=& \int_U \int_U |f(x)| |f(y)| \left(\int_0^{\infty}(p(s,x,y)-p^D(s,x,y))\mu(s)ds\right) \, dx\, dy\\
&=& \int_U \int_U |f(x)| |f(y)|(J^X(x,y)-J^{Y^D}(x,y))\, dx\, dy\\
&\le & j^X(\mathrm{dist}(U,\partial D))\int_U \int_U |f(x)| |f(y)|\, dx\, dy <\infty\, .
\end{eqnarray*}
Together with \eqref{e:df-of-XU}--\eqref{e:dom-df-of-XU} and \eqref{e:df-of-YDU}--\eqref{e:dom-df-of-YDU},  this immediately implies the claim. \qed

We now give some other expressions for $\kappa^X_U$ and $\kappa^{Y^D}_U$. First, it is easy to see that
$$
\kappa^X_U(x)=\int_{(0,\infty)}(1-P_t 1_U(x))\mu(t)\, dt\, .
$$
Further,
\begin{eqnarray*}
\kappa^{Y^D}_U(x)&=&\int_{(0,\infty)}(1-P_t^D1(x))\mu(t)dt+\int_{D\setminus U}\int_{(0,\infty)}p^D(t,x,y)\mu(t)dt dy\\
&=&\int_{(0,\infty)}(1-P_t^D1(x))\mu(t)dt+\int_{(0,\infty)}\left(\int_{D\setminus U}p^D(t,x,y)dy\right)\mu(t)dt\\
&=&\int_{(0,\infty)}(1-P_t^D1(x))\mu(t)dt+\int_{(0,\infty)}P_t^D 1_{D\setminus U}(x)\mu(t)ft\\
&=&\int_{(0,\infty)}(1-P_t^D 1_U(x))\mu(t)dt\, .
\end{eqnarray*}
Let 
\begin{eqnarray*}
q_U(x)&:=&\kappa^{Y^D}_U(x)-\kappa^X_U(x)=\int_0^{\infty}(P_t 1_U-P_t^D 1_U)(x)\mu(t)\, dt\\
&=&\int_0^{\infty}\left(\P_x(W_t\in U)-\P_x(W_t^D\in U)\right)\mu(t)\, dt\, .
\end{eqnarray*}
Clearly, $q_U(x)\ge 0$ and
\begin{eqnarray*}
q_U(x)
&=&\int_0^{\infty}\int_U \E_x[p(t-\tau_D,W_{\tau_D}, y), t>\tau_D]dy\,  \mu(t)\, dt\\
&=&\int_U \E_x\int_{\tau_D}^{\infty}p(t-\tau_D, W_{\tau_D}, y)\mu(t)\, dt \, dy\\
&=&\int_U(J^X(x,y)-J^{Y^D}(x,y))dy\, .
\end{eqnarray*}

For $x,y\in D$, $x\ne y$, let 
$$
F(x,y):=\frac{J^{Y^D}(x,y)}{J^X(x,y)}-1=\frac{J^{Y^D}(x,y)-J^X(x,y)}{J^X(x,y)}\, ,
$$
and $F(x,x)=0$.  We also define $F(x,\partial)=0$, where $\partial$ denotes the cemetery point (this is where all killed processes end up). Then $-1<F(x,y)\le 0$.  Note that
$$
\int_U F(x,y) J^X(x,y)dy=-q_U(x)\, .
$$

The next lemma shows that  the absolute value of $F$ can be bounded by $1/2$ on balls of small radius sufficiently away from the boundary of $D$.
\begin{lemma}\label{l:estimate-of-F}
Assume that {\bf (A1)}--{\bf (A3)} hold.    If $d=2$, we assume that 
 {\bf (A4)} also holds.
There exists $b=b(\phi,d)>2$  such that for all $x_0\in D$ and all $r\in (0,1/b)$ satisfying $B(x_0, (b+1)r)\subset D$, we have that
$$
\sup_{x,y\in B(x_0,r)}|F(x,y)| \le \frac12\, .
$$
\end{lemma}
\pf 
First note that if $B(x_0, (b+1)r)\subset D$, then for every $y\in B(x_0,r)$ it holds that $\delta_D(y)>br$. Hence by \eqref{e:J-difference}, for $x, y\in B(x_0,r)$, $J^X(x,y)-J^{Y^D}(x,y)\le j^X(br)$. Therefore
$$
|F(x,y)|=\frac{J^X(x,y)-J^{Y^D}(x,y)}{J^X(x,y)}\le 
\frac{j^X(br)}{j^X(|x-y|)}\le \frac{j^X(br)}{j^X(2r)}\, .
$$
By \cite[Proposition 2.6]{KM14}, there exists $c_1=c_1(\phi,d)$ such that
\begin{equation}\label{e:estimate-of-j7}
c_1^{-1}\frac{\phi'(r^{-2})}{r^{d+2}}\le j^X(r) \le c_1\frac{\phi'(r^{-2})}{r^{d+2}}\, ,\quad 0<r<1\, .
\end{equation}
Hence for $b>2$ and $r\in (0,1)$ such that $br<1$, 
\begin{equation}\label{e:estimate-for-F}
|F(x,y)|\le c_1^2 \frac{\frac{\phi'(b^{-2}r^{-2})}{(br)^{d+2}}}{\frac{\phi'(2^{-2}r^{-2})}{(2r)^{d+2}}}\le c_1^2 \left(\frac{2}{b}\right)^{d+2}\frac{\phi'(b^{-2}r^{-2})}{\phi'(2^{-2}r^{-2})}\, .
\end{equation}
If $d\ge 3$, we use the 
fact that $\lambda\mapsto \lambda^2\phi'(\lambda)$ is increasing, 
cf.~Lemma \ref{l:pbf}(b), to conclude that for $b\ge 2$,
$$
|F(x,y)|\le c_1^2 \left(\frac{2}{b}\right)^{d-2} \frac{(b^{-2}r^{-2})^2\phi'(b^{-2}r^{-2})}{(2^{-2}r^{-2})^2\phi'(2^{-2}r^{-2})}\le c_1^2 \left(\frac{2}{b}\right)^{d-2} \, .
$$
Now choose $b>2$ such that $c:=c_1^2 (2/b)^{d-2} \le 1/2$. 

When $d=2$ we use \eqref{e:estimate-for-F} and  {\bf (A4)} to get
$$
|F(x,y)|\le c_1^2  \frac{1}{\sigma'}\left(\frac{2}{b}\right)^{4-\delta'}\, ,
$$
where $\sigma'$ and $\delta'$ are the constants from {\bf (A4)}. Again, we choose $b>2$ such that  the last expression is smaller than $1/2$.
\qed

Let $b>2$ be the constant from Lemma \ref{l:estimate-of-F}. For $r<1/b$, let $U\subset D$ be such that $\mathrm{diam}(U)\le r$ and $\mathrm{dist}(U,\partial D)\ge (b+2)r$. Then there exists a ball $B(x_0,r)$ such that $U\subset B(x_0,r)$ and $B(x_0, (b+1)r)\subset D$. By Lemma \ref{l:estimate-of-F} we see that 
\begin{align}
\label{e:Fbound}
|F(x,y)|\le 1/2 \qquad \text{ for all }x,y\in U.
\end{align}
 Hence we can define the non-local multiplicative functional
\begin{align*}
K_t^U
&:=\exp\left(\sum_{0<s\le t}\log(1+F(X_{s-}^U,X_s^U))-\int_0^t 
\int_U F(X_s^U,y)J^X(X_s^U,y)dy\, ds-\int_0^t q_U(X_s^U)ds\right)\\
&=\exp\left(\sum_{0<s\le t}\log(1+F(X_{s-}^U,X_s^U))\right)\, .
\end{align*}

Let 
$$
T^U_tf(x):=\E_x[K^U_t f(X_t^U)]\, .
$$
By \cite[(4.5) and Theorem 4.8]{CS03a}, $(T^U_t)_{t\ge0}$ is a strongly continuous semigroup on $L^2(U,dx)$ with the associated quadratic form $(\QQ, \DD(\EE^{X^U}))$ where $$
\QQ(f,f)=\EE^{X^U}(f,f)-\int_U\int_U f(x)f(y)F(x,y)J^X(x,y)\, dy\, dx\, .
$$
\begin{lemma}\label{l:eqofforms} 
Assume that {\bf (A1)}--{\bf (A3)} hold.    If $d=2$, we assume that 
 {\bf (A4)} also holds.
For $r<1/b$, let $U\subset D$ be such that $\mathrm{diam}(U)\le r$ and $\mathrm{dist}(U,\partial D)\ge (b+2)r$. Then
$$
(\QQ, \DD(\EE^{X^U}))=(\EE^{Y^{D,U}}, \DD(\EE^{Y^{D,U}}))\, .
$$
\end{lemma}

\pf
Note that
\begin{eqnarray*}
\QQ(f,f)&=& \frac12 \int_U\int_U (f(x)-f(y))^2(J^{Y^D}(x,y)-F(x,y)J^X(x,y))\, dy\, dx 
+\int_U f(x)^2\kappa^X_U(x)\, dx \\
& & - \frac12  \int_U\int_U f(x)f(y)F(x,y)J^X(x,y)\, dy\, dx\\
&=& \frac12 \int_U\int_U (f(x)-f(y))^2 J^{Y^D}(x,y)\, dy\, dx +\int_U f(x)^2\int_U F(x,y)J^X(x,y)\, dy\, dx \\
& &+\int_U f(x)^2\kappa^X_U(x)\, dx\\
&=& \frac12 \int_U\int_U (f(x)-f(y))^2 J^{Y^D}(x,y)\, dy\, dx +\int_U f(x)^2  q_U(x) \, dx+\int_U f(x)^2\kappa^X_U(x)\, dx\\
& =&\EE^{Y^{D,U}}(f,f)\, .
\end{eqnarray*}
It was shown in Lemma \ref{l:domains-equal} that  $\DD(\EE^{X^U})=\DD(\EE^{Y^{D,U}})$, which finishes the proof.
\qed

For $x,y\in U$, $x\neq y$, let
$$
u^U(x,y):=\E_x^y [K^U_{\tau_U^X}]
$$
be the conditional gauge function for $K^U_t$. Here $\E_x^y$ denotes the expectation with respect to the conditional probability $\P_x^y$ defined via Doob's $h$-transform with 
respect to the Green function $G^X_U(\cdot, y)$ starting from $x\in D$.
Since $F\le 0$, we have $\log(1+F)\le 0$, hence $K^U_{\tau_U^X}\le 1$. Therefore, $u^U(x,y)\le 1$.
Define
$$
V^U(x,y)=u^U(x,y)G_U^X(x,y), \quad x,y\in U.
$$ 
It follows from \cite[Lemma 3.9]{che02} that $V^U(x,y)$ is the Green function for the semigroup $(T^U_t)_{t\ge 0}$. 
Combining this with Lemma \ref{l:eqofforms}
we can conclude that $V^U$ is equal to the Green function $G^{Y^D}_U$ of $Y^{D,U}$. Therefore,
\begin{equation}\label{e:gf-XU-YDU-comparable}
G^{Y^D}_U(x,y)=u^U(x,y)G_U^X(x,y), \quad x,y\in U.
\end{equation}
Our next goal is to show that for $C^{1,1}$ open sets $U$, the conditional gauge function $u^U$ is bounded below by a 
strictly positive constant uniform in the diameter of $U$.
Together with the above equality this proves that the Green function 
of $Y^{D,U}$ is comparable to that of $X^U$.

In the remainder of the section, we assume that 
\textbf{(A1)}--\textbf{(A3)} hold. If the constant $\delta$ in {\bf (A3)} satisfies $0<\delta \le 1/2$,  we assume that {\bf (A7)} also holds.
 If $d=2$, we assume that 
 {\bf (A4)} and  {\bf (A5)} also hold. Without loss of generality we assume $\phi(1)=1$. 

For $r>0$ define
$$
\phi^r(\lambda):=\frac{\phi(\lambda r^{-2})}{\phi(r^{-2})}\, ,\quad \lambda >0\, .
$$
Then $\phi^{r}$ is a Bernstein function with potential density $u^r$ and  L\'evy density $\mu^r$ given by
$$
u^r(t)=r^2\phi(r^{-2})u(r^2 t)\, ,
\qquad \mu^r(t)=r^2\phi(r^{-2})^{-1}\mu(r^2 t)\, ,\quad t>0\, ,
$$
cf.~\cite[(2.7), (2.8)]{KSV14}. In particular, $\phi^r$ satisfies {\bf (A1)}. Further, if $0<r\le 1$, then for $\lambda\ge 1$ we have that $\lambda r^{-2}\ge 1$, and thus for $t\ge 1$,
$$
\frac{(\phi^r)'(\lambda t)}{(\phi^r)'(\lambda)}=\frac{\phi'(\lambda t r^{-2})}{\phi'(\lambda r^{-2})}\le \sigma t^{-\delta}\, .
$$
Hence, {\bf (A3)} holds for $\phi^r$ with the same constants. Similarly, one can see that $\phi^r$ satisfies
  {\bf (A4)} (if $d=2$) and {\bf (A7)} (if $0<\delta \le 1/2$)
  with the same constants. It is immediate that $\phi^r$ also satisfies {\bf (A5)}. It remains to check {\bf (A2)}. Suppose first that $t>1$ is such that $r^2 t\le 1$ 
    (so that $r^2 (t+1)\le 2$). 
  Then by using \eqref{e:lower-estimate-mu} and \eqref{e:upper-estimate-mu} in the first inequality and the fact that $\phi'$ is decreasing in the last, we get
$$
\frac{\mu^r(t+1)}{\mu^r(t)}=\frac{\mu(r^2(t+1))}{\mu(r^2 t)}\ge c_1\frac{\phi'(r^{-2}(t+1)^{-1})}{r^4 (t+1)^2}\, \frac{r^4 t^2}{\phi'(r^{-2}t^{-1})} \\
\ge \frac{c_1}{2}\, .
$$
If $t>1$ and $r^2 t >1$, then by {\bf (A2)}, we get
$$
\frac{\mu^t(t+1)}{\mu^r(t)}=\frac{\mu(r^2(t+1))}{\mu(r^2 t)}\ge \frac{\mu(r^2t+1)}{\mu(r^2 t)}\ge 
c_2.
$$
Thus $\phi^r$ also satisfies {\bf (A2)} with the constant independent of $a$.

Let $S^r=(S^r_t)_{t\ge 0}$ be a subordinator with Laplace exponent $\phi^r$ independent of the Brownian motion  $W$. Let $X^r=(X^r_t)_{t\ge 0}$ be defined by $X^r_t=W_{S^r_t}$. Then $X^r$ is an isotropic L\'evy process with characteristic exponent $\phi^r(|\xi|^2)=\phi(|\xi|^2 r^{-2}) / \phi(r^{-2})$, $\xi\in \R^d$, which shows that $X^r$ is identical in law to the process 
$\{r^{-1}X_{t/ \phi(r^{-2})}\}_{t\ge 0}$.

Let $V\subset \R^d$ be a bounded $C^{1,1}$ open set such that $0\in \partial V$.  For $r\in (0,1]$, let $V^r=\{rx:\, x\in V\}$ be the scaled version of $V$. Denote by $G_V^{X^r}$ (respectively $G_{V^r}^X$) the Green function of $V$ with respect to the process $X^r$ (respectively the Green function of $V^r$ with respect to the process $X$). Then by scaling,
\begin{equation}\label{e:scaling-for-G}
G_{V^r}^X(x,y)=r^{-d}\phi(r^{-2})^{-1}G_V^{X^r}(x/r,y/r)\, , \quad x,y\in V^r\, .
\end{equation}
For any open set $U\subset \R^d$, we let
\begin{align}
\label{e:gU}
g^r_U(x,y):=\left(1\wedge \frac{\phi^r(|x-y|^{-2})}{\sqrt{\phi^r(\delta_U(x)^{-2})\phi^r(\delta_U(y)^{-2})}}\right)\, \frac{(\phi^r)'(|x-y|^{-2})}{|x-y|^{d+2}\phi^r(|x-y|^{-2})^2}\, , \quad x,y\in U\,,
\end{align}
and $g_U(x,y):=g^1_U(x,y)$.
\begin{prop}\label{p:scale-inv-gfe}
Let $V\subset \R^d$ be a bounded $C^{1,1}$ open set with characteristics $(R,\Lambda)$ such that $0\in \partial V$ and $\mathrm{diam}(V)\le 1$. There exists a constant $C=C(R,\Lambda,\phi,d)\ge 1$ such that for every $r\in (0,1]$,
$$
C^{-1}g_{V^r}(x,y)\le G_{V^r}^X(x,y)\le C g_{V^r}(x,y)\, , \qquad x,y\in V^r\, .
$$
The dependence of $c$ on $\phi$ is only through the constants in assumptions {\bf (A1)}--{\bf (A5)} and {\bf (A7)}.
\end{prop}
\pf Let $r\in (0,1]$. By \cite[Theorem 1.2]{KM14} applied to the process $X^r$, there exists a constant $c=c(\mathrm{diam}(V), R, \Lambda, \phi^r, d)$ such that 
\begin{equation}\label{e:scale-inv-gfe1}
c^{-1}g^r_V(x,y)\le G_V^{X^r}(x,y) \le c g^r_V(x,y)\, ,\quad x,y\in V\, .
\end{equation}
Moreover, although not explicitly mentioned, the dependence of $c$ on $\phi^r$ is only through the constants in 
{\bf (A1)}-{\bf (A5)} and {\bf (A7)} 
for $\phi^r$. But as shown above, those constants are independent of $r\in (0,1]$. Since 
$\mathrm{diam}(V)\le 1$, 
we conclude that $c=c(R, \Lambda, \phi, d)$ where dependence of $c$ on $\phi$ is only through the constants in 
{\bf (A1)}-{\bf (A5)} and {\bf (A7)} for $\phi$.

By a straightforward computation we get that
\begin{equation}\label{e:scale-inv-gfe2}
g^r_V(x/r,y/r)=r^d \phi(r^{-2})g_{V^r}(x,y)\, , \quad x,y\in V^r\, .
\end{equation}
The claim of the proposition now follows by combining \eqref{e:scaling-for-G}, \eqref{e:scale-inv-gfe1} and \eqref{e:scale-inv-gfe2}. \qed

We note that the assumption in Proposition \ref{p:scale-inv-gfe}
that $0\in \partial V$ is irrelevant. For any $z\in \partial V$  we could use the scaling $V^r=\{r(x-z)+z;\, x\in V\}$ and obtain the same result.
\begin{lemma}\label{l:3G}
Let $R>0$ and $\Lambda >0$.
There exists $C=C(R,\Lambda,\phi,d)>0$ such that for every $r\in (0,1]$ and every $C^{1,1}$ open set $U\subset D$ with characteristics $(rR,\Lambda/r)$ and $\mathrm{diam}(U)\le r$ satisfying $\mathrm{dist}(U,\partial D)\ge r$, 
\begin{equation}\label{e:3G}
\int_U \int _U \frac{G_U^X(x,z)G_U^X(w,y)}{G_U^X(x,y)}|F(z,w)|J^X(z,w)dz\, dw\le C\, , \quad \textrm{for all }x,y\in U\, .
\end{equation}
\end{lemma}
\pf Since $\mathrm{dist}(U,\partial D)\ge r$, we have 
by Lemma \ref{l:J-difference} 
 that 
 \begin{align}
 \label{e:FJj}
 |F(z,w)|J^X(z,w)\le j^X(r) \qquad \text{for all }z,w\in U.
 \end{align}
Let $V^r=U$ and define $V=r^{-1}V^r=\{r^{-1}x;\ x\in V^r\}$. Then $V$ is a $C^{1,1}$ open set with characteristics $(R,\Lambda)$ and $\mathrm{diam}(V)\le 1$. By Proposition \ref{p:scale-inv-gfe} there exists $c_1=c_1(R,\Lambda,\phi,d)>0$ such that 
\begin{align}
\label{e:gUGU}
c_1^{-1}g_{U}(x,y)\le G_{U}^X(x,y)\le c_1 g_{U}(x,y)\, , \qquad x,y\in U\, .
\end{align}

Next note that
$$
\int_U \int_U \frac{G_U^X(x,z)G_U^X(w,y)}{G_U^X(x,y)} dz\, dw =\frac{\E_x(\tau_U^X)\E_y(\tau_U^X)}{G_U^X(x,y)}\, .
$$
By using \eqref{e:gUGU} and \cite[Lemma 7.1]{CKS}
we get that
\begin{eqnarray*}
\E_x(\tau_U^X)&=&\int_U G_U^X(x,y)\, dy \le 2 c_1  \phi(\delta_U(x)^{-2})^{-1/2} \int_U\frac{\phi'(|x-y|^{-2})}{|x-y|^{d+2}\phi(|x-y|^{-2})^{3/2}}\, dy\\
&\le& 2c_1 c(d)\phi(\delta_U(x)^{-2})^{-1/2}\int_0^{2r}t^{d-1}\frac{\phi'(t^{-2})}{t^{d+2}\phi(t^{-2})^{3/2}}\ dt\\
&=&2c_1 c(d)\phi(\delta_U(x)^{-2})^{-1/2}\int_0^{2r}d\left(\phi(t^{-2})^{-1/2}\right)\\
&=& 2c_1 c(d)\phi(\delta_U(x)^{-2})^{-1/2} \phi((2r)^{-2})^{-1/2}\\
&\le &c_2 \phi(\delta_U(x)^{-2})^{-1/2}\phi(r^{-2})^{-1/2}\, ,
\end{eqnarray*}
where $c_2=c_2(R,\Lambda,\phi,d)>0$.

If $\phi(|x-y|^{-2})\le \phi(\delta_U(x)^{-2})^{1/2}\phi(\delta_U(y)^{-2})^{1/2}$, then
\begin{eqnarray*}
\frac{\E_x(\tau_U^X)\E_y(\tau_U^X)}{G_U^X(x,y)}&\le &c_2^2 c_1 \frac{\phi(r^{-2})^{-1}\phi(\delta_U(x)^{-2})^{-1/2}\phi(\delta_U(y)^{-2})^{-1/2}}{\frac{\phi'(|x-y|^{-2})}{\phi(\delta_U(x)^{-2})^{1/2}\phi(\delta_U(y)^{-2})^{1/2}|x-y|^{d+2}\phi(|x-y|^{-2})}}\\
&=&c_2^2 c_1	\phi(r^{-2})^{-1}|x-y|^{d-2}\frac{\phi(|x-y|^{-2})}{(|x-y|^{-2})^2\phi'(|x-y|^{-2})}\\
&\le &c_2^2 c_1	\phi(r^{-2})^{-1} (2r)^{d-2} \frac{\phi((2r)^{-2})}{(2r)^{-4}\phi'((2r)^{-2})}\\
&\le &2^{d+2} c_2^2 c_1 r^{d+2}\phi'(r^{-2})^{-1}\, .
\end{eqnarray*}
In the third line we have used Lemma \ref{l:pbf}(b) (with $a=1$).

Suppose now that  $\phi(|x-y|^{-2})> \phi(\delta_U(x)^{-2})^{1/2}\phi(\delta_U(y)^{-2})^{1/2}$.  By using the obvious inequality $\phi(\delta_U(x)^{-2})^{-1/2}\le \phi(r^{-2})^{-1/2}$ and Lemma \ref{l:pbf}(b) again (this time with $a=2$), we get
\begin{eqnarray}
\frac{\E_x(\tau^X_U)\E_y(\tau^X_U)}
{G_U^X(x,y)}&\le &c_2^2 c_1  \phi(r^{-2})^{-2}|x-y|^{d-2}\frac{|x-y|^{4} \phi(|x-y|^{-2})^2}{\phi'(|x-y|^{-2})}\nn\\
&\le & 2^{d+2}c_2^2 c_1  r^{d+2}\phi'(r^{-2})^{-1}\, . \label{e:EEG}
\end{eqnarray}
Using  the  inequalities \eqref{e:estimate-of-j7} and \eqref{e:FJj},  and applying \eqref{e:EEG},  we conclude that
$$
\int_U \int _U \frac{G_U^X(x,z)G_U^X(w,y)}{G_U^X(x,y)}|F(z,w)|J^X(z,w)dz\, dw \le 2^{d+2} c_2^2 c_1  
r^{d+2}\phi'(r^{-2})^{-1}j^X(r) \le C
$$
with a constant $C=C(R,\Lambda,\phi,d)$.  \qed
\begin{lemma}\label{l:3G-application}
Let $R>0$ and $\Lambda >0$. Then for every $r\in (0,1/b)$ and every $C^{1,1}$ open set $U\subset D$ with characteristics $(rR,\Lambda/r)$ and $\mathrm{diam}(U)\le r$ satisfying $\mathrm{dist}(U,\partial D)\ge (b+2)r$, we have
$$
\E_x^y\sum_{0<s\le \tau_U^X}|F(X_{s-}^U, X_s^U)|\le C\, ,
$$
where $C$ is the constant from Lemma \ref{l:3G} and $b>2$ the constant from Lemma \ref{l:estimate-of-F}.
\end{lemma}
\pf By \cite[Proposition 3.3]{CS03a} we have
\begin{eqnarray*}
\E_x^y\sum_{0<s\le 
\tau^X_U}|F(X_{s-}^U, X_s^U)|&=&\E_x\int_0^{\tau_U^X}\int_U \frac{|F(X_s,w)|G_U^X(w,y)}{G_U^X(x,y)}J^X(X_s,w)\, dw\, ds\\
&=&\int_U \int _U \frac{G_U^X(x,z)G_U^X(w,y)}{G_U^X(x,y)}|F(z,w)|J^X(z,w)dz\, dw \, .
\end{eqnarray*}
The claim is now a consequence of Lemma \ref{l:3G}. \qed
\begin{lemma}\label{l:bounded-gauge}
Let $R>0$ and $\Lambda >0$. There exists $C=C(R,\Lambda, \phi,d)\in (0,1)$ such that for every $r\in (0,1/b)$ and every $C^{1,1}$ open set $U\subset D$ with characteristics $(rR,\Lambda/r)$ and $\mathrm{diam}(U)\le r$ satisfying $\mathrm{dist}(U,\partial D)\ge (b+2)r$, we have
$$
C\le u^U(x,y) \le 1\, , \quad x,y\in U, x\neq y\, .
$$
\end{lemma}
\pf The upper bound is clear. By \eqref{e:Fbound}, $F\ge -1/2$ on $U \times U$. 
By using Jensen's inequality in the first inequality, the fact that $\log(1+t)\ge t$ for $t\in [-1/2,0)$ in the second, and Lemma \ref{l:3G-application} in the third, we have
\begin{eqnarray*}
u^U(x,y)&\ge & \exp\left\{\E_x^y \left(\sum_{0<s\le \tau_U^X} \log\left(1+F(X_{s-}^U, X_s^U)\right)\right)\right\}\\
&\ge &\exp\left\{ \E_x^y \left(\sum_{0<s\le 
\tau_U^X} F(X_{s-}^U, X_s^U)\right)\right\}\,\ge \, \exp\{-C\}\, .
\end{eqnarray*}
\qed

Combining this lemma with \eqref{e:gf-XU-YDU-comparable} we arrive at
\begin{prop}\label{p:green-comparable}
Let $R>0$ and $\Lambda >0$. There exists $C=C(R,\Lambda, \phi,d)\in (0,1)$ such that for every $r\in (0,1/b)$ and every $C^{1,1}$ open set $U\subset D$ with characteristics $(rR,\Lambda/r)$ and $\mathrm{diam}(U)\le r$ satisfying $\mathrm{dist}(U,\partial D)\ge (b+2)r$, we have 
$$
CG_U^X(x,y)\le G^{Y^D}_U(x,y)\le G_U^X(x,y)\, ,\quad x,y\in U\, .
$$
\end{prop}

A similar method was used to establish comparability of the Green functions of some one-dimensional processes in \cite{Wag}.

\noindent
\textbf{Proof of Theorem \ref{t:bhp-2}:} 
Let $R=\delta_D(Q)\wedge 1$. 
Choose
a $C^{1,1}$-function $\varphi:\R^{d-1}\mapsto \R$ satisfying $\varphi(\wt 0)=0$, $\nabla \varphi(\wt 0)=(0,\dots, 0)$, $\| \nabla \varphi \|_{\infty}\le \Lambda$, $|\nabla \varphi(\wt y)-\nabla \varphi(\wt w)|\le \Lambda |\wt y-\wt w|$, and an orthonormal coordinate system $CS_Q$ with its origin at $Q$ such that
$$
B(Q,R)\cap E=\{y=(\wt{y},y_d)\in B(0,R) \textrm{ in } CS_Q:\, y_d>\varphi(\wt{y})\}\, .
$$
Define $\rho_Q(x)=x_d-\varphi(\wt{x})$, where $(\wt{x},x_d)$ are the coordinates in $CS_Q$. 
Recall that $\kappa=(1+(1+\Lambda)^2)^{-1/2}$. For $r>0$, define the box
$$
E_Q(r)=\{y\in E:\, 0<\rho_Q(y)<\kappa r/2, |\wt{y}|<\kappa r/2\},
$$
so that $\mathrm{diam}(E_Q(r))\le r$.
There exists $L=L(R,\Lambda, d)$ such that for every $r\le \kappa R$ one can find a $C^{1,1}$ domain $V_Q(r)$ with characteristics $(rR/L, \Lambda L/r)$ such that 
$E_Q(r/2)\subset V_Q(r)\subset E_Q(r)$ (see \cite[Lemma 2.2]{So}).
 Furthermore, $B(Q,\kappa^2 r/4)\cap E \subset E_Q(r/2)$ and 
 $E_Q(r)\subset B(Q, r)\cap E$ (see \eqref{e:Uzr}).
Since $r<(\delta_D(Q) \wedge 1)/(b+2)$,
we conclude that $\mathrm{dist}(V_Q(r),\partial D)\ge (b+2)r$.
Now it follows from Proposition \ref{p:green-comparable} (with $(R/L, \Lambda L)$ instead of $(R,\Lambda)$) that 
\begin{align}
\label{e:GV}
 c_0 \, G_{V_Q(r)}^X(x,y)\le G^{Y^D}_{V_Q(r)}(x,y)\le G_{V_Q(r)}^X(x,y)\, ,\quad x,y\in V_Q(r)\, ,
\end{align}
where $c_0=c_0(\delta_D(Q) \wedge 1,\Lambda, \phi,d)$ is independent of $r$.

Let $B=B(Q,2 r)$. Then $V_Q(r)\subset B\subset D$ and $\mathrm{dist}(B,\partial D)\ge br/2$. Recall that by Lemma \ref{l:ilb4jk}
there exists a constant $c_1=c_1(b)\in (0,1)$ such that for every 
$r<(\delta_D(Q) \wedge 1)/(b+2)$,  
\begin{align}
\label{e:JXYD}
c_1 J^X(w,y)\le J^{Y^D}(w,y)\le J^X(w,y)\, ,\qquad w,y\in B\, .
\end{align}

Let $f$ be regular harmonic in $E\cap B(Q, r)$ with respect to $Y^D$ and vanish on $E^c\cap B(Q, r)$. Then for $x\in V_Q(r)$,
\begin{equation}\label{e:f-decomposition}
f(x)=\E_x[f(Y^D_{\tau_{V_Q(r)}})
; Y^D_{\tau_{V_Q(r)}}\in B]+\E_x[f(Y^D_{\tau_{V_Q(r)})}
; Y^D_{\tau_{V_Q(r)}}\in D\setminus B]=:f_1(x)+f_2(x)\, .
\end{equation}
We first estimate $f_1(x)$. 
By using \eqref{e:GV} and \eqref{e:JXYD} we get
\begin{eqnarray*}
f_1(x)&=&\int_{B\setminus V_Q(r)}\int_{V_Q(r)} G^{Y^D}_{V_Q(r)}(x,w)J^{Y^D}(w,y) 
f(y)\, dy\, dw\\
&\asymp &\int_{B\setminus V_Q(r)}\int_{V_Q(r)} 
G_{V_Q(r)}^X(x,w)J^X(w,y) 
f(y)\, dy\, dw\\
&=&\E_x[f(X_{\tau_{V_Q(r)}}); X_{\tau_{V_Q(r)}}\in B]=:h(x)\, ,
\end{eqnarray*}
where the comparison constants in the second line depend only on 
$\delta_D(Q) \wedge 1,\Lambda, \phi,d$.
Since the function $h$ is regular harmonic in $B(Q,\kappa^2 r/4) \cap E$ with respect to $X$ and vanishes on 
$B(Q,\kappa^2 r/4) \cap E^c$,
 we can use the factorization from \cite[Lemma 5.4]{KM14} to conclude that
$$
h(x)\asymp \E_x[\tau_{B(Q,\kappa^2 r/4) \cap E}^X] \int_{B(Q, \kappa^2 r/8)^c}J^X(Q,y)h(y)\, dy\, ,\quad x\in B(Q,\kappa^2 r/8)
$$
(with the comparison constants depending on $\phi, d$). Hence, 
\begin{equation}\label{e:f1}
 f_1(x) \asymp  \E_x[\tau_{B(Q,\kappa^2 r/4) \cap E}^X] \int_{B(Q, \kappa^2 r/8)^c}J^X(Q,y)h(y)\, dy\, ,\quad x\in B(Q,\kappa^2 r/8)\, .
\end{equation}

In order to estimate $f_2(x)$ we use Proposition \ref{p:estimate-of-J-away} in the second line below to conclude that
\begin{eqnarray}\label{e:f2}
f_2(x)&=&\E_x\int_0^{\tau_{V_Q(r)}^{Y^D}}\int_{D\setminus B} f(y)J^{Y^D}(Y^D_s,y)\, dy\, ds \nonumber\\
&\asymp &\E_x\int_0^{\tau_{V_Q(r)}^{Y^D}}\int_{D\setminus B} f(y)J^{Y^D}(Q,y)\, dy\, ds \nonumber\\
&=& \E_x[\tau_{V_Q(r)}^{Y^D}]\int_{D\setminus B} f(y)J^{Y^D}(Q,y)\, dy
\end{eqnarray}
(where the comparison constants depend on the constants from {\bf (B1)} and {\bf (B2)}).  
Now, by Proposition \ref{p:scale-inv-gfe} and \eqref{e:GV}  we have that 
for $x\in B(Q, \kappa^4 r/(32))$,
\begin{align}
\label{e:compX}
  \E_x[\tau_{V_Q(r)}^{Y^D}] 
 \asymp \phi(r^{-2})^{-1/2}\phi(\delta_{V_Q(r)}^{-2}(x))^{-1/2}=\phi(r^{-2})^{-1/2}\phi(\delta_{E}^{-2}(x))^{-1/2}\, ,
\end{align}
\begin{align}
\label{e:compY1}
 &\E_x[\tau_{B(Q,\kappa^2 r/4) \cap E}^X]   \le  \E_x[\tau_{V_Q(r)}^X] \nn\\
 & \le c_2 \phi(r^{-2})^{-1/2}\phi(\delta_{V_Q(r)}^{-2}(x))^{-1/2}=c_2\phi(r^{-2})^{-1/2}\phi(\delta_{E}^{-2}(x))^{-1/2}\, ,
\end{align}
and
\begin{align}
\label{e:compY2}
& \E_x[\tau_{B(Q,\kappa^2 r/4) \cap E}^X]   \ge  \E_x[\tau_{2^{-2}\kappa^2V_Q(r)}^X] 
\nn\\
&\ge c_3\phi(r^{-2})^{-1/2}\phi(\delta_{2^{-2}\kappa^2V_Q(r)}^{-2}(x))^{-1/2}=c_3\phi(r^{-2})^{-1/2}\phi(\delta_{E}^{-2}(x))^{-1/2}\, .
\end{align}
By combining \eqref{e:f-decomposition}-\eqref{e:compY2}, we get that for $x\in B(Q, \kappa^4 r/(32))$, 
$$
f(x)\asymp \phi(\delta_{E}^{-2}(x))^{-1/2}\phi(r^{-2})^{-1/2}\left(\int_{B(\kappa^2 r/8)^c}J^X(Q,y)h(y)\, dy+\int_{D\setminus B} f(y)J^{Y^D}(Q,y)\, dy\right)\, .
$$
This approximate factorization of the regular harmonic function $f$ immediately implies the claim. \qed

\bigskip
\noindent
{\bf Acknowledgements:} 
We thank the referee for helpful comments on the first version of this paper.
\bigskip
\noindent

 \vspace{.1in}
\begin{singlespace}


\end{singlespace}

\end{doublespace}
\vskip 0.1truein

\parindent=0em

{\bf Panki Kim}

Department of Mathematical Sciences and Research Institute of Mathematics,

Seoul National University, Building 27, 1 Gwanak-ro, Gwanak-gu Seoul 151-747, Republic of Korea

E-mail: \texttt{pkim@snu.ac.kr}

\bigskip

{\bf Renming Song}

Department of Mathematics, University of Illinois, Urbana, IL 61801,
USA, and\\
School of Mathematical Sciences, Nankai University, Tianjin 300071, PR China

E-mail: \texttt{rsong@math.uiuc.edu}

\bigskip

{\bf Zoran Vondra\v{c}ek}

Department of Mathematics, Faculty of Science, University of Zagreb, Zagreb, Croatia, and \\
Department of Mathematics, University of Illinois, Urbana, IL 61801,
USA

Email: \texttt{vondra@math.hr}
\end{document}